\input amstex
\documentstyle{amsppt}

\define\r{{\bold R}}

\define\codim{\text{codim~}}
\define\M{{\Cal M}}
\define\R{{\Cal R}}
\define\E{{\Cal E}}
\define\La{{\Cal L}}
\define\LT{{\Cal LT}}
\define\Aa{{\goth A}}
\define\bb{{\goth B}}
\define\cc{{\goth C}}

\define\nem{\overline {NE}}
\define\Tt{{\goth T}}
\define\p{\bold P}

\NoBlackBoxes

\magnification1200

\define\verte{\hbox to15pt{\hfill \vbox to30pt{
\vfill\nointerlineskip
\hbox to15pt{\hfill $\bigcirc$ \hfill}
\nointerlineskip
\vfill}\hfill}}

\define\vertex#1{\hbox to15pt{\hfill \vbox to30pt{
\hbox to15pt{\hfill $#1$ \hfill}
\nointerlineskip\vfill
\hbox to15pt{\hfill $\bigcirc$ \hfill}
\nointerlineskip
\vskip10pt}\hfill}}

\define\vertexb#1{\hbox to20pt{\hfill \vbox to30pt{
\hbox to20pt{\hfill $#1$ \hfill}
\nointerlineskip\vfill
\hbox to20pt{\hfill $\bullet$ \hfill}
\nointerlineskip
\vskip13pt}\hfill}}

\define\darr{\hbox to30pt{\hfill
\vbox to30pt{\vfill \nointerlineskip
\hbox to25pt{\rightarrowfill} \nointerlineskip
\hbox to25pt{\leftarrowfill} \nointerlineskip
\vfill}
\hfill}}

\define\darrup#1{
\hbox to30pt{\hfill
\vbox to30pt{\vfill
\hbox to25pt{\hfill $#1$ \hfill}
\nointerlineskip
\hbox to25pt{\rightarrowfill} \nointerlineskip
\hbox to25pt{\leftarrowfill} \nointerlineskip
\vskip8pt}
\hfill}}

\define\darrdown#1{
\hbox to30pt{\hfill
\vbox to30pt{\vskip8pt \nointerlineskip
\hbox to25pt{\rightarrowfill} \nointerlineskip
\hbox to25pt{\leftarrowfill} \nointerlineskip
\hbox to25pt{\hfill $#1$ \hfill}
\vfill}
\hfill}}

\define\darrupdown#1#2{
\hbox to30pt{\hfill
\vbox to30pt{
\hbox to25pt{\hfill $#1$ \hfill} \nointerlineskip
\vskip2pt\vfill
\hbox to25pt{\rightarrowfill} \nointerlineskip
\hbox to25pt{\leftarrowfill} \nointerlineskip\vfill
\hbox to25pt{\hfill $#2$ \hfill}
\vfill}
\hfill}}

\define\horlines{\hbox to60pt{\hfil \vrule height30pt \hskip-3pt
\hbox to55pt{\hfil  \vbox to30pt{
\hrule width50pt
\nointerlineskip\vfil
\hrule width50pt
\nointerlineskip\vfil
\hrule width50pt
\nointerlineskip\vfil
\hrule width50pt
\nointerlineskip\vfil
\hrule width50pt
\nointerlineskip\vfil
\hrule width50pt
\nointerlineskip\vfil
\hrule width50pt
\nointerlineskip\vfil
\hrule width50pt
\nointerlineskip\vfil
\hrule width50pt
\nointerlineskip\vfil
\hrule width50pt
\nointerlineskip\vfil
\hrule width50pt}\hfil}
\hskip-6pt \vrule height30pt \hfil}}

\define\verlines{\hbox to60pt{\hfil
\vbox to30pt{
\hrule width50pt
\hbox to 50pt{\vrule height30pt
\hfil \vrule height30pt \hfil \vrule height30pt \hfil \vrule height30pt
\hfil \vrule height30pt \hfil \vrule height30pt \hfil \vrule height30pt
\hfil \vrule height30pt \hfil \vrule height30pt \hfil \vrule height30pt}
\hrule width50pt
}
\hfil}}

\topmatter

\title
Diagram Method for $3$-folds and its application to K\"ahler Cone
and Picard Number of Calabi-Yau $3$-folds. I
\endtitle

\author
Viacheslav V. Nikulin \footnote{Partially supported by
Grant of Russian Fund of Fundamental Research; Grant of AMS;
Grant of ISF M16000.}
\endauthor

\address
Steklov Mathematical Institute,
ul. Vavilova 42, Moscow 117966, GSP-1, Russia.
\endaddress

\email
slava\@nikulin.mian.su
\endemail

\abstract
We prove the general diagram method theorem valid for the quite large
class of $3$-folds with $\Bbb Q$-factorial singularities
(see Basic Theorems  1.3.2 and 3.2 and also Theorem 2.2.6).
This gives the generalization of our results about Fano $3$-folds
with $\Bbb Q$-factorial terminal singularities
(Preprint alg-geom/9311007).

As an application, we get the following result about $3$-dimensional
Calabi-Yau manifolds $X$: Assume that
the Picard number $\rho (X) > 40$. Then one of two cases
(i) or (ii) holds:
(i) There exists a small extremal ray on $X$.
(ii) There exists a $nef$ element $h$ such that $h^3=0$
(thus, the $nef$ cone $NEF(X)$ and the cubic intersection hypersurface
${\Cal W}_X$ have a common point; here, we don't claim that $h$
is rational!).

As a corollary, we get: Let $X$ be a $3$-dimensional Calabi-Yau
manifold. Assume that the $nef$ cone $NEF(X)$ is finite polyhedral
and $X$ does not have a small extremal ray. Then there exists
a rational $nef$ element $h$ with $h^3=0$ if $\rho(X)>40$.

To prove these results on Calabi-Yau manifolds, we also use
one result by V.V. Shokurov on the length of divisorial extremal rays
(see Appendix by V.V. Shokurov). Thus, one should consider the
results about Calabi-Yau 3-folds above
as our common results with V.V. Shokurov.

   We also discuss generalization of results above to so called
$\Bbb Q$-factorial models of Calabi-Yau 3-folds,
which sometimes permits to involve non-polyhedral case
and small extremal rays to the game.

\endabstract

\endtopmatter

\document

\rightheadtext{Diagram Method and Calabi-Yau 3-folds}

\leftheadtext{Viacheslav V. Nikulin}

\centerline{\bf With Appendix by Vyacheslav V. Shokurov:}
\centerline{\bf Anticanonical Boundedness for Curves.}
\vskip30pt

\head
0. Introduction
\endhead

We consider algebraic projective varieties over the field $\Bbb C$ of
complex numbers.

In our paper \cite{N8}, we developed for
Fano $3$-folds $X$ with terminal $\Bbb Q$-factorial singularities
so called Diagram Method (for divisorial case).
As an application, we proved that if Picard number $\rho (X)>7$,
then $X$ either has a small extremal ray or a $nef$ rational element
$h$ with $h^3=0$.

Here, we generalize this method for arbitrary $3$-folds (also for
divisorial case). It was very surprising for us that this is possible.
See Basic Theorems 1.3.2 and 3.2.

As an application, we use this method for Calabi-Yau manifolds.
We refer to very important papers of
P.M.H. Wilson \cite{W1}, \cite{W2}
about terminology and basic results on Calabi-Yau $3$-folds we have
to use (we cite these results in Section 4.1).
Using one result of V.V. Shokurov about length of divisorial
extremal rays for log-terminal situation (see Appendix by V.V. Shokurov),
we get the following result which one should
consider as a common result of the author and V.V. Shokurov:

\proclaim{Theorem 0.1 (by the author and V.V. Shokurov)}
Let $X$ be a $3$-dimensional Calabi-Yau
manifold and the Picard number $\rho (X) > 40$.

Then one of two cases (i) or (ii) below holds:

(i) There exists a small extremal ray on $X$.

(ii) There exists a $nef$ element $h$ such that $h^3=0$
(thus, the $nef$ cone $NEF(X)$ and the cubic intersection hypersurface
${\Cal W}_X$ have a common point; here, we don't claim that $h$
is rational!).
\endproclaim

As a corollary, we get

\proclaim{Theorem 0.1' (by the author and V.V. Shokurov)}
Let $X$ be a $3$-dimensional Calabi-Yau
manifold. Assume that the $nef$ cone $NEF(X)$ (equivalently,
Mori cone $\nem (X)$) is finite polyhedral
and $X$ does not have a small extremal ray.

Then there exists a rational $nef$ element $h$ with
$h^3=0$ if $\rho(X)>40$.
\endproclaim

\demo{Proof} See Theorems 4.3.1 and 4.3.1' for more exact
statements and the proof.
\enddemo

It seems that existence of a $nef$ rational element $h$
with $h^3=0$ is very important for Calabi-Yau $3$-folds.
See \cite{D-Gro}, \cite{Gro}, \cite{Gra}, \cite{Hu}, \cite{O},
\cite{W1}, \cite{W2}, \cite{W3}.
By I.I. Piatetsky-Shapiro and I.R. Shafarevich
\cite{P\u S-\u S}, an algebraic $K3$ surface has a
$nef$ rational element with the square zero if its Picard number $\ge 5$.
See Sect. 6 where we discuss these results and their connection
with our results.

In Sect. 5 we consider one possibility to extend our results for
cases when either Mori cone $\nem (X)$ is not finite polyhedral
or $X$ has small extremal rays. It is connected with
considering of so called {\it $\Bbb Q$-factorial models $Y$
of a Calabi-Yau manifold $X$} which one gets
as a sequence of contractions
of divisorial extremal rays and flops
in small extremal rays. These models have $\Bbb Q$-factorial
canonical singularities, and we can apply Diagram Method to these
models obtaining results similar to Theorems 0.1 and 0.1'
(with replacing of the constant $40$ by another constant).
For some very special class of $\Bbb Q$-factorial
models (we name them {\it very good})
we prove results similar to
Theorems 0.1 and 0.1' (with the constant $163$ instead $40$).
See Theorem 5.5 and Corollaries 5.6 and 5.9.
We conjecture that analogs
of Theorems 0.1 and 0.1' are valid for arbitrary
$\Bbb Q$-factorial models (one should replace $40$ by
another absolute constant). See Conjecture 5.3.

In particular, this preprint contains
results we announced in \cite{N10}.

I am grateful to I. Dolgachev, A. Grassi, M. Gross,
M. Reid, D. Morrison, V.V. Shokurov
for useful discussions.
I am grateful to Professor I.R. Shafarevich for his support of
these my studies.

\head
1. One effective variant of the diagram method for the divisorial case
\endhead

\subhead
1.1. Reminding
\endsubhead
First, we recall one combinatorial result of \cite{N8}.

Let $X$ be a projective algebraic variety with $\Bbb Q$-factorial
singularities over an
algebraically closed field. Let
$N_1(X)$ be the $\Bbb R$-linear space generated by all
algebraic  curves on  $X$ by the numerical equivalence,
and let $N^1(X)$ be the  $\Bbb  R$-linear space generated
by all Cartier (or Weil) divisors  on X by the  numerical
equivalence. Linear spaces $N_1(X)$ and $N^1(X)$ are dual to one
another by the intersection pairing.
Let $NE(X)$ be a convex cone in $N_1(X)$ generated by
all effective curves on $X$.
Let $\nem (X)$ be the closure of the
cone $NE(X)$  in $N_1(X)$. It  is called
{\it Mori cone ({\rm or} polyhedron)} of $X$.
A non-zero element  $x \in N^1(X)$ is called  $nef$ if
$x\cdot \overline {NE}(X) \ge 0$. Let $NEF(X)$  be  the
set of all $nef$ elements of $X$ and the zero.  It  is the
convex cone in $N^1(X)$  dual to Mori cone $\overline
{NE}(X)$.  A ray  $R\subset \overline {NE}(X)$ with origin 0
is called {\it extremal}  if  from $C_1 \in \overline
{NE}(X)$, $C_2\in \overline {NE}(X)$  and $C_1+C_2\in R$
it  follows that $C_1  \in R$ and $C_2    \in R$.

We consider a condition (i) for
a set $\Cal R$ of extremal rays on X.

(i) \it If $R \in \Cal R$, then all curves $C \in R$ fill out
an irreducible  divisor $D(R)$  on $X$. We call this extremal ray
divisorial.

\rm
In this case, we can correspond to $\Cal R$ (and subsets of $\Cal
R$) an
oriented graph $G({\Cal R})$ in the following way: Two different
rays $R_1$ and
$R_2$ are joined by an arrow $R_1R_2$ with the beginning in $R_1$
 and the
end  in $R_2$ if $R_1\cdot D(R_2)>0$. Here and in what follows,
for an extremal
ray  $R$ and a divisor $D$ we  write $R\cdot D>0$  if
$r\cdot D>0$
for $r\in R$ and $r\not= 0$. (The same for the symbols $\le$, $\ge$
 and $<$.)

A set $\Cal E$ of extremal rays is called {\it extremal}
if it is contained in  a
face of ${\overline {NE}(X)}$. Equivalently, there exists a nef
element $H \in N^1(X)$ such that ${\Cal E}\cdot H=0$. Evidently,
a subset of an extremal set is extremal too.

We consider the following condition (ii) for extremal sets
$\E$ of extremal rays.

(ii) \it An extremal set
$\E = \{ R_1,...,R_n \}$ satisfies the condition
(i), and for any real numbers $m_1\ge 0,....,m_n\ge 0$ which are not
all equal to $0$,
there exists a ray $R_j\in \Cal E$ such that
$R_j\cdot (m_1D(R_1)+m_2D(R_2)+...+m_nD(R_n))<0$. In particular,
the effective divisor $m_1D(R_1)+m_2D(R_2)+...+m_nD(R_n)$ is not
$nef$.

\rm
A set $\La $ of extremal rays is called {\it $E$-set} (extremal in a
different sense) if the $\Cal L$ is
not extremal but every proper subset of $\Cal L$ is extremal.
Thus, $\La$ is a minimal non-extremal set of extremal rays.
Evidently, an $E$-set $\La$ contains at least two elements.

We consider the following condition (iii) for $E$-sets
$\La $.

(iii) \it Any proper subset of an $E$-set
$\La=\{ Q_1,...,Q_m\}$ satisfies the
condition (ii), and
there exists a non-zero effective nef divisor
$D({\Cal L})=a_1D(Q_1)+a_2D(Q_2)+...+a_mD(Q_m)$.

\rm
We have the following statement:

\proclaim{Lemma  1.1.1}
An $E$-set $\La$ satisfying the condition (iii) is  connected in the
following sense: For any decomposition
${\Cal L}={\Cal L}_1\coprod {\Cal L}_2$, where
${\Cal L}_1$ and ${\Cal L}_2$ are non-empty, there  exists an
arrow  $Q_1Q_2$ such that $Q_1 \in {\Cal L}_1$ and $Q_2 \in
{\Cal L}_2$.

If ${\Cal L}$ and ${\Cal M}$ are two different
$E$-sets satisfying the condition (iii), then there exists an arrow
 $LM$ where
$L\in {\Cal L}$ and $M\in {\Cal M}$.
\endproclaim

\demo{Proof} See \cite{N8, Lemma 1.1}.
\enddemo

Let $NEF(X)=\overline {NE}(X)^\ast \subset N^1(X)$ be the cone of nef
elements of X and ${\Cal M}(X)=NEF(X)/{\Bbb R}^+$ its projectivization.
We use usual relations of orthogonality between subsets of
$\M(X)$ and $\nem(X)$. So, for
$U\subset \M(X)$ and $V\subset \nem(X)$ we write
$U\perp V$ if
$x\cdot y =0$ for any ${\Bbb R}^+x\in U$ and any $y\in V$. Thus,
for $U \subset \M(X)$, $V\subset \nem(X)$ we denote
$$
U^\perp = \{ y\in \nem(X) \mid U\perp y\},\ \
V^\perp = \{ x\in \M(X)\mid x\perp V\}.
$$

A subset $\gamma \subset \M(X)$ is called
a {\it face} of $\M(X)$ if there
exists a non-zero element $r \in \nem(X)$ such that
$\gamma =r^\perp$.
Similarly, a subset $\alpha \subset \nem (X)$ is called
a {\it face} of $\nem (X)$ if there
exists a non-zero element $h \in \M (X)$ such that
$\alpha =h^\perp$.

A convex set is called {\it a closed polyhedron} (equivalently,
{\it finite polyhedral}) if it is a
convex hull of a finite set of points. A convex closed polyhedron is called
{\it simplicial} if all its proper faces are simplexes. A convex
closed polyhedron is called {\it simple} (equivalently,
it has {\it simplicial angles}) if it is dual to a simplicial one.
In other words, any its face of codimension
 $k$ is contained exactly in $k$ faces of
the highest dimension. Evidently, $\M (X)$ is simple if the Mori cone
$\nem (X)$ is polyhedral (has a finite set of extremal rays)
and is simplicial (all proper faces of $\nem (X)$ are cones over simplexes).
Evidently, the last property is equivalent to the fact that any extremal set
of extremal rays on $X$ is linear independent.

Let $A, B$ are two vertices of an oriented graph $G$. The
{\it distance} $\rho (A,B)$ in $G$ is a length (the number of links)
of a shortest oriented path of
the graph $G$ with the beginning in $A$ and
the end in $B$. The distance is $+\infty $ if this path does not
exist.  The {\it diameter} diam $G$ of an oriented graph $G$ is
the maximum distance between
ordered pairs of its
vertices. By the Lemma 1.1.1, the diameter of an $E$-set is a finite
number if this set satisfies the condition (iii).

We have the following

\proclaim{Theorem 1.1.2}
Let $X$ be a projective algebraic variety
with $\Bbb Q$-factorial singularities and $\dim~X \ge2$.
Let us suppose that
${\Cal M}(X)$ is closed and simple (equivalently,
Mori cone $\nem (X)$ is a finite polyhedral simplicial cone).

Assume that all extremal ray on $X$ are divisorial
(satisfies the condition (i) above), each extremal subset of extremal rays
satisfies the condition (ii), and each $E$-set of
extremal rays satisfies the
condition (iii). Assume
that there are some
constants $d, C_1, C_2$ such that the conditions (a) and (b)
below hold:

(a)
$$
\text{diam}~G({\Cal L})\le d.
$$
for any $E$-set of extremal rays on $X$.

(b)
$$
\sharp \{ (R_1, R_2)\in \E \times \E
\mid 1 \le \rho (R_1,R_2)\le d\} \le C_1 \sharp \E ;
$$
and
$$
\sharp \{ (R_1, R_2)\in \E \times \E
\mid d+1\le \rho (R_1,R_2) \le 2d+1\} \le C_2 \sharp \E .
$$
for any extremal set $\E$ of extremal rays on $X$.

Then $\dim N_1(X)=\dim \nem (X)\le (16/3)C_1+4C_2+6$.

\endproclaim

\demo{Proof} This is a particular case of  \cite{N8, Theorem 1.2}.
\enddemo

\subhead
1.2. General results on divisorial extremal rays for 3-folds
\endsubhead

In fact, the most part of results here was contained in \cite{N8}.

We restrict considering normal
projective 3-folds $X$ with $\bold Q$-factorial
singularities.

Let $R$ be an extremal ray of Mori polyhedron $\overline {NE}(X)$
 of $X$.
A morphism $f:X\rightarrow Y$ on a normal projective
variety $Y$ is called the {\it contraction} of the ray $R$ if for
 an irreducible curve $C$ of $X$ the
image $f(C)$ is a point iff $C\in R$. The contraction $f$
is defined by a linear system $H$ on $X$ ($H$ gives
 the nef element of $N^1(X)$, which we denote by $H$ also). It
follows that an irreducible curve $C$ is contracted iff
$C\cdot H=0$. We assume that the contraction $f$ has properties:
$f_\ast {\Cal O}_X={\Cal O}_Y$ and the sequence
$$
0\rightarrow {\bold R}R\rightarrow N_1(X)\rightarrow N_1(Y)
\rightarrow 0
\tag1-2-1
$$
is exact where the arrow $N_1(X)\rightarrow N_1(Y)$ is $f_\ast$.
An extremal ray $R$ is called {\it contractible} if there exists
its contraction $f$ with these properties.

The number $\kappa (R)=\dim Y$ is called {\it Kodaira dimension}
of the contractible extremal ray $R$.

We recall (see above) that a
subset $\gamma$ of $\nem (X)$ is called a {\it face} if there
exists a non-zero $nef$ element $r\in \nem (X)$ such that
$\gamma =r^\perp$.
A face $\gamma $ of $\nem (X)$ is called {\it contractible}
if there exists a morphism $f:X\rightarrow Y$ on a normal projective
variety $Y$ such that $f_\ast \gamma =0$,
$f_\ast {\Cal O}_X={\Cal O}_Y$ and $f$  contracts  curves
lying in $\gamma $ only. The $\kappa (\gamma )=\dim Y$ is called
{\it Kodaira dimension of} $\gamma $.

Let $H$ be a general nef element orthogonal to a face $\gamma $
of Mori polyhedron. {\it Numerical Kodaira dimension of} $\gamma$ is
defined by the formula
$$
\kappa_{num}(\gamma )=
\cases
3, &\text{if $H^3>0$;}\\
2, &\text{if $H^3=0$ and $H^2\not\equiv 0$;}\\
1, &\text{if $H^2\equiv 0$.}
\endcases
$$
It is obvious that for a contractible face $\gamma $ we have
$\kappa_{num}(\gamma )\ge \kappa (\gamma )$. In particular,
$\kappa_{num}(\gamma )=\kappa (\gamma )$ for a contractible face
$\gamma $ of Kodaira dimension $\kappa (\gamma )=3$.

We will use the following statement which (in different variants)
is standard:

\proclaim{Proposition 1.2.1}
Let $X$ be a projective 3-fold with
$\bold Q$-factorial singularities,
\linebreak
$D_1,...,D_m$ irreducible
divisors on $X$ and $f:X\rightarrow Y$ a surjective morphism such
 that
$\dim X=\dim Y$ and $\dim f(D_i)<\dim D_i$.
 Let $y\in f(D_1)\cap ...\cap f(D_m)$.

Then there are $a_1>0,\ ...,\ a_m>0$ and an open $U$, $y\in U\subset
f(D_1)\cup ...\cup f(D_m)$, such that
$$
C\cdot (a_1D_1+...+a_mD_m)<0 $$
if a curve $C\subset D_1\cup ...\cup D_m$ belongs to a non-trivial
algebraic family of curves on $D_1\cup ...\cup D_m$
and $f(C)=point\  \in U$.
\endproclaim

\demo{Proof} See Proposition 2.2.6 in  \cite{N8}
\enddemo

By this Proposition, we have

\proclaim{Lemma 1.2.2}  Let $R$ be a contractible extremal ray of
Kodaira dimension 3 and $f:X\rightarrow Y$ its contraction.

Then there are three possibilities:

(I) All curves $C\in R$ fill an irreducible Weil divisor $D(R)$,
the contraction $f$ contracts $D(R)$ in a point and $R\cdot D(R)<0$.

(II) All curves $C\in R$ fill
an irreducible Weil divisor $D(R)$, the
contraction $f$ contracts $D(R)$ on an irreducible curve and
$R\cdot D(R)<0$.

(III) (small extremal ray)
All curves $C\in R$ give a finite set of irreducible curves
and the contraction $f$ contracts these curves in points.
\endproclaim

\demo{Proof} Assume that some curves of $R$ fill
an irreducible divisor $D$. Then $R\cdot D<0$,
by Proposition 1.2.1. Suppose
that $C\in R$ and $D$ does not contain $C$. It follows that
$R\cdot D \ge 0$. We get a contradiction. It follows the Lemma.
\enddemo

According to the Lemma 1.2.2, we say that an extremal ray $R$ has the
{\it type (I), (II) or (III) (small)} if it is contractible of Kodaira
dimension 3 and the statements (I), (II) or (III) respectively hold.
Extremal rays of the type (I) and (II) we also call {\it divisorial}.

For a divisor $D$ on $X$ let

$$
\overline {NE}(X, D)=(\hbox{image\ } \nem (D)) \subset \nem (X).
$$

\proclaim{Lemma 1.2.3} Let $R$ be a divisorial extremal ray of the type (I).

Then $\nem (X, D(R))$ is the ray $\nem (X, D(R))=R$.

Let $R$ be a divisorial extremal ray of the type (II),
and $f$ its contraction.

Then $\nem (X, D(R))$ is an angle $\nem (X,D(R))=R + {\bold R}^+S$  with
edges $R$ and $S$, where
$f_\ast {\bold R}^+S={\bold R}^+(f(D)).$
\endproclaim

\demo{Proof}
This follows at once from the exact sequence (1-2-1).
\enddemo

Using Lemma 1.2.3, we get Lemmas 1.2.4 and 1.2.5 below.

\proclaim{Lemma 1.2.4}
Let $R_1$ and $R_2$ are two different
extremal rays of the type (II) such that the divisors $D(R_1)=D(R_2)$.

Then for
$D=D(R_1)=D(R_2)$ we have:
$$
\overline {NE}(X, D)=R_1+R_2.
$$
In particular, do not exist three different extremal rays of
the type (II) such that their divisors are coincided.
\endproclaim

If for two different extremal rays $R_1, R_2$ of the
type (II), $D(R_1)=D(R_2)$
(thus, we have the case of the Lemma 1.2.4 above), we
say that the set $\{ R_1, R_2\}$ of extremal rays has
{\it the type $\bb_2$}.

\proclaim{Lemma 1.2.5}
The divisors $D(R_1)$ and $D(R_2)$ of two different extremal rays of
the type (I) do not intersect one another.

The divisors of an extremal ray of the type (I) and a pair
of the type $\bb_2$ do not intersect one another.

The divisors of two different pairs of the type
$\bb_2$ do not intersect one
another.
\endproclaim

The next Lemma was proved in \cite{N8, Theorem 2.3.3}, but we give
the proof since this statement is very important.

\proclaim{Lemma 1.2.6}
Suppose that a pair $\{ R_1, R_2\}$ has the type $\bb_2$. Then
$$
\nem (X, D(R_1)) = \nem (X, D(R_2)) = R_1+R_2
$$
is a 2-dimensional
face of Mori polyhedron of the numerical Kodaira dimension $3$ and
such that $(R_1+R_2)^\perp$ is a face of the $NEF(X)$ of the
codimension $2$.
\endproclaim

\demo{Proof} Since the rays $R_1,R_2$ are extremal of
Kodaira dimension $3$, there are $nef$ elements $H_1, H_2$ such that
$H_1\cdot R_1=H_2\cdot R_2=0, {H_1}^3>0, {H_2}^3>0$. Let
$0\not=C_1\in R_1$ and $0\not=C_2\in R_2$. Let $D$ be a divisor of the rays
$R_1$ and $R_2$. Let us consider a map
$$
(H_1, H_2) \to H=
\tag1-2-2
$$
$$
=(-D\cdot C_2)(H_2\cdot C_1)H_1+(-D\cdot C_1)(H_1\cdot C_2)H_2+
(H_2\cdot C_1)(H_1\cdot C_2)D.
$$

For a fixed $H_1$, we get a linear map $H_2\to H$ of the set of $nef$
elements $H_2$ orthogonal to $R_2$ into the set of $nef$ elements $H$
orthogonal to $R_1$ and $R_2$. This map has a one dimensional kernel,
generated by
$(-D\cdot C_2)H_1+(H_1\cdot C_2)D$. It follows that
$R_1+R_2$ is a $2$-dimensional face of $\nem(X)$.

For a general $nef$ element $H=a_1H_1+a_2H_2+bD$ orthogonal to this face,
where $a_1, a_2, b>0$, we have
$H^3=(a_1H_1+a_2H_2+bD)^3\ge (a_1H_1+a_2H_2+bD)^2\cdot
(a_1H_1+a_2H_2)=(a_1H_1+a_2H_2+bD)\cdot (a_1H_1+a_2H_2+bD)\cdot
(a_1H_1+a_2H_2)\ge (a_1H_1+a_2H_2)^2\cdot (a_1H_1+a_2H_2+bD) \ge
(a_1H_1+a_2H_2)^3>0$, since
$a_1H_1+a_2H_2+bD$ and $a_1H_1+a_2H_2$ are $nef$. It follows that the face
$R_1+R_2$ has the numerical Kodaira dimension $3$.

The last statement follows from construction.
\enddemo

By Proposition 1.2.1, we have

\proclaim{Lemma 1.2.7} Let $\E = \{R_1,...,R_n\} $ be a
set of divisorial extremal rays of the type (I) or (II) on $X$
and the $\E$ is contained in a face of $\nem (X)$ of
Kodaira dimension 3,

Then there are real $a_1,..., a_n$ such that
$$
R_i\cdot (a_1D(R_1)+\cdots a_nD(R_n))<0
\tag1-2-3
$$
for all $R_i\in \E$.

It follows that $\E$ is linear independent if all divisors
$D(R_1),..., D(R_n)$ are different.
\endproclaim

\demo{Proof} We only need proving the last statement.
Let us assume that we have a linear dependence condition
$c_1R_1+\cdots c_sR_s+c_{s+1}R_{s+1}+\cdots+c_nR_n=0$.
We can suppose that $c_1,..., c_s$ are positive and
$c_{s+1},...,c_n$ are non-positive. By (1-2-3) and
our condition,
$$
(c_1R_1+\cdots + c_sR_s)\cdot (a_1D(R_1)+\cdot + a_sD(R_s))<0
$$
and
$$
(c_{s+1}R_{s+1}+ \cdots +c_nR_s)\cdot (a_1D(R_1)+\cdot + a_sD(R_s))\le 0.
$$
We get the contradiction.
\enddemo

We have the following inverse statement based on standard arguments
connected with Perron-Frobenius Theorem.

\proclaim{Lemma 1.2.8} We have the following inverse statement
to the previous one:
Let $\E = \{R_1,...,R_n\} $ be a
set of divisorial extremal rays of the type (I) or (II) on $X$
and all divisors $D(R_1),..., D(R_n)$ are
different. Assume that there are positive $a_1,..., a_n$ such that
$$
R_i\cdot (a_1D(R_1)+\cdots a_nD(R_n))<0
$$
for all $R_i\in \E$.

Then for any
$b_1\ge 0,..., b_n\ge 0$ which are not all equal to 0,
there exists $1\le j \le n$ such that
$$
R_j \cdot (b_1D(R_1)+ \cdots +b_nD(R_n))<0.
$$
For each  $1\le k \le n$, there are non-negative
$u_{1k},...,u_{nk}$ such that
$R_k\cdot (u_{1k}D(R_1)+ \cdots + u_{nk}D(R_n))<0$ and
$R_j\cdot (u_{1k}D(R_1)+ \cdots +u_{nk}D(R_n)) = 0$ if
$j\not =k$ and $1\le j \le n$.
It particular, the $\E$ is linear independent.
\endproclaim

\demo{Proof} We can find $\lambda >0$ such that
$\lambda a_i \ge b_i$
for all $1\le i \le n$ and one of these inequalities is
equality. Suppose that
$\lambda a_j=b_j$ for $1\le j\le n$. Then
$$
\split
R_j\cdot (b_1D(R_1)+ \cdots +b_nD(R_n))&=
R_j\cdot \lambda(a_1D(R_1)+ \cdots +a_nD(R_n))\\
&+R_j\cdot ((b_1-\lambda a_1)D(R_1)+ \cdots +(b_n-\lambda a_n)D(R_n))<0.
\endsplit
$$

Evidently,
$R_i\cdot (a_1D(R_1)+\cdots a_{n-1}D(R_{n-1}))<0$
for $1\le i \le n-1$. Let us choose generators $r_i\in R_i$.
Since $r_i\cdot D(R_n)\ge 0$ for $1\le i \le n-1$, using
induction, we can find $u_1\ge 0,...,u_{n-1}\ge 0$ such that
$-r_i\cdot D(R_n)=r_i\cdot (u_1D(R_1)+\cdots u_{n-1}D(R_{n-1}))$
for all $1\le i \le n-1$. Thus,
$u_1D(R_1)+\cdots u_{n-1}D(R_{n-1})+D(R_n)$ is orthogonal to
$R_1,...,R_{n-1}$. By the statement proved above, then
$R_n\cdot (u_1D(R_1)+\cdots u_{n-1}D(R_{n-1})+D(R_n))<0$.

This finished the proof.
\enddemo

Using Lemma 1.2.8, we get

\proclaim{Lemma 1.2.9}
Let $\E = \{R_1,...,R_n\} $ be a
set of divisorial extremal rays of the type (I) or (II) on $X$
and all divisors $D(R_1),..., D(R_n)$ are
different. Assume that there are positive $a_1,..., a_n$ such that
$$
R_i\cdot (a_1D(R_1)+\cdots +a_nD(R_n))<0
$$
for all $R_i\in \E$. Let us
additionally suppose that
$C\cdot D(R_i) \ge 0$
for any curve $C\subset D(R_i),\ C\notin R_i$ and all $1\le i \le n$.

Then $R_1+\cdots +R_n$ is a face of the dimension $n$
(thus, it is a cone over $n$-dimensional simplex)  and of the numerical
Kodaira dimension $3$ of $\nem (X)$ and such that
$(R_1+\cdots +R_n)^\perp$ is a face of the cone $NEF(X)$
of the codimension $n$.
\endproclaim

\demo{Proof} Let $H$ be a $nef$ element on $X$. By
Lemma 1.2.8, there are non-negative linear functions
$b_1(H),...,b_n(H)$
such that
$H^\prime = H+\sum_{i=1}^{i=n}{b_i(H)D(R_i)}$
is orthogonal to $R_1,...,R_n$. By additional condition,
$H^\prime$ is $nef$.
The map
$H\to H^\prime$ gives a linear map from
the set of $nef$ elements on $X$ to the set of
$nef$ elements orthogonal to $R_1,...,R_n$.
This map has $n$-dimensional kernel generated by
$D(R_1),..., D(R_n)$. If follows that
the extremal rays $R_1,..., R_n$ belong to
a face of $\nem (X)$ of dimension $\le n$.
By Lemma 1.2.8, $R_1,..., R_n$ are linear
independent. If follows that $R_1,..., R_n$
belong to a face $\gamma$ of $\nem (X)$
of the dimension $n$.
By induction, we can suppose that any $n-1$-element
subset of $\{ R_1,...,R_n\}$ generates a face of
$\nem (X)$ which is a cone over $n-1$-dimensional simplex.
It follows that the $\gamma =R_1+ \cdots +\R_n$ is a cone
over $n$-dimensional simplex.

Like above, one can prove that $(H^\prime )^3 \ge H^3>0$
for an ample $H$. Thus, the face
$R_1+\cdots +\R_n$ has Kodaira dimension $3$.
We get the last property by the construction.
\enddemo

Using Perron-Frobenius Theorem, we get

\proclaim{Lemma 1.2.10} Let $\{ R_1,..., R_n\}$ be  a set of
divisorial extremal rays of the type (I) or (II) and with
different divisors $D(R_1),...,D(R_n)$.
Let us suppose that any its proper subset satisfies conditions of
Lemma 1.2.8 but the set  $\{ R_1,..., R_n\}$ itself does not.

Then there exist positive $a_1,...,a_n$ such that
$$
R_i\cdot (a_1D(R_1)+\cdots+a_nD(R_n))\ge 0
\tag1-2-4
$$
for all $1\le i\le n$.
Additionally, we have one of cases:

(a)
$R_i\cdot (a_1D(R_1)+\cdots +a_nD(R_n))=0$ for all
$1\le i\le n$. Then the set of positive $(a_1,...,a_n)$ with
the property (1-2-4), is defined uniquely up to multiplication on
$\lambda >0$.

(b) There exists $1\le j \le n$ such that  the inequality (1-2-4)
is strong for this $j$:
$$
R_j \cdot (a_1D(R_1)+\cdots +a_nD(R_n)) > 0.
$$
\endproclaim

\demo{Proof} Let us suppose that after changing numeration, for
some $1 \le m < n$, we have
$R_i\cdot D(R_j)=0$ for all $1\le i \le m < j\le n$.
By our condition, we can find
$a_1>0,...,a_m>0$ and $a_{m+1}>0,...,a_n>0$ such that
$$
R_i\cdot (a_1D(R_1)+\cdots +a_mD(R_m))<0
$$
for all $1\le i \le m$, and
$$
R_j\cdot (a_{m+1}D(R_{m+1})+\cdots +a_nD(R_n))<0
$$
for all $m+1\le j \le n$.
Evidently, for small $\epsilon >0$, we then have
$$
R_k\cdot \epsilon (a_1D(R_1)+\cdots +a_mD(R_m))+
(a_{m+1}D(R_{m+1})+\cdots +a_nD(R_n))<0
$$
for all $1\le k\le n$.
We get the contradiction.

Thus, the subdivision above is impossible. Then, by Perron--Frobenius
Theorem, there are positive $a_1,...,a_n$ with the property (1-2-4) above.

Assume that we have the case (a). Let us suppose that there are positive
$b_1,...,b_n$ such that $R_i\cdot (b_1D(R_1)+\cdots +b_nD(R_n))\ge 0$ for
all $1\le i \le n$. There exists $\lambda >0$ such that
$b_i-\lambda a_i\ge 0$ and one of these inequalities is equality. Then
$$
\split
R_i\cdot (b_1D(R_1)+\cdots +b_nD(R_n))&=\\
\lambda R_i\cdot (a_1D(R_1)+\cdots +a_nD(R_n))&+
R_i\cdot ((b_1-\lambda a_1)D(R_1)+\cdots +(b_n-\lambda a_n)D(R_n))\ge 0.
\endsplit
$$
for all $1\le i \le n$.
If at least one $b_i-\lambda a_i>0$, we then get a contradiction with
Lemma 1.2.8.
\enddemo

 From Lemma 1.2.10, we get

\proclaim{Lemma 1.2.11}
Let $\{ R_1,..., R_n\}$ be  a set of
divisorial extremal rays of the type (I) or (II) and with
different divisors $D(R_1),...,D(R_n)$.
Let us suppose that any its proper subset satisfies the condition of
Lemma 1.2.8 but the set  $\{ R_1,..., R_n\}$ itself does not.
Let us
additionally suppose that
$C\cdot D(R_i) \ge 0$
for any curve $C\subset D(R_i),\ C\notin R_i$ and all $1\le i \le n$.

Then, in notation of Lemma 1.2.10, the element
$H=a_1D(R_1)+\cdots +a_nD(R_n)$ is $nef$.
The set $\{ R_1,...,R_n\}$ is not contained in a face of $\nem (X)$ of
Kodaira dimension $3$. For the case (a) of Lemma 1.2.10, the set
$\{ R_1,...,R_n\}$ is extremal.

\endproclaim

\demo{Proof}
By additional condition,  the $H$ is $nef$.
The set $\{ R_1,...,R_n\}$ is not contained in a face
of $\nem (X)$ of Kodaira dimension $3$ by Proposition 1.2.1 and
Lemma 1.2.8.

If we have the case (a) of
Lemma 1.2.10, then $\{ R_1,...,R_n\}$ is contained in the face of
$\nem (X)$ orthogonal to $H$. This finishes the proof.
\enddemo

As a particular case, we get

\proclaim{Lemma 1.2.12}
Assume that all finite polyhedral faces $\gamma$ of $\nem (X)$ with
the property $\text{codim}~\gamma ^\perp=\dim \gamma$ are contractible
and their numerical Kodaira dimension is equal to Kodaira dimension
(here the $\gamma ^\perp$ is the corresponding face of $NEF(X)$).
Let $\La =\{ R_1,..., R_n\}$ be an $E$-set of divisorial
extremal rays of the type (I) or (II).
Let us
additionally suppose that
$C\cdot D(R_i) \ge 0$
for any curve $C\subset D(R_i),\ C\notin R_i$ and all $1\le i \le n$.

Then there are non-negative $a_1,..., a_n$ which are not all
equal to zero such that \linebreak
$a_1D(R_1)+\cdots +a_nD(R_n)$ is $nef$
\endproclaim

\demo{Proof} By Proposition 1.2.1 and Lemma 1.2.9,
there exists a minimal subset
$\La^\prime \subset \La$ such that for $\La^\prime$ conditions of the
Lemma 1.2.11 hold. By Lemma 1.2.11, we get the statement.
\enddemo

Using Theorem 1.1.2 and Lemmas above, we get

\proclaim{Theorem 1.2.13}
Let $X$ be a projective 3-fold
with $\Bbb Q$-factorial singularities. Let us suppose
that Mori cone $\nem (X)$ is finite polyhedral and any its
face has Kodaira dimension $3$.
Assume that all extremal rays on $X$ are divisorial of the type
(I) or (II) and $X$ does not have a pair of extremal rays of
the type $\bb_2$.

Then extremal and $E$-sets of extremal rays on $X$
satisfy the conditions (i), (ii)
and (iii) of Section 1.1, and any proper face of $\nem (X)$ is a cone
over simplex.

Assume
that there are some
constants $d, C_1, C_2$ such that the conditions (a) and (b)
below hold:

(a)
$$
\text{diam}~G(\La )\le d
$$
for any $E$-set $\La$ of extremal rays on $X$.

(b)
$$
\sharp \{ (R_1, R_2)\in \E \times \E
\mid 1 \le \rho (R_1,R_2)\le d\} \le C_1 \sharp \E ;
$$
and
$$
\sharp \{ (R_1, R_2)\in \E \times \E
\mid d+1\le \rho (R_1,R_2) \le 2d+1\} \le C_2 \sharp \E .
$$
for any extremal set $\E$ of extremal rays on $X$ (here we use distance in
the graph $G(\E ))$.

Then we have the inequality:
$$
\dim~N_1(X)=\dim~\nem (X) \le (16/3)C_1+4C_2+6.
$$
\endproclaim

\subhead
1.3. Basic Theorem
\endsubhead

\definition{Definition 1.3.1}
We say that an algebraic 3-fold $X$ belongs to the class $\LT$
if $X$ has $\Bbb Q$-factorial singularities;
each face $\gamma$ of $\nem (X)$ generated by a finite set of
divisorial extremal rays of the type (I) or (II) and
with the property
$\dim \gamma =\text{codim}~\gamma ^\perp$ for the face $\gamma ^\perp$
of $NEF(X)$, and of numerical Kodaira dimension $3$
is contractible and has Kodaira dimension $3$;
the contraction of
any sequence of extremal rays of the type (I) or (II) starting
from  $X$ gives a 3-fold with $\Bbb Q$-factorial singularities
and with the properties above.

For a 3-fold $X$ from the class $\LT$ we say that it has constants
$q(X)$, $d(X)$, $C_1(X)$, $C_2(X)$ if we have:
$$
\# \E \le q(X)
$$
for any extremal set $\E$ of Kodaira dimension $3$ (i.e. $\E$ is
contained in a face of $\nem (X)$ of Kodaira dimension $3$)
of extremal
rays the type (II) such that the graph $G(\E)$ is full
(i.e. any two rays $R_1, R_2 \in \E$ are joint by both arrows
$R_1R_2$ and $R_2R_1$); the diameter
$$
\text{diam\ }G(\La ) \le d(X)
$$
for any $E$-set $\La$ of extremal rays of the type (I) or (II)
such that any proper subset of
$\La$ is extremal of Kodaira dimension $3$
and $\La$ has the property (iii) of Sect. 1.1;
$$
\sharp \{ (R_1, R_2)\in \E \times \E
\mid 1 \le \rho (R_1,R_2)\le d(X) \} \le C_1(X) \sharp \E ;
$$
and
$$
\sharp \{ (R_1, R_2)\in \E \times \E
\mid d(X)+1\le \rho (R_1,R_2) \le 2d(X)+1\} \le C_2(X) \sharp \E .
$$
for any extremal set $\E$ of Kodaira dimension $3$
of extremal rays of the type (I) or (II) with different divisors.
\enddefinition

We want to prove the following basic  result:

\proclaim{Basic Theorem 1.3.2}
Let $X$ be a 3-fold from the class $\LT$. We assume that
Mori cone $\nem (X)$ is finite polyhedral and the
conditions (A) and  (B) below hold:

(A) The $\nem (X)$ does not have a face of Kodaira dimension $1$ or $2$;

(B) The $\nem (X)$ does not have a small extremal ray (i.e. all
extremal rays on $X$ are divisorial of the type (I) or (II)).

Then we have the following statements about $X$ with the constants
$q(X)$, $d(X)$, $C_1(X)$, $C_2(X)$ above:

(1) The $X$ does not have a pair of extremal rays of the type $\bb_2$
and Mori cone $\nem (X)$ is simplicial;

(2) The number of extremal rays of the type (I) on $X$
is not greater than $q(X)$;

(3) The $\rho (X)=\dim N_1(X) \le (16/3)C_1(X) + 4C_2(X) + 6$.

\endproclaim

\demo{Proof} Let us prove (1).
We need the following analog of \cite{N8, Lemma 2.5.7}.

\proclaim{Lemma 1.3.3} Let $X$ be a 3-fold from the class $\LT$ and
Mori cone $\nem (X)$ is finite polyhedral.

Let $\E$ be the set of all extremal rays of a proper
face $[\E]$ of $\nem (X)$.
Let
$$
\{ R_{11}, R_{12}\} \cup ...\cup \{ R_{t1}, R_{t2}\}
$$
be a set of different pairs of extremal rays of the type $\bb_2$.
Assume that $R\cdot D(R_{i1})=0$ for all $R\in \E$ and all
$i$, $1\le i \le t$.

Then there are extremal rays
$Q_1,...,Q_r$ such that the following statements hold:

(a) $r \le t$;

(b) For any $i$, $1\le i \le r$, there exists $j$, $1\le j \le t$,
such that $Q_i\cdot D(R_{j1})>0$ (in particular, $Q_i$ is different
from extremal rays of pairs of extremal rays
$\{R_{u1},R_{u2}\}$ of the type $\bb_2$);

(c) For any $j$, $1\le j\le t$, there exists an extremal
ray $Q_i$, $1\le i \le r$, such that
$$
Q_i \cdot D(R_{j1})>0;
$$

(d) The set
$\E\cup \{ Q_1,...,Q_r\}$
is extremal, and extremal rays
$\{ Q_1,...,Q_r\}$ are linearly independent.
\endproclaim

\demo{Proof} If $t=0$, we can take $r=0$.
Thus, we assume that $t\ge 1$.

Since $R_{ij}\cdot D(R_{ij})<0$,
$1\le i \le t, 1\le j\le 2$,
the set $\E$ does not contain the rays $R_{ij}$.
Let $H$ be a general
$nef$ element orthogonal to $[\E]$. Since $t\ge 1$,
there exists $a>0$
such that $H^\prime = H+aD(R_{11})$ is $nef$ and $H^\prime$ is
orthogonal to $\E$ and one of the rays $R_{11}, R_{12}$. Let this ray
be $R_{11}$. Then the set $E\cup \{ R_{11}\}$ is extremal and is contained
in a (proper) face of $\nem (X)$.
It follows,
$\dim [\E]<\dim [\E\cup \{ R_{11}\}]<\dim \nem (X)$, and
$\dim [\E] < \dim \nem (X)-1$.
Let us consider a linear subspace
$V(\E)\subset N_1(X)$ generated by all extremal rays
$\E$. By our condition, $V(\E)$ is a linear envelope of the
face $[\E]$ of $\nem (X)$.

Let us consider the factorization map
$\pi : N_1(X)\to N_1(X)/V(\E)$. Since the cone $\nem (X)$ is
polyhedral, the cone $\pi (\nem (X))$ is generated
by images of extremal rays $T$ such that the set $\E\cup \{ T \}$
is contained in  a face $[\E \cup \{ T\} ]$ of $\nem (X)$ of the dimension
$\dim [\E ]+1$. In particular, since $\dim [\E] < \dim N_1(X)-1$, the
face $[\E \cup \{ T\}]$ is proper, and the set
$\E \cup \{T \}$ is extremal.

There exists a curve $C$ on
$X$ such that $C\cdot D(R_{11})>0$. This curve $C$ (as
any element $x\in \nem (X)$) is a linear combination of extremal rays $T$
with non-negative coefficients and extremal rays from $\E$ with
real coefficients. We have $R\cdot D(R_{11})=0$ for any
extremal ray $R\in \E$. Thus,
there exists an extremal ray $T$ above
such that $T\cdot D(R_{11})>0$.
It follows that $T$ is different from extremal rays of pairs
of the type $\bb_2$. We take $Q_1=T$. By our construction,
the set $\E \cup \{ Q_1\}$ is extremal.
If $Q_1 \cdot D(R_{j1})>0$ for any $j$ such that
$1 \le j \le t$, then $r=1$, and the set $\{Q_1\}$ gives the set
we were looking for. Otherwise, there exists a minimal $j$ such that
$2 \le j \le t$ and $Q_1\cdot D(R_{j1})=0$. Then we
replace $\E$ by the set $\E_1$ of all extremal rays in the face
$[\E \cup \{Q_1\} ]$ of the dimension
$\dim [\E_1]=\dim [\E] +1$, and the set
$$
\{ R_{11}, R_{12}\} \cup ...\cup \{ R_{t1}, R_{t2}\}
$$
by
$$
\{ R_{j1}, R_{j2} \mid 1\le j \le t,\  Q_1\cdot D(R_{j1})=0\},
$$
and repeat this procedure.
\enddemo

Also, we need the following Lemma:

\proclaim{Lemma 1.3.4} Let $X$ be a 3-fold from the class $LT$ and
Mori cone $\nem (X)$ is polyhedral.
Assume that extremal rays on $X$ have the type (I) or (II).
Let $R_{11},R_{12}$ be a pair of extremal rays of the type $\bb_2$ on $X$
and $\rho (X)\ge 3$.

Then there exists an extremal ray $R$ on $X$ such that $R$ does not
belong to a pair of extremal rays of the type $\bb_2$ and
the sets $\{ R_{11}, R\}$ and  $\{ R_{12}, R\}$ generate 2-dimensional
faces of $\nem (X)$.
\endproclaim

\demo{Proof} Let  $R_{21}, R_{22}$ be another pair of extremal rays of
the type $\bb_2$. By Lemmas 1.2.5 and 1.2.6,
the extremal rays $R_{11},R_{12}$
generate a 2-dimensional face of $\nem (X)$ and
divisors $D(R_{11})$ and $D(R_{21})$ are disjoint (have empty
intersection).
Applying Lemma 1.3.3 to
$\E=\{ R_{11}, R_{12}\}$ and the set $\{ \{ R_{21}, R_{22}\}\}$
of pairs of the type $\bb_2$, we find an extremal ray $Q$ such that
$Q\cdot D(R_{21})>0$. It follows that
$D(Q)\cap D(R_{21})$ is a non-empty curve $C$. By Lemma 1.2.5, the $Q$
has the type (II).
By Lemma 1.2.6, $\nem (X,D(R_{21}))=R_{21}+R_{22}$
is a $2$-dimensional face of $\nem (X)$. By Lemma 1.2.3,
then the $2$-dimensional angle $\nem (X, D(Q))=Q+\r^+C$.
Since $\nem (X,D(R_{11}))=R_{11}+R_{12}$
is another
$2$-dimensional face of $\nem (X)$ which does not have a common
ray with $R_{21}+R_{22}$, the angle
$\nem (X,D(Q))$ does not have a common ray with the angle
$\nem (X,D(R_{11})$. Thus, $D(Q)\cap D(R_{11})=\emptyset$.
By Lemma 1.2.5, the $Q$ does not belong to a pair of the type $\bb_2$.

Let $H$ be a general $nef$ element orthogonal to the 2-dimensional face
$R_{11}+R_{12}$. Then there exists $\alpha >0$ such that the $nef$
element
$H^\prime = H+\alpha D(Q)$
is orthogonal to the set of extremal rays $R_{11},R_{12},Q$.
If an extremal ray $R$ is different from these three extremal rays,
evidently, $H^\prime\cdot R>0$ since all extremal rays on $X$ are
divisorial and $Q$ does not belong to a pair of the type $\bb_2$.
It follows that $R_{21}+R_{22}+Q$ is a simplicial 3-dimensional
face of $\nem (X)$. It follows that $Q$ is
the extremal ray we are looking for.

Thus, we can suppose that the pair $R_{11},R_{12}$ is the only pair of
extremal rays of the type $\bb_2$ on $X$.

If there exists an extremal ray $Q$ such that
$D(Q)\cap D(R_{11})=\emptyset$,
like above, $Q$ is the required extremal ray. Thus, we can suppose
that any extremal ray $Q$ which is different from $R_{11}, R_{12}$,
has the property $D(Q) \cap D(R_{11})\not= \emptyset$. In particular,
$D(R_{11})\cap D(Q)$ is a non-empty curve.

Since $\nem (X)$ is
polyhedral, there exists an extremal ray
$Q_1\notin \{R_{11},R_{12}\}$ such that
$R_{11}+Q_1$ is a 2-dimensional face of $\nem (X)$.
If $Q_1+R_{12}$ is a 2-dimensional face of $\nem (X)$,
the $Q_1$ is the desired extremal ray. Thus, we can suppose that
$Q_1+R_{12}$ is not a 2-dimensional face of $\nem (X)$.
Similarly, we can find an extremal ray $Q_2 \notin \{R_{11}, R_{12} \}$
such that $R_{12}+Q_2$ is a $2$-dimensional face of $\nem (X)$ but
$R_{11}+Q_2$ is not. Then $Q_1\not=Q_2$. We recall that besides,
we suppose that $C_1=D(R_{11})\cap D(Q_1)$ and
$C_2=D(R_{11})\cap D(Q_2)$ are non-empty curves.

We will normalize the generator $C\in T$ of a divisorial extremal ray $T$
by the condition $C\cdot D(T)=-2$.

Let $r_{ij}$ be the generator of $R_{ij}$ and $q_i$ of $Q_i$.
Let $t=q_1\cdot D(R_{11})$ and $t_1=r_{11}\cdot D(Q_1)$ and
$t_2=r_{12}\cdot D(Q_1)$.

Since $R_{11}+Q_1$ is a $2$-dimensional face of $\nem (X)$ and
all faces of $\nem (X)$ have Kodaira dimension $3$,
by Proposition 1.2.1, there are positive $a_1,a_2$ such that
$r_{11}\cdot (a_1 D(R_{11})+a_2 D(Q_1))=-2a_1+t_1a_2<0$ and
$q_1 \cdot (a_1 D(R_{11})+a_2 D(Q_1))=ta_1-2a_2<0$. Thus,
$tt_1<4$.
We claim that $t_1<t_2$.
Let us assume that $t_2\le t_1$. By inequality above,
$tt_2 < 4$.
Let $H$ be a general $nef$ element orthogonal to $R_{12}$.
Then
$$
H^\prime=H+((H \cdot q_1)/(2-(tt_2)/2))((t_2/2)D(R_{11})+D(Q_1))
$$
is a $nef$ element which is orthogonal to extremal rays
$R_{12}, Q_1$ only. We very use here the inequality $t_2\le t_1$
to check that $R_{11}\cdot H^\prime>0$.
Thus, $R_{12}+Q_1$ is a $2$-dimensional face of $\nem (X)$.
We get a contradiction. Thus, we have proved the claim: $t_1<t_2$.

Let us consider the curve $C_1=D(R_{11})\cap D(Q_1)$.
Then $C_1=u_1r_{11}+u_2r_{12}$.
Let $r_{11}\cdot r_{12}=m$ (equivalently, the
$m$ is the degree of the maps
$f_{11} \mid r_{12}:r_{12}\to f_{11}(D(R_{11}))$ and
$f_{12}\mid r_{11}:r_{11}\to f_{12}(D(R_{12}))$
where $f_{11}$ and $f_{12}$
are the contractions of $R_{11}$ and $R_{12}$ respectively).
Then
$$
(t_2:t_1)=(r_{12}\cdot D(Q_1):r_{11}\cdot D(Q_1))=
(r_{12}\cdot C_1:r_{11}\cdot C_1)=(u_1m:u_2m).
$$
Thus,
$$
C_1=u_1r_{11}+u_2r_{12},\
(u_1:u_2)=(t_2:t_1),\ t_2>t_1.
\tag1-2-5
$$

Similarly, for $Q_2$ we introduce
$s_1=r_{11}\cdot D(Q_2), s_2=r_{12}\cdot D(Q_2)$ where
$s_2<s_1$. And we have for
$C_2=D(R_{11})\cap D(Q_2)$:
$$
C_2=v_1r_{11}+v_2r_{12},\
(v_1:v_2)=(s_2:s_1),\ s_2<s_1.
\tag1-2-6
$$

By (1-2-5) and (1-2-6), $C_1\cdot D(Q_2)=u_1s_1+u_2s_\ge u_1s_1>0$,
Thus, divisors $D(Q_1)$ and $D(Q_2)$ have a common curve.
It follows that 2-dimensional angles
$\nem (X,D(Q_1))=\r^+C_1+Q_1$ and
$\nem (X,D(Q_2))=\r^+C_2+Q_2$ should have a common ray.
By (1-2-5) and (1-2-6), we then get that
rays $R_{11},R_{12},Q_1,Q_2$ generate a $3$-dimensional space $V$,
and are extremal rays of the 3-dimensional polyhedral convex cone
$\nem (X)\cap V$ with $2$-dimensional faces $R_{11}+R_{12}$,
$R_{11}+Q_1$, $R_{12}+Q_2$. Besides, angles $\r^+C_1+Q_1$ and
$\r^+C_2+Q_2$ have a common ray. One should draw the picture to
see a contradiction with (1-2-5) and (1-2-6).
This finishes the proof of Lemma.
\enddemo

Now we can prove the statement (1). Thus, suppose that $X$ satisfies the
conditions of Theorem and has a pair $\{ R_{11}, R_{12}\}$ of the type
$\bb_2$. By Lemma 1.3.4, we can find an extremal ray $Q$ which does not
belong to a pair of extremal rays of the type $\bb_2$ and
$Q+R_{11},\ Q +R_{12}$ are $2$-dimensional faces of $\nem (X)$.
Let us contract the extremal ray $Q$. By our conditions, we get a
3-fold  $X^\prime$ from the class $LT$,
without small extremal rays, with polyhedral $\nem (X^\prime)$,
and the image of the pair $R_{11}, R_{12}$
of the type $\bb_2$ will be a pair of the type $\bb_2$ again. Thus,
using the Lemma 1.3.4 $\rho (X)-2$ times,
we get a 3-fold $Y$ which satisfies
the condition of the Theorem, has $\rho (Y)=2$, and has a pair of extremal
rays $\{R_{11},R_{12}\}$ of the type $\bb_2$. Then evidently
$\nem (Y)=R_{11}+R_{12}$. We have $R_{11}\cdot D(R_{11})<0$ and
$R_{12}\cdot D(R_{12})<0$. Thus, any curve of $Y$ has negative
intersection with the effective divisor $D(R_{11})=D(R_{12})$.
We get the contradiction.

This proves the statement (1).

\smallpagebreak

Now let us prove (2): $X$ does not have more than $q(X)$ extremal rays of
the type (I).

By Proposition 1.2.1 and Lemma 1.2.5,
divisors of different extremal rays of
the type (I) do not have a common point and their number is finite.
By Lemma 1.2.9, the set of extremal rays of
the type (I) generates a simplicial face of $\nem (X)$ of Kodaira
dimension 3.
Let
$$
\{ R_1,...,R_s\}
$$
be the hole set of extremal rays of the type (I) on $X$.
We should prove that $s\le q(X)$.

We say that two divisorial extremal rays are joint (more formally:
divisorially joint) if their divisors
have a common point. It defines connected components of a set of
divisorial extremal rays.

Let $\E$ be a maximal extremal set of extremal rays on $X$
containing the set $\{ R_1,..., R_s\}$ and such that each
connected component of $\E$
contains at least one of extremal rays $R_1,...,R_s$.
Let $T$ be a connected component of the $\E$. Assume that $T$
has two different extremal rays $R_i,R_j$ from $R_1,...,R_s$.
After contracting all extremal rays
of $T$ different from $R_i,R_j$,
we get a 3-fold $Y$ from the class $\LT$, and images of extremal rays
$R_i, R_j$ give different extremal rays of the type (I) on $Y$ such that
their divisors are not disjoint. We get the contradiction with
Lemma 1.2.5.

Thus, $\E$ has exactly $s$ connected components
$T_1,..., T_s$ such that $T_i$ contains the extremal ray $R_i$.
Evidently, the maximal $\E$ does exist.

By Proposition 1.2.1, for $1 \le i \le s$, there exists
an effective divisor $D(T_i)$ which is a linear
combination of divisors of rays from $T_i$ with positive coefficients
and $R\cdot D(T_i)<0$ for any $R\in T_i$.
Evidently, the contraction of $T_i$ contracts all $T_i$ to a point.
Thus, any curve of divisors of $T_i$
belongs to the sum of extremal rays of $T_i$ with positive coefficients.
Thus, for this curve $C$ we also have $C\cdot D(T_i)<0$.
Using the divisors $D(T_i)$, similarly to Lemma 1.3.4,
we can find extremal rays
$$
\{Q_1,...,Q_r\}
$$
with properties:

(a) $r \le s$;

(b) For any $i$, $1\le i \le r$, there exists $j$, $1\le j \le t$,
such that $Q_i\cdot D(T_j)>0$ (in particular, $Q_i$ is different
from extremal rays of $\E$ and does not have the type (I));

(c) For any $j$, $1\le j\le s$, there exists an extremal
ray $Q_i$, $1\le i \le r$, such that
$$
Q_i \cdot D(T_j)>0;
$$

(d) The set
$\{ Q_1,...,Q_r\}$
of extremal rays
is extremal.

By our conditions, all extremal rays on $X$ are divisorial. Thus,
by (b), the extremal rays $Q_1,...,Q_r$ have the type (II).

Let us take the ray $Q_i$, and let $Q_i\cdot D(T_j)>0$.
By Lemma 1.2.9, the set
$T_j$ generates a simplicial face $\gamma_j$
of $\nem (X)$. We have mentioned above that
each curve of divisors of rays from $T_j$ belongs to this face.
It follows that $\nem (X, D(Q_i))$ is a 2-dimensional angle
bounded by the ray $Q_i$ and a ray from the face $\gamma_j$
since the divisor $D(Q_i)$ evidently has a common curve with
one of divisors $D(R),\ R\in T_j$. Since any two sets from
$T_1,..., T_s$
do not have a common extremal ray, the faces
$\gamma_1,...,\gamma_s$ do not have a common ray (not necessarily
extremal).
It follows that the angle $\nem (X, D(Q_i))$ does not have a common ray
with the face $\gamma_k$ for $k \not=j$.
Thus, the divisor $D(Q_i)$ does not have a
common point with divisors of rays $T_k$. It follows that
$r=s$ and we can choose
numeration $Q_1,...,Q_s$ such that $Q_i\cdot D(T_i)>0$ but
$D(Q_i)$ do not have a common point with divisors of extremal rays
$T_j$ if $j\not=i$.

Let us fix $i$, $1\le i \le s$.
By our construction, the set
$\E\cup \{Q_i\}$ has
connected components
$$
T_1,..., T_{i-1},\  T_i \cup \{ Q_i\} ,\
T_{i+1},...,T_s.
$$
By definition of $\E$, then the $\E \cup\{ Q_i\}$ is not extremal. Thus,
it contains an $E$-set (minimal non-extremal) $\La_i$ which
contains $Q_i$. By Lemmas 1.2.12 and 1.1.1, the $\La_i$ is
connected.
Thus, $\{Q_i\} \subset \La_i \subset T_i\cup \{ Q_i\}$.
Let us consider the sets $\La_1,...,\La_s$. By Lemma 1.1.1,
the $\La_i, \La_j$ are joint by arrows in both directions.
By our construction,
it follows that $Q_i, Q_j$ are joint by arrows $Q_iQ_j$ and
$Q_jQ_i$ for any $1\le i < j \le s$.
Thus, $s\le q(X)$ because the set $\{ Q_1,...,Q_s\}$ is extremal.

This finishes the proof of the statement (2).

 From the (1) and Theorem 1.2.13, the statement (3) follows.
This finishes the proof of Basic Theorem 1.3.2.
\enddemo

\smallpagebreak

In fact, Basic Theorem 1.3.2 is the generalization of similar
theorem about Fano $3$-folds with terminal $\Bbb Q$-factorial
singularities proved in \cite{N8}.
We will apply Basic Theorem 1.3.2 to
Calabi-Yau $3$-fold in Section 3.

To apply Theorem 1.3.2, one should work with extremal and $E$-sets of
extremal rays of the type (I) or (II). For Fano $3$-folds this is
done in \cite{N8}. The next Section 2 is, in fact, devoted to
generalization of this part of \cite{N8} for $3$-folds from the
class $\LT$ (see Theorem 2.2.6). Roughly speaking, we should reduce
the problem to sets of divisorial extremal rays which have
non-single arrows only. The last case is very similar to surfaces.

\head
2. Sets of divisorial extremal rays
\endhead

\subhead
2.1. Special divisorial extremal rays
\endsubhead

Here we continue geometrical study started in Section 1.2
of configurations of divisorial extremal rays of the types (I) or (II)
for $3$-folds.

Let $T$ be a set of divisorial extremal rays of the types
(I) or (II). We recall (see Section 1.1)
that we correspond
to $T$ a graph $G(T)$ as follows: we draw an arrow $R_1R_2$ from
$R_1$ to $R_2$ if $R_1\cdot D(R_2)>0$.
We draw an extremal ray of the type (I)
(respectively (II))
as a black (respectively white)
vertex.
A pair of the type $\bb_2$ we
draw as a pair of white vertices connected by a dotted line.
For this pair $\{ R_1,R_2\}$ we have $D(R_1)=D(R_2)$ and
$R_1\cdot D(R_1)<0, R_2\cdot D(R_1)<0$.
Thus, two divisorial extremal rays $R_1,R_2$
give one of pictures of
Figure 1. One should
consider a single arrow $R_1R_2$ (thus,
$R_1\cdot D(R_2)>0$, but $R_2\cdot D(R_1)=0$) as a "weak orthogonality" of
extremal rays. They are completely orthogonal (are not
joint by arrows) if their
divisors do not have a common point: this
follows from Lemma 1.2.2.

\vskip50pt


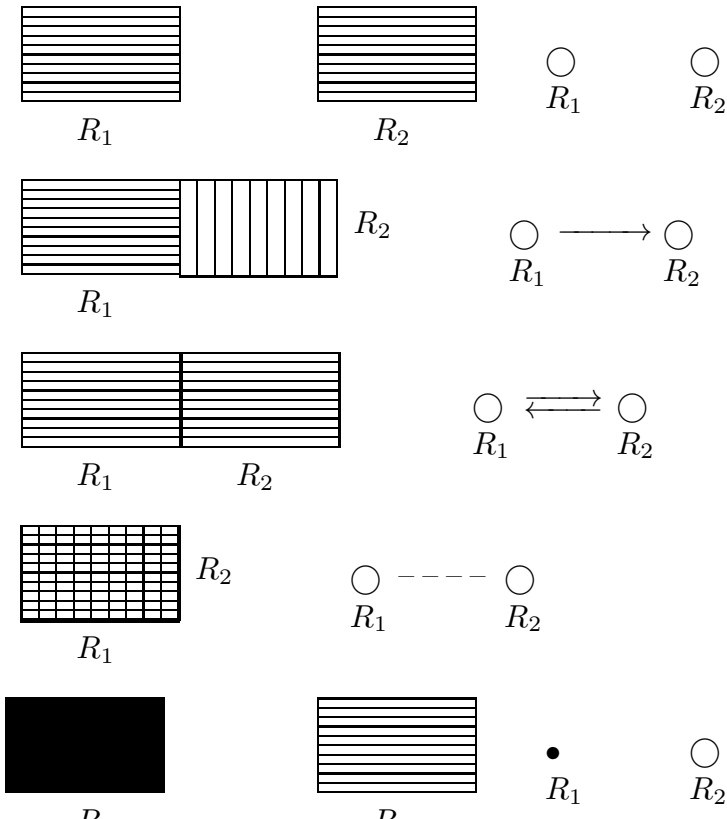
\line{\hskip10pt
\vbox{\horlines \hbox to60pt{\hfil $R_1$ \hfil}}
\hskip30pt
\vbox{\horlines\hbox to60pt{\hfil $R_2$ \hfil}}
\hskip10pt
\raise10pt\hbox{
\vbox{\vskip-10pt \hbox{$\bigcirc$} \hbox{$R_1$}}
\hskip30pt
\vbox{\vskip-10pt \hbox{$\bigcirc$} \hbox{$R_2$}}}
\hfil}

\vskip10pt

\line{\hskip10pt
\vbox{\horlines \hbox to60pt{\hfil $R_1$ \hfil}}
\hskip-12pt
\vbox{\verlines\hbox to60pt{\hfil}}
\hskip-5pt
\raise25pt\hbox{$R_2$}
\hskip30pt
\raise10pt\hbox{
\vbox{\vskip-10pt \hbox{$\bigcirc$} \hbox{$R_1$}}
\hskip0pt
\raise13pt\hbox to30pt{\rightarrowfill}
\hskip0pt
\vbox{\vskip-10pt \hbox{$\bigcirc$} \hbox{$R_2$}}}
\hfil}

\vskip10pt

\line{\hskip10pt
\vbox{\horlines \hbox to60pt{\hfil $R_1$ \hfil}}
\hskip-13pt
\vbox{\horlines\hbox to60pt{\hfil $R_2$ \hfil}}
\hskip30pt
\raise10pt\hbox{
\vbox{\vskip-10pt \hbox{$\bigcirc$} \hbox{$R_1$}}
\hskip0pt
\darr
\hskip0pt
\vbox{\vskip-10pt \hbox{$\bigcirc$} \hbox{$R_2$}}}
\hfil}

\vskip10pt

\line{\hskip10pt
\vbox{\horlines \hbox to60pt{\hfil $R_1$ \hfil}}
\hskip-62pt
\vbox{\verlines\hbox to60pt{\hfil}}
\hskip-5pt
\raise25pt\hbox{$R_2$}
\hskip30pt
\raise10pt\hbox{
\vbox{\vskip-10pt \hbox{$\bigcirc$} \hbox{$R_1$}}
\hskip0pt
\raise13pt\hbox to30pt{\hfil -- -- -- -- \hfil }
\hskip0pt
\vbox{\vskip-10pt \hbox{$\bigcirc$} \hbox{$R_2$}}}
\hfil}

\vskip10pt

\line{\hskip10pt
\vbox to42pt{\hrule height30pt width50pt depth0pt
\vfil \hbox to60pt{\hfil $R_1$ \hfil}}
\hskip30pt
\vbox{\horlines\hbox to60pt{\hfil $R_2$ \hfil}}
\hskip10pt
\raise10pt\hbox{
\vbox{\vskip-10pt \hbox{$\bullet$} \hbox{$R_1$}}
\hskip30pt
\vbox{\vskip-10pt \hbox{$\bigcirc$} \hbox{$R_2$}}}
\hfil}

\line{\hskip10pt
\vbox to42pt{\hrule height30pt width50pt depth0pt
\vfil \hbox to60pt{\hfil $R_1$ \hfil}}
\hskip-18pt
\vbox{\horlines\hbox to60pt{\hfil $R_2$ \hfil}}
\hskip30pt
\raise10pt\hbox{
\vbox{\vskip-10pt \hbox{$\bullet$} \hbox{$R_1$}}
\hskip0pt
\darr
\hskip0pt
\vbox{\vskip-10pt \hbox{$\bigcirc$} \hbox{$R_2$}}}
\hfil}

\vskip10pt

\line{\hskip10pt
\vbox to42pt{\hrule height30pt width50pt depth0pt \vfil
\hbox to60pt{\hfil $R_1$ \hfil}}
\hskip30pt
\vbox to42pt{\hrule height30pt width50pt depth0pt \vfil
\hbox to60pt{\hfil $R_2$ \hfil}}
\hskip10pt
\raise10pt\hbox{
\vbox{\vskip-10pt \hbox{$\bullet$} \hbox{$R_1$}}
\hskip30pt
\vbox{\vskip-10pt \hbox{$\bullet$} \hbox{$R_2$}}}
\hfil}

\vskip20pt

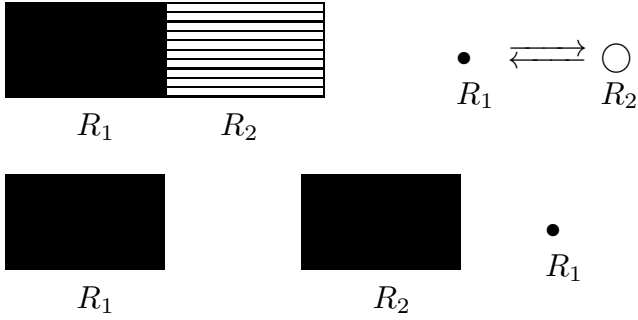
\centerline{{\bf Figure 1.}}

\vskip30pt

\definition{Definition 2.1.1}
We say that extremal rays $Q_1, Q_2,..., Q_m, \ m\ge 2,$  define
the configuration $\cc_m$ if all these extremal rays
have the type (II), extremal rays $Q_2,...,Q_m$ are divisorially
disjoint (intersection of any two divisors
$D(Q_2), D(Q_3),..., D(Q_m)$ is empty) and $Q_iQ_1$ is a single
arrow
(thus, $Q_i\cdot D(Q_1)>0,\ Q_1\cdot D(Q_i)=0$) for all
$i=2,...,m$.
\enddefinition

The following Lemma was proved in \cite{N8}. Since it is
important, we give the proof.

\proclaim{Lemma 2.1.2}
Assume that divisorial extremal rays
$Q_1,Q_2,...,Q_m$ define the cofiguration of the type $\cc_m,\ m\ge 2$.
Then $Q_i$ does not belong to a pair of
the type $\bb_2$ and $\nem (X, D(Q_i)) = Q_1+Q_i$
is the $2$-dimensional face of the numerical
Kodaira dimension $3$ of Mori polyhedron for any $2\le i \le m$.

Besides, $Q_1+ Q_2+\cdots +Q_m$ is  $m$-dimensional face of
$\nem (X)$ of numerical Kodaira dimension $3$ and such that
the face $(Q_1+ Q_2+\cdots +Q_m)^\perp $ of $NEF(X)$ has
codimension $m$.
\endproclaim

\demo{Proof} Let $2\le i \le m$.
Evidently, the curve $D(Q_i)\cap D(Q_1)$ belongs to $Q_1$.
By Lemma 1.2.3, it follows that $\nem (X, D(Q_i))=Q_1+Q_i$.
If $Q_i$ belongs to a pair of the type $\bb_2$, the angle
$Q_1+Q_i$ contains another extremal ray (different from
$Q_1$ and $Q_i$), which is impossible.

Let $H$ be a $nef$ element orthogonal to
the ray $Q_1$. Let $0\not=C_i\in Q_i$. Let us consider a map
$$
H \to H^\prime = H +
\sum_{i=2}^m {(-(H\cdot C_i)/(C_i\cdot D(Q_i)))D(Q_i)} .
\tag2-1-1
$$
It is a linear map of the set of $nef$ elements $H$ orthogonal to $Q_1$
into the set of $nef$ elements $H^\prime$ orthogonal to the rays
$Q_1, Q_2,..., Q_m$. Here the element $H^\prime$ is $nef$ because
$C\in Q_1+Q_i$ for any curve $C\subset D(Q_i)$ if $2\le i \le m$.
The kernel of the map (2-1-1) has the dimension $m-1$. It
follows that the rays $Q_1,Q_2,...,Q_m$ belong to a face of
$\nem (X)$ of a dimension $\le m$. On the other hand, multiplying
rays $Q_1,...,Q_m$ on the divisors $D(Q_1),...,D(Q_m)$, one can see very
easily that the rays $Q_1,...,Q_m$ are linearly independent.
Thus, they generate a $m$-dimensional face of $\nem (X)$. Let us show that
this face is $Q_1+Q_2+...+Q_m$.
To prove this, we show that every $m-1$ subset of $\E$ is contained in a face
of $\nem (X)$ of a dimension $\le m-1$.

If this subset contains the ray $Q_1$, this subset has the type $\cc_{m-1}$.
By induction, we can suppose that
this subset belongs to a face of $\nem (X)$ of the
dimension $m-1$. Let us consider the subset $\{Q_2, Q_3,..., Q_m\}$. Let
$H$ be an ample element of $X$. For the element $H$, the map (2-1-1) gives an
element $H^\prime$ which is orthogonal
to the rays $Q_2,...,Q_m$, but is not orthogonal
to the ray $Q_1$. It follows that the
set $\{Q_2,...,Q_m\}$ belongs to a face of the Mori polyhedron of the
dimension $< m$. Like above (see the proof of Lemma 1.2.6),
one can see that for a general $H$ orthogonal to
$Q_1$ the element $H^\prime$ has $(H^\prime )^3 \ge H^3>0$.

The last property follows from the construction.
\enddemo

We use the next statement quite often:

\proclaim{Lemma 2.1.3} Let $X$ be a $3$-fold from the
class $\LT$ (see Definition 1.3.1).
Let $R$ be a divisorial extremal ray of the
type (II), and $D_1, D_2$ are two different
irreducible Weil divisors different from
$D(R)$.

Then the curve $D_1\cap D_2$ does not belong to $R$.
\endproclaim

\demo{Proof} Let $f:X\to X^\prime$
be the contraction of the ray $R$. If
$D_1\cap D_2$ belongs to $R$, then intersection of Weil divisors
$f(D_1)$ and $f(D_2)$ is zero-dimensional. It is impossible if
$X^\prime$ is $\Bbb Q$-factorial.
\enddemo

We denote by $\p C$ the projectivization of a cone $C$ with
the beginning at zero.
Let $R$ be an extremal ray of the type (II). By Lemma 1.2.3,
the projectivization
$\p\nem (X, D(R))$ is an interval with one of its endpoints
$\p R$.

Considering lines generated by these intervals we use
the following well-known and elementary

\proclaim{Proposition 2.1.4} Let $S$ be a set of lines such that
any two lines of $S$ have a common point.

Then there are two cases:

(a) There exists a 2-dimensional plane $\Pi$ such that
each line of $S$ belongs to $\Pi$.

(b) There exists a point $P$ which belongs to each line of $S$.

The $\Pi$ and $P$ are unique if $S$ has at least two different
lines.
\endproclaim

Applying this statement, we get

\proclaim{Lemma 2.1.5} Let $\Cal S$ be a set of divisorial
extremal rays of the type (I) or (II) such that
divisors $D(R_i)$, $D(R_j)$ have a common point for
any two extremal rays from $\Cal S$ (equivalently, $R_i$ and $R_j$ are
connected by non-single, single arrow or dotted line). Let us suppose that
$\Cal S$ contains at least 3 elements. Then
$S$ has one of the following types:

(a) $\dim [{\Cal S}] =3$;

(b) Angles $\nem (X, D(R_i)),\ R_i\in {\Cal S}$, have a common ray
(not necessarily extremal) $Q$.
\endproclaim

Let $P$ be a set of divisorial extremal rays. We say that $P$ is
{\it divisorially connected}
if there does not exist a decomposition $P=P_1\cup P_2$ such
that both $P_1$ and $P_2$ are non-empty and for any $R\in P_1$ and
any $Q\in P_2$ divisors $D(R)$ and $D(Q)$ do not have a common point.
It defines {\it divisorially connected components}
of a set of extremal rays.
Also, we can say what does it mean that two sets $P_1$ and $P_2$
of divisorial extremal rays are {\it divisorially joint}: this means that
there exist extremal rays $Q_1\in P_1$ and $Q_2\in P_2$ such that
divisors $D(Q_1)$ and $D(Q_2)$ have a common point (in particular,
this divisors or even extremal rays $Q_1$, $Q_2$ may coincide).

\definition{Definition 2.1.6} A divisorial extremal ray $R$ is called
{\it special} if $R$ satisfies one of conditions:

(1) $R$ has the type (I);

(2) $R$ belongs to a pair of extremal rays of the type $\bb_2$;

(3) There exists an extremal ray $Q$ such that exactly one
pair $RQ$ or $QR$ is an arrow: thus, either $R\cdot D(Q)>0$ and
$Q\cdot D(R)=0$ or $Q\cdot D(R)>0$ and $R\cdot D(Q)=0$ (equivalently,
the set  $\{ R, Q\} $ has the type $\cc_2$).
\enddefinition

We have the following description of the special set:

\proclaim{Theorem 2.1.7} Assume that a
$3$-fold $X$ belongs to the class $\LT$.

Then the set of special divisorial
extremal rays is finite.
Its divisorially connected component is one of
the following
:

$\Aa_1$: One extremal ray $Q$ of the type (I).

$\bb_2$: Two different extremal rays $Q_{11}, Q_{12}$
of the type (II) with the same divisor
$D(Q_{11})$ $=D(Q_{12})$. Then $Q_{11}+Q_{12}$ is the $2$-dimensional
face of $\nem (X)$ of Kodaira dimension $3$.

$\cc_n,\ n\ge 2$:
$n$ extremal rays $Q_1,...,Q_n$ of the type (II),
such that divisors $D(Q_2),...,$ $D(Q_n)$ do not have a
common point, and $Q_i$ and $Q_1$ are joint by the single arrow $Q_iQ_1$ .
For this case, $Q_1+\cdots +Q_n$ is the $n$-dimensional
face of $\nem (X)$ of Kodaira dimension $3$.

$\bb_2\cc_n,\ n\ge 1$:  $n+1$ extremal rays
$\{Q_{11}, Q_2,...,Q_n, Q_{12}\}$ of
the type (II), such that $Q_1=Q_{12},...,Q_n$
define the configuration
$\cc_n$ (in notation above), and $Q_{11}, Q_{12}$ is the pair of the
type $\bb_2$. The extremal rays
$Q_i,Q_{12}$ are connected by non-single arrows for $i\ge 2$.
For this case, $Q_{11}+ \cdots +Q_n$ is the face of
$\nem (X)$ of dimension $n$ and Kodaira dimension $3$,
and $Q_{11}+Q_{12}$ is the $2$-dimensional
face of $\nem (X)$ of Kodaira dimension $3$.

$\Tt_3$ (triangle): Three extremal rays $Q_1, Q_2, Q_3$ of the type (II)
connected by single arrows $Q_1Q_2$, $Q_2Q_3$ and non-single
arrows $Q_1Q_3$ and $Q_3Q_1$.
For this case, $\nem (X,D(Q_1))=Q_1+Q_2$ and
$\nem (X,D(Q_2))=Q_2+Q_3$ are $2$-dimensional faces of $\nem (X)$ of
Kodaira dimension $3$. And $\nem (X,D(Q_3))=T+Q_3$ where
$T\subset \nem (X,D(Q_1))=Q_1+Q_2$.

$\Tt_3^\prime$ (special triangle):
Three extremal rays $Q_1, Q_2, Q_3$ of the type (II)
connected by single arrows $Q_1Q_2$, $Q_2Q_3$ and $Q_3Q_1$.
For this case, $\nem (X,D(Q_1))=Q_1+Q_2$,
$\nem (X,D(Q_2))=Q_2+Q_3$ and $\nem (X,D(Q_3))=Q_3+Q_1$
are $2$-dimensional faces of $\nem (X)$ of
Kodaira dimension $3$.

\endproclaim


\demo{Proof} Assume that extremal rays $Q_1, Q_2$ of the type (II) are
joint by the single arrow $Q_1Q_2$, and extremal
rays $R_1,R_2$ of the type (II) are joint by the single arrow
$R_1R_2$. Let us suppose that the sets $\{ Q_1, Q_2\}$ and
$\{ R_1, R_2\}$ are divisorially joint. We consider possible
cases below.

Suppose that $Q_1,R_1$ are divisorially disjoint.

If $Q_2=R_2$, we get the configuration $\cc_3$.
We prove that the only case is possible.

Assume that $Q_2\not=R_2$.
Since the curve $D(R_1)\cap D(R_2)$ belongs to $R_2$, we have
$R_2\cdot D(Q_1)=0$ (because $D(R_1) \cap D(Q_1)=\emptyset$).
Similarly, $Q_2\cdot D(R_1)=0$.

Let us suppose that $D(Q_2)=D(R_2)=D$.
Then any curve of the divisor $D$  belongs to
$R_2+Q_2$ by Lemma 1.2.4. Since $Q_2\cdot D(R_1)=R_2\cdot D(R_1)=0$,
divisors $D(R_2)$ and $D(R_1)$ do not have a common curve.
We get a contradiction, because, $R_1R_2$ is an arrow.

Thus, we can suppose that $D(Q_2)\not= D(R_2)$.

Let us suppose that $R_1Q_2$ is an arrow. We have proved that
$Q_2R_1$ is not an arrow. It follows that
the divisor $D(R_1)$ contains curves of 3 extremal rays:
$R_1, R_2, Q_2$ which is impossible. Thus, $R_1Q_2$ and
$Q_2R_1$ are not arrows. It follows that divisors
$D(R_1)$ and $D(Q_2)$ are disjoint. Similarly, the divisors
$D(R_2)$ and $D(Q_1)$ are disjoint.
Since the curve $D(R_1)\cap D(R_2)$ belongs to $R_2$, we get that
$R_2\cdot D(Q_2)=0$. Similarly, $Q_2\cdot D(R_2)=0$.
It follows that divisors $D(R_2)$ and $D(Q_2)$ are
divisorially disjoint.
It follows that sets $\{ R_1, R_2\}$ and $\{ Q_1, Q_2\}$
are divisorially disjoint. We get a contradiction.

Thus, we proved that if the rays $Q_1$ and $R_1$ are divisorially
disjoint, then $Q_2=R_2$, and we have the configuration $\cc_3$.

Now, assume that $Q_1$ and $R_1$ are divisorially joint.
Then, $D(Q_1)\cap D(R_1)$ is non-empty.

If $D(Q_1)=D(R_1)$, evidently, sets $\{Q_1,Q_2\}$ and
$\{ R_1,R_2\}$ are equal because the divisor $D(Q_1)=D(R_1)$
cannot contain $3$ extremal rays by Lemma 1.2.3.

Thus, we assume that $D(Q_1)\cap D(R_1)$ is a non-empty curve.
By Lemma 2.1.2,
$\nem (X,$ $D(Q_1)) = Q_1+Q_2$ and $\nem (X, D(R_1)) = R_1+R_2$
are 2-dimensional faces of $\nem (X)$. It follows that faces
$Q_1+Q_2$ and $R_1+R_2$ have a common extremal ray.

Thus,
sets $\{Q_1, Q_2\}$ and $\{R_1,R_2\}$ have a common extremal ray.

Assume that $Q_2=R_1$ (the case $Q_1=R_2$ is similar). We denote
$Q_3=R_2$. Thus, we have 3 extremal rays $Q_1, Q_2, Q_3$ where
$Q_1,Q_2$ are joint by a single arrow $Q_1Q_2$ and
$Q_2, Q_3$ are joint by a single arrow $Q_2Q_3$.
The curve $D(Q_1)\cap D(Q_2)$ belongs to the ray $Q_2$ and
$Q_2\cdot D(Q_3)>0$. it follows that divisors $D(Q_1)$ and
$D(Q_3)$ have a non-empty common curve.
The extremal rays $Q_1$ and $Q_3$ cannot be joint by a single arrow
$Q_1Q_3$ since then $D(Q_1)$ contains curves of 3 different extremal rays
$Q_1,Q_2,Q_3$. Thus, the extremal rays $Q_1,Q_3$ are joint either by
non-single arrows $Q_1Q_3$ and $Q_3Q_1$ or by the single arrow $Q_3Q_1$.
Thus, we get cases $\Tt_3$ and $\Tt_3^\prime$.

Now suppose that $Q_2=R_2$. We denote $R_1=Q_3$. Thus, we have
3 extremal rays $Q_1,Q_2,Q_3$ such that $Q_1,Q_2$ are joint by a
single arrow $Q_1Q_2$, and $Q_3,Q_2$ are joint by a single arrow
$Q_3Q_2$. If $Q_1$ and $Q_3$ are divisorially disjoint, we get the
case $\cc_3$.

Suppose that $D(Q_1)\cap D(Q_3)$ is a non-empty curve.
By Lemma 2.1.2,
$\nem (X, D(Q_1)) = Q_1+Q_2$ and $\nem (X, D(Q_3)) = Q_3+Q_2$
are 2-dimensional faces of $\nem (X)$. The set of common points of
this faces is the extremal ray $Q_2$. It follows that the curve
$D(Q_1)\cap D(Q_3)$ belongs to the ray $Q_2$.  By Lemma 2.1.3, this
is impossible.

Thus, we had proved that divisors $D(Q_1)$ and $D(Q_3)$ do not have a common
point.

As a result, we have proved that any configuration of two single arrows
which are divisorially joint
satisfies the statement of Theorem. Thus this is
either $\cc_3$ or $\Tt_3$ or $\Tt_3^\prime$.

Suppose that there exists a single arrow $R_1R_2$ which is divisorially
joint with a configuration $\Tt_3$ or $\Tt_3^\prime$
of extremal rays and different from
arrows of this configurations. By the result above, the set
$\{R_1,R_2\}$ has a common extremal ray with both sets $\{Q_1,Q_2\}$
and $\{Q_2,Q_3\}$. Thus, either $Q_2=R_1$ or $Q_2=R_2$. If
$Q_2=R_1$, the divisor $D(R_1)$ contains 3 extremal rays: $R_1, R_2, Q_3$
which is impossible. If $Q_2=R_2$, then $R_1,R_2,Q_3$ give a configuration
$\Tt_3$ or $\Tt_3^\prime$. Then the angle
$\nem (X, D(Q_3))$ has a side
which simultaneously is a ray of the $Q_1+Q_2$ and
$R_1+Q_2$. Thus, this side is $Q_2$. Then $D(Q_2)=D(Q_3)$
which is impossible.

 From above, it follows that a divisorially connected component $S$ of
a finite non-empty configuration of single arrows of extremal rays of
the type (II) has the type
either $\cc_n, n\ge2$ or $\Tt_3$ of $\Tt_3^\prime$.
Using Lemma 2.1.2, one can easily see that
if an extremal ray $Q$ of this connected component
$S$ belongs to a pair $Q,R$ of the type $\bb_2$, then $S$ has the type
$\cc_n, n\ge2$ and $Q=Q_1$ in notation of Theorem.
Thus, $S\cup \{R\}$ gives rise a configuration of the type $\bb_2\cc_n$.

Now suppose that we have two sets $\{Q_1,Q_2\}$
and $\{ R_{11},R_{12}\}$ of extremal
rays such that $Q_1Q_2$ is a single arrow
and $R_{11}R_{12}$ is a dotted line, and this two sets are divisorially
joint. If these two sets have a common extremal ray, then they give
a configuration $\bb_2\cc_2$ since above considerations.

Thus, let us assume that these two sets do not have a common extremal
ray. By Lemmas 2.1.2 and 1.2.6,
$\nem (X, D(Q_1))=Q_1+Q_2$ and $\nem (X, D(R_{11}))=R_{11}+R_{12}$
are 2-dimensional faces of $\nem (X)$. It follows that if
divisors $D(Q_1), D(R_{11})$ have a common point, then
sets $\{Q_1,Q_2\}$
and $\{ R_{11},R_{12}\}$ have a common extremal ray.
But we assume that this is not the case. Thus, we can suppose that
divisors $D(Q_1)$ and $D(R_{11})$ do not have a common point.
The curve $D(Q_1)\cap D(Q_2)$ belongs to the ray $Q_2$. Since
$D(Q_1)\cap D(R_{11})=\emptyset$, then
$Q_2\cdot D(R_{11})=0$. Thus, $R_{11}Q_2$ is a single arrow.
This is impossible by Lemma 2.1.2. Thus, we get a contradiction.

Thus, we proved that a finite configuration of
single arrows and dotted lines has divisorially
connected components of the type $\bb_2,
\cc_n, \bb_2\cc_n, \Tt_3, \Tt_3^\prime$ of Theorem.

Now let $R$ be an extremal ray of the type (I). Let $Q_1Q_2$ be
a single arrow. If $D(Q_1)\cap D(R)$ is non-empty, then
$D(Q_1)$ contains 3 extremal rays $Q_1,Q_2$ and $R$ which is impossible.
The curve $D(Q_1)\cap D(Q_2)$ belongs to $Q_2$. From above,
$Q_2\cdot D(R)=0$. Thus, $RQ_2$ is a single arrow which is impossible for
an extremal ray $R$ of the type (I). Thus, $\{ R\}$ and $\{ Q_1, Q_2\}$
are divisorially disjoint.

By Lemma 1.2.5 and considerations above, we proved that any finite
set of special extremal rays has connected components of the
Theorem. Let us suppose that this set $S$ has $n$ extremal rays.
We claim that then $n\le 2\rho (X)$ where $\rho (X)=\dim N_1(X)$.
 From the description of connected components of $S$, for any connected
component $S_i$ of $S$ with $m_i$ elements there exists at least
$\max \{1, m_i-1\}\ge m_i/2$
extremal rays with disjoint divisors. Since any
extremal ray $R\in S$ has $R\cdot D(R)<0$,
it follows the inequality above.

All other statements of Theorem follow from Lemmas 1.2.6 and
2.1.2.

This finishes the proof of Theorem.

\enddemo

Now let us consider a divisorial extremal ray $S$ which is not a special one
but the divisor $D(S)$ has a common point with the divisor of one of special
divisorial extremal rays. Thus, $S$ is joint by non-single arrows with a
special divisorial extremal ray.
These extremal rays also look very specially. We
describe them below.

\proclaim{Theorem 2.1.8} Let $X$ be a $3$-fold from the class $\LT$.
Let $S$ be a non-special divisorial extremal
ray of $X$.
Let $Z_S$ be the set of
special divisorial extremal rays which are divisorially
joint with $S$ (by definition, $S$ has the type (II) and
is joint by non-single arrows with all extremal rays from $Z_S$).
Assume that $Z_S\not=\emptyset$.

Then $Z_S$ is contained in one of divisorially
connected components $C=C_S$ of the set of special extremal rays, and
depending to their types of Theorem 2.1.7, we have the following
description for
$Z_S\subset C_S$
:

($\Aa_1(S)$) $C_S$ has the type $\Aa_1$, $Z_S=Q$;

($\cc_n^{(1)}(S)$)  $C_S$ has the type $\cc_n$, $Z_S=\{Q_i\}$,
for one of $i>1$;

($\cc_2^{(2)}(S)$) $C_S$ has the type $\cc_2$, $Z_S=\{ Q_1, Q_2\}$;

($\bb_2\cc_n(S)$) $C_S$ has the type $\bb_2\cc_n$,
$n>1$, $Z_S=\{Q_i\}$ for one of $i>1$;

($\bb_2(S)$) $C_S$ has the type $\bb_2$, $Z_S=\{Q_{11}, Q_{12}\}$.

We remark that it follows that $C_S$ cannot have the type $\Tt_3$
and $\Tt_3^\prime$.

\endproclaim


\demo{Proof} We consider one case only.
Similarly, one can consider all other cases.

Let us assume that $S$ is divisorially joint with a divisorially
connected component of the type $\cc_n,\ n\ge 2$, of the set of
special extremal rays.  Thus, $D(S)$ has a common
curve with one of divisors $D(Q_i),\ i=1,...,n$
(we use notation of Theorem 2.1.7).

Let us suppose that $D(S)$ has a common curve with the divisor $D(Q_1)$.
Then $S$ and $Q_1$ are joint by non-single arrows. Thus, $Q_1\cdot D(S)>0$.
Since $D(Q_i)\cap D(Q_1)$ belongs to $Q_1$,  we get that
$D(S)$ has a common curve with all  divisors $D(Q_i)$, $i=2,...,n$, too.
By Lemma 2.1.2, $\nem (X,D(Q_i))=Q_1+Q_i$ is
a 2-dimensional face of $\nem (X)$.
We then get that the 2-dimensional angle $\nem (X,D(S))$ has
the second side (different from $S$) which belongs to the
face $Q_1+Q_i$, $i=2,...,n$. If $n>2$, we then get that
the second side of $\nem (X, D(S))$ is generated by the
curve $D(S)\cap D(Q_i)$ which belongs to $Q_1$.
This is
impossible by Lemma 2.1.3. Thus, $n=2$, and we get the case
$\cc_2^{(2)}(S)$.

Now, assume that $S$ is not divisorially joint with $Q_1$,
thus $D(S)\cap D(Q_1)=\emptyset$.
Then $D(S)\cap D(Q_i)$ is not empty for one of $i>1$. Then, like
above, $\nem (X, D(S))$ is a 2-dimensional angle with the second side
generated by the curve $D(S)\cap D(Q_i)$ which
belongs to the 2-dimensional face
$\nem (X, D(Q_i))=Q_i+Q_1$ of $\nem (X)$. This curve cannot belong to
$Q_1$ since $D(S)\cap D(Q_1)=\emptyset$.
Thus, the second side of the angle $\nem (X, D(S))$ is the ray of the
$Q_i+Q_1$ different from $Q_1$. It follows that the extremal ray $Q_i$,
$i=2,...,n$, such that $D(S)\cap D(Q_i)\not= \emptyset$ is unique.
Thus, we get the case $\cc_n^{(1)}(S)$ if we additionally prove that
$S$ is divisorially disjoint with all special extremal rays different
from $Q_1,...,Q_n$.

By Theorem 2.1.7, any special extremal ray $Q$ of the type (II)
together with another special extremal ray $Q^\prime$ of the same
component of the set of special extremal rays generate a 2-dimensional
face  $Q+Q^\prime$ of $\nem (X)$, and this 2-dimensional face is either
$\nem (X,D(Q))$ or $\nem (X,D(Q^\prime))$ (for the extremal ray
$Q$ of the type (I), one should
use that $\nem (X, D(Q))=Q$). Using this property, like
above, one can see that $D(S)\cap D(Q)$ is empty for
any special extremal ray $Q$ different from extremal rays
$Q_1,...,Q_n$ above.
   One can similarly consider cases when $S$ is divisorially
joint with other types of connected components
of the set of special extremal rays.
\enddemo

Now, we have the following description for cases when $S\not= S^\prime$
but $Z_S\cap Z_{S^\prime}\not=\emptyset$. These cases are very rare.

\proclaim{Theorem 2.1.9}
Let $X$ be a $3$-fold from the class $\LT$.
Let $S, S^\prime$ are non-special divisorial
extremal rays such that $Z_S\cap Z_{S^\prime}\not=\emptyset$.

Then $S=S^\prime $ except cases below (in notation of Theorems
2.1.7 and 2.1.8):

(1) Case $C_S=C_{S^\prime }=\cc_2,\ Z_S=Z_{S^\prime}=\{ Q_1,Q_2\}$,
$S,S^\prime$ are joint by non-single arrows,
$\dim~[Q_1,Q_2,S,S^\prime ]=3$.

(2) Case $C_S=C_{S^\prime }=\cc_2,\ Z_S=\{ Q_1,Q_2\},
\ Z_{S^\prime}=\{ Q_2 \}$,
$S,S^\prime$ are joint by non-single arrows,
$\dim~[Q_1,Q_2,S,S^\prime ]=3$.

(3) Case $C_S=C_{S^\prime}=\bb_2,\ Z_S=Z_{S^\prime}=\{Q_{11}, Q_{12}\}$,
$S, S^\prime$ are joint by non-single arrows,
$\dim~[Q_{11}, Q_{12},S, S^\prime]=3$

(4) Case $C_S=C_{S^\prime}=\bb_2,\ Z_S=Z_{S^\prime}=\{Q_{11}, Q_{12}\}$,
$S, S^\prime$ are joint by non-single arrows,
there exists a ray $T\subset Q_{11}+Q_{12}$ such that
$\nem (X,D(S))=T + S$ and $\nem (X,D(S^\prime ))=T + S^\prime $.
\endproclaim

\demo{Proof} Let $S^\prime \not= S$ but
$Z_{S^\prime } \cap Z_S \not= \emptyset$.
Using Theorem 2.1.8, we consider all possible cases
(we use notation of Theorems 2.1.8 and 2.1.7).

The case ($\Aa_1(S)$), $Z_S=Z_{S^\prime}=\{ Q\}$.
First, remark that
$D(S)\cap D(S^\prime)\not=\emptyset$ since $D(S), D(S^\prime)$
both have non-empty intersection with $D(Q)$ and
$Q\cdot D(S)>0$ and $Q\cdot D(S^\prime )>0$. Besides,
$\nem (X,D(S))=S+Q$ and $\nem (X,D(S^\prime))=S^\prime +Q$
are $2$-dimensional angles with the common edge $Q$. If the
curve $D(S)\cap D(S^\prime )$ has a component which
does not belong to $Q$,
then one of this angles is contained in another one which is
impossible since $S$ and $S^\prime$ are different extremal rays.
Thus, the curve $D(S)\cap D(S^\prime )$ belongs to the extremal ray $Q$.
This is impossible by Lemma 2.1.3.

The case ($\cc_n^{(1)}(S)$),
$C_S$ has the type $\cc_n$, $Z_S=Z_{S^\prime} = \{Q_i\}$, for one of $i>1$.

Let $C=D(S)\cap D(Q_i)$. We have $C\cdot D(Q_1)=0$ because
$D(S)\cap D(Q_1)=\emptyset$. Besides,
$C\in Q_i+Q_1$ and $Q_i\cdot D(Q_1)>0,\ Q_1\cdot D(Q_1)<0$.
Thus, the ray $\r^+C\subset Q_i+Q_1$ is defined uniquely by
the property $\r^+C\cdot D(Q_1)=0$.

Now assume that $S\not=S^\prime$ and $Z_{S^\prime}=Z_S=\{ Q_i\}$.
Like above, $D(S^\prime )\cap D(Q_i)=C^\prime$ and
$\r^+C^\prime=\r^+C$. By Lemma 2.1.2,
$\nem (X,D(Q_i))=Q_i+Q_1$ is a
$2$-dimensional face of $\nem (X)$. It follows that
$\nem (X,D(S))=\r^+C + S$ and $\nem (X,D(S^\prime ))=\r^+C + S^\prime$.
If the curve $D(S)\cap D(S^\prime)$ has a component which does not
belong to $D(Q_i)$, we then get that one of angles
$\nem (X,D(S)),\ \nem (X,D(S^\prime))$ contains another, which is
impossible since $S, S^\prime$ are different extremal rays of $\nem (X)$.
Let us consider
the contraction $f:X\to X^\prime$ of the extremal ray $Q_i$.
Then the curve $C=f(D(S)) \cap f(D(S^\prime))$ belongs to the extremal
ray $f(Q_1)$ of $X^\prime$. This is impossible by Lemma 2.1.3.

The case ($\cc_2^{(2)}(S)$), $C_S$ has the type
$\cc_2$, $Z_S=\{ Q_1, Q_2\}$.
For this case, $D(S)\cap D(Q_1)$ contains a component which does not
belong to $Q_1$ (otherwise, $SQ_1$ is a single arrow).
It follows that the $2$-dimensional angles $\nem (X,D(S))$ and
$\nem (X,D(Q_1))$ have a common ray different from the edge $Q_1$ of
the $\nem (X,D(Q_1))$.
Besides, the curve $D(S)\cap D(Q_2)$ is not empty (since $SQ_2$ is arrow).
Thus, the $2$-dimensional
angle $\nem (X,D(S))$ has a common ray with the angle
$\nem (X,D(Q_2))$. It follows that the angle $\nem (X,D(S))$ is contained in
the $3$-dimensional space $V$ generated by angles $\nem (X,D(Q_1))$ and
$\nem (X,D(Q_2))$ with the common edge $Q_1$.
The same is valid for any $S^\prime$ with
$Z_{S^\prime}=Z_S=\{ Q_1,Q_2\}$. This gives the case (1) of
the Theorem.

By Theorem 2.3.8, $Z_{S^\prime}=\{ Q_2\}$ if
$Z_{S^\prime}\cap Z_S\not=\emptyset$ and
$Z_{S^\prime}\not=Z_S$.
Let $C^\prime =D(S^\prime)\cap D(Q_2)$. Since $\nem (X,D(Q_2))=Q_2+Q_1$,
$C^\prime \in Q_2+Q_1$. Since $Q_2\cdot D(S)>0, Q_1\cdot D(S)>0$, it then
follows that $C^\prime \cdot D(S)>0$ and
$D(S^\prime)\cap D(S)\not=\emptyset$.
Thus, angles $\nem (X,D(S))$ and $\nem (X,D(S^\prime))$ have a common
ray. If $\nem (X,D(S)) \cap \nem (X,D(S^\prime))=\r^+C^\prime$,
like above, considering the contraction of
the extremal ray $Q_2$, we get the contradiction with Lemma 2.1.3.
Thus, $\nem (X,D(S))\cap \nem (X,D(S^\prime))$ contains a
ray which is different from $\r^+C^\prime$.
It follows that the angle
$\nem (X,D(S^\prime))$ is contained in the
$3$-dimensional space $V$ above
containing extremal rays $Q_1,Q_2,S$. This gives the case (2) of
the Theorem.

The case ($\bb_2\cc_n(S)$), $C_S$ has the type $\bb_2\cc_n$,
$n>1$, $Z_S=Z_{S^\prime}=\{ Q_i\}$ for one of $i>1$. Like
above for the case $\cc_n$, this case is impossible for
$S\not=S^\prime$.

The case ($\bb_2(S)$), $C_S$ has the type $\bb_2$,
$Z_S=\{ Q_{11},Q_{12}\}$: By Lemma 1.2.4,
$\nem (X,$ $D(Q_{11}))=Q_{11}+Q_{12}$ is a face of $\nem (X)$.
It follows that $\nem (X,D(S))=T+S$ where $T$ is a ray of the angle
$Q_{11}+Q_{12}$. Let $S^\prime $ be another extremal ray with
$Z_{S^\prime}=\{ Q_{11},Q_{12}\}$. Since $Q_{11}\cdot D(S^\prime)>0,
Q_{12}\cdot D(S^\prime)>0$, we have $T\cdot D(S^\prime)>0$.
It follows that $S,S^\prime$ are connected by non-single arrows.
Applying Lemma 2.1.5, we get cases (3) and (4) of the Theorem.
This finishes the proof.
\enddemo

\subhead
2.2. Extremal and $E$-sets of divisorial extremal rays
\endsubhead

We want to apply results of Section 2.1 to describe
extremal and $E$-sets of extremal rays of the type
(I) and (II) for $3$-folds of the class $\LT$. But it is
important for future studies considering more general
elliptic, connected parabolic and Lanner subsets of
divisorial extremal rays which had in fact appeared in
Section 1.2.

\definition{Definition 2.2.1}
Let  $\E=\{R_1,...,R_n\}$ be a set of extremal rays of the type (I) or (II).
Then $\E$ is called {\it elliptic} if there are
$a_1>0,...,a_n>0$ such that
$$
R_i\cdot (a_1D(R_1)+\cdots +a_nD(R_n))<0
$$
for all $1\le i \le n$.
\enddefinition

\definition{Definition 2.2.2}
Let  ${\Cal P}=\{R_1,...,R_n\}$ be a set of
extremal rays of the type (I) or (II).
Then $\Cal P$ is called
{\it connected parabolic} if each proper subset of $\Cal P$ is
elliptic and there are
$a_1>0,...,a_n>0$ such that
$$
R_i\cdot (a_1D(R_1)+\cdots +a_nD(R_n))=0
$$
for all $1\le i \le n$.

A set $Q$ of extremal rays of the type (I) or
(II) is called {\it parabolic}
if each divisorially connected component of $Q$
is connected parabolic.
\enddefinition

\definition{Definition 2.2.3}
Let  $\La=\{ R_1,...,R_n\}$ be a set of
extremal rays of the type (I) or (II).
Then $\La$ is called {\it Lanner} (respectively, {\it quasi-Lanner})
if each proper subset of $\Cal L$ is
elliptic (respectively, either elliptic or parabolic)
and there are
$a_1>0,...,a_n>0$ such that
$$
R_i\cdot (a_1D(R_1)+\cdots +a_nD(R_n))\ge 0
$$
for all $1\le i \le n$ and there exists $j,\ 1\le j \le n$, such that
$$
R_j\cdot (a_1D(R_1)+\cdots a_nD(R_n)) >0.
$$
\enddefinition

The following statement will be useful.

\proclaim{Proposition 2.2.4}
Let  $\E=\{R_1,...,R_n\}$ be a set of extremal rays of the type (I) or (II)
with different divisors $D(R_1),...,D(R_n)$.
Then $\E$ is elliptic if and only if for any non-negative
$b_1,...,b_n$ not all equal to $0$ there exists $R_i,\ 1\le i \le n,$
such that
$$
R_i\cdot (a_1D(R_1)+\cdots +a_nD(R_n))<0.
$$
\endproclaim

\demo{Proof} If $\E$ is elliptic, then $\E$ satisfies the condition above
by Lemma 1.2.10.

Now assume that $\E$ satisfies the condition of Proposition.
If $\E$ is not elliptic, there exists a minimal non-elliptic
subset $\E^\prime \subset \E$. The subset $\E^\prime$ is not
empty since any one element subset of $\E$ is elliptic.
Let $\E^\prime =\{R_1,...,R_t\}$
where $1 \le t\le n$. By Lemma 1.2.10, there are positive
$c_1,...,c_t$ such that
$$
R_i\cdot (c_1D(R_1)+ \cdots + c_tD(R_t))\ge 0
$$
for any $1\le i \le t$.
The same inequality is evidently true for $t+1 \le i \le n$.
We get a contradiction with the condition of Proposition.
\enddemo

\definition{Definition 2.2.5} A sequence
$$
C=\{ R_1,...,R_n\},
$$
$n\ge 1,$ of extremal rays of the type (I) or (II) is called
{\it chain} if
$R_i, R_j$ are divisorially disjoint for $j-i>1$, and
$R_iR_{i+1}$ is a non-single arrow (thus, $R_{i+1}R_i$ is
an arrow too) for any $1\le i<n$.
Extremal rays $R_1, R_n$ are called terminal for the chain.
\enddefinition

\proclaim{Theorem 2.2.6} Let $X$ be a $3$-fold from the class
$\LT$. Let $T$ be a divisorially connected
set of extremal rays of the type
(I) or (II) such that any proper subset of $T$
is elliptic. (In particular, this is valid if either $T$ itself
is elliptic or is contained in a face of $\nem (X)$ of Kodaira dimension $3$
or $T$ is connected parabolic or Lanner or $T$ is an
$E$-set such that each proper subset of $T$
is contained in a face of $\nem (X)$ of Kodaira dimension $3$).
Let us assume that each two extremal rays of $T$ are different from
a pair of the type $\bb_2$. Then $T$ has one of the types
:

(A)  All extremal rays of $T$ have the type (II) and $T$ does not
have a single arrow.

(B) $T=\{ R\}\cup C_1\cup ... \cup C_k$
has only extremal rays of the type (II), where
$$
C_1=\{ R_{11},...,R_{1n_1}\},\  C_2=\{ R_{21},...,R_{2n_2}\},...,
C_k=\{ R_{k1},...,R_{kn_k}\}
$$
are divisorially disjoint to one another chains, and all arrows between
$R$ and extremal rays of these
chains are single arrows $R_{j1}R,\ j=1,...,k$; besides
each extremal ray
of the chains $C_1,...,C_k$ does not belong to a pair of the type $\bb_2$.

(C) $T=\{R_1,R_2,R_3,...,R_n\}, \ n\ge 3,$ where
all extremal rays of $T$ have the type (II) and $R_2R_1$ is a single
arrow, $R_1R_2, R_2R_1$ and $R_2R_3,R_3R_2$ are non-single arrows,
$R_3,...,R_n$ is a chain such that $R_4,...,R_n$ are
divisorially disjoint with extremal rays $R_1$ and $R_2$; besides,
each extremal ray $R_1,R_2, R_3,...,R_n$ does not belong
to a pair of the type $\bb_2$.

(D) $T=\{ R_1,R_2,...,R_i,...,R_n\}$ is a chain such that $R_1$
has the type (I), and $R_i$ has the type (II) for all $i>1$;
besides, each extremal ray $R_1,R_2,...,R_n$ does not belong to a pair
of the type $\bb_2$.

(E) $T=\{Q_1,Q_2,Q_3\}$ has type of triangle: thus,
$Q_1Q_2$ and $Q_2Q_3$ are single arrows, and $Q_1Q_3, Q_3Q_1$ are
non-single arrows; besides, each extremal ray $Q_1,Q_2,Q_3$ does not
belong to a pair of the type $\bb_2$.

(E') $T=\{ Q_1,Q_2,Q_3\}$ has type of special triangle: thus,
$Q_1Q_2,\ Q_2Q_3,\ Q_3Q_1$ are single arrows; besides, each extremal
ray $Q_1,Q_2,Q_3$ does not belong to a pair of the type $\bb_2$.

Moreover, for the case (B), the $T$ is contained in
a face of $\nem (X)$ of Kodaira dimension $3$. Thus, the case (B)
is impossible when $T$ is not contained in a face of $\nem (X)$
of Kodaira dimension $3$.
\endproclaim



\demo{Proof} If $T$ does not have an extremal ray of
the type (I) and a single arrow, we get the case (A).

Thus, we can suppose that $T$ has either an extremal ray of the
type (I) or a single arrow.

Let $T_1$ be the set which contains
all extremal rays of the type (I) and all extremal rays of single
arrows of $T$. If $T=T_1$, then $T_1$ is divisorially connected, and
by Theorem 2.1.7, we get one of following cases: (B) with
one element chains, (D) with one extremal ray $R_1$,
(E) or (E').

Let us suppose that $T\not=T_1$. By Theorem 2.1.7,
$T_1$ has divisorially connected components of
the types $\Aa_1$, $\cc_n$, $\Tt_3$ or
$\Tt^\prime_3$.
Since $T$ is divisorially connected, for each connected component of
$T_1$ there exists an extremal rays $Q\in T-T_1$
such that $Q$ is joined by
non-single arrows with an extremal ray
of this connected component. By Theorem 2.1.8, then
this connected component should have the type $\Aa_1$ or $\cc_n$.
It follows that all connected components of $T_1$ have the type
$\Aa_1$ or $\cc_n$, and $T_1$ has a connected component
of one of these types.

We consider induction on $\# T$. For $\# T=1, 2$ the statement is
clear. Let us suppose that $\#T =n>2$. We consider several cases.

Let us suppose that there exists an extremal ray $R\in T$ of the
type (I). Let $Q\in T$ and $Q\not= R$. We claim that then
$Q+R$ is a $2$-dimensional face of Kodaira dimension $3$ of $\nem (X)$.
At first, suppose that  $D(R)\cap D(Q)=\emptyset$.
Let $H$ be a $nef$ element orthogonal to $Q$ (the set of these
elements $H$ defines a face of codimension one in $NEF(X)$
by the exact sequence (1-2-1)).
There exists a linear function $\alpha(H) \ge 0$ such that the
map  $H\to H^\prime =H + \alpha (H)D(R)$ is
linear with one-dimensional kernel,
and $H^\prime$ is orthogonal to extremal rays $Q$ and $R$.
The $H^\prime $ is evidently $nef$ since $R$ has the type (I) and $C\in R$
for any curve $C\subset D(R)$.
Thus, $Q+R$ is a $2$-dimensional
face of $\nem (X)$, and, by construction, the face $(Q+R)^\perp$ of
$NEF(X)$ has codimension $2$. Evidently, $(H^\prime)^3 \ge H^3>0$
(see the proof of Lemma 1.2.6). Thus, by definition of the class $\LT$,
this face has Kodaira dimension $3$. (Further, we will not be so
formal, and just show that the corresponding face $\alpha$ of
$\nem (X)$ has numerical
Kodaira dimension $3$; automatically, by our construction, it will have
the property $\codim \alpha^\perp=\dim \alpha$ and will have
Kodaira dimension $3$ by definition of the class $\LT$).
Now, suppose that $D(Q)\cap D(R)\not=\emptyset$. Since $\# T =n>2$,
the set $\{ Q, R\}$ is elliptic.
Thus, there are positive $a_1, a_2$ such that
$Q\cdot (a_1D(Q)+a_2D(R))<0,\ R\cdot (a_1D(Q)+a_2D(R))<0$.
Let $H$ be an ample element. By Lemma 2.1.8, there are
positive $c_1,c_2$ such that
$H^\prime = H + c_1D(Q) + c_2D(R)$ is orthogonal to both $Q, R$.
Since $R$ has the type (I) and $D(Q)\cap D(R)\not=\emptyset$,
we then get $\nem (X,D(Q))=Q+R$. It follows that
$H^\prime \cdot C\ge 0$ for any curve $C$. Thus, $H^\prime$ is
$nef$. Evidently, $H^\prime $ is orthogonal to extremal rays
$Q, R$ only and $H^\prime \ge H^3>0$.
It follows that $Q+R$ is a face of dimension $2$
and Kodaira dimension $3$ of $\nem (X)$.
Thus, we proved that $Q+R$ is a $2$-dimensional face of $\nem (X)$
for any $Q\in T$ and $Q\not= R$.

Let us consider the contraction
$f:X\to X^\prime$ of the extremal ray $R$. Then the image $f(Q)$ is
an extremal ray of $X^\prime$ for $Q\not=R$ since the
claim we proved above.
Evidently, if $D(Q)\cap D(R)=\emptyset$, the extremal ray
$f(Q)$ has the same type
((I) or (II)) as $Q$. If $D(Q)\cap D(R)\not= \emptyset$, the
$f(Q)$ has the type (I). By Theorem 2.1.9, the
extremal ray $Q$ with this property is unique. Since $T$ is
connected, the $Q$ does exist.
By the projection formula and Proposition 2.2.4,
any proper subset of the set
$f(T)=f(T-\{ R\})$ is elliptic. It follows that
the set $f(T)$ has the same properties as $T$ but has one
element less. By induction, we then get that
$f(T)$ has the type (D).
Considering preimages, one can easily see that then
$T$ has the same type (D) and has desirable properties.

Now, assume  that $R_2R_1$ is a single arrow
for $R_2,R_1 \in T$.
At first, let us suppose that there does not exist $Q\in T$ such
that $Q$ is joined by non-single arrows with the terminal $R_1$
of the arrow $R_2R_1$. Thus,
$R_1\cdot D(Q)=0$ for all $Q\in T,\ Q\not=R_1$.
By Lemma 2.1.2,
$\nem (X,D(R_2))=R_1+R_2$ is a $2$-dimensional face of $\nem (X)$ of
Kodaira dimension $3$.
It follows that $\nem (X,D(Q))=P+Q,\ P\subset R_1+R_2$
for an extremal ray $Q$ which is
joined by non-single arrows with $R_2$. By Lemma 2.1.9, this extremal
ray $Q$ is unique if it does exist.
Then, like above, we can prove that $R_2, Q$ generate a $2$-dimensional
face  $R_2+Q$ of Kodaira dimension $3$ of $\nem (X)$
for each $Q\in T,\  Q\not=R_2$ (for $Q=R_1$, this follows from
Lemma 2.1.2).
Let us consider the contraction $f:X \to X^\prime $
of the extremal ray $R_2$. Then $f(Q)$ is a divisorial
extremal ray of $X^\prime$ for $Q\in T,\ Q\not=R_2$.
If there exists $Q$ such that $QR_2$ is a non-single arrow, then
$f(Q)f(R_1)$ is evidently a single arrow because $Q$ and $R_1$
are divisorially disjoint by Lemma 2.1.8.
If there does not exist $Q$ such that $QR_2$ is a non-single
arrow then
there exists an extremal ray $S\in T$ such that
$SR_1$ is a single arrow and $S$ is divisorially disjoint with $R_2$.
This follows from Lemmas 2.1.7, 2.1.8 and conditions on $T$ above
(otherwise, $T=\{ R_2, R_1\}$). Then $f(S)f(R_1)$ is a
single arrow. Thus, the image $f(T)$ is a
divisorially connected set of extremal rays of the type (I) and (II)
which contains a single arrow and like for the case above
any proper subset of $f(T)$ is
elliptic. By induction, the $f(T)$ has the type
(B) and only the extremal ray $f(R_1)$ may belong to a pair of the
type $\bb_2$.
One can see easily using our construction, that then $T$ has
also the type (B) and only $R_1$ may belong to a pair of the type
$\bb_2$.

Now we consider the case when $R_2R_1$ is a single arrow for
$R_2,R_1\in T$ and there exists $R_3\in T$ such that
$R_3R_1,R_1R_3$ are non-single arrows. This is the case
$\cc_2^{(2)}(R_3)$ of Theorem 2.1.8. By Theorem 2.1.8, then
$R_3R_2,R_2R_3$ are non-single arrows and $R_2R_1$ is the only
single arrow in $T$ with the terminal $R_1$.

Since $\# T \ge 3$ and any proper subset of $T$ is elliptic,
there are
$a_2>0,a_3>0$ such that $R_2\cdot (a_2D(R_2)+a_3D(R_3))<0$ and
$R_3\cdot (a_2D(R_2)+a_3D(R_3))<0$. Let $H$ be an ample element.
By Lemma 1.2.8, there exist $b_2>0,\ b_3>0$ such that
$H^\prime =H+b_2D(R_2)+b_3D(R_3)$ is orthogonal to both $R_2,R_3$.
By Lemma 2.1.2, $\nem (X, D(R_2))=R_1+R_2$. It follows that
$C\cdot H^\prime \ge 0$ for any curve $C\subset D(R_2)$.
We have $D(R_2)\cap D(R_3)$ is a non-empty curve.
By Lemma 2.1.2 $\nem (X,D(R_2))=R_1+R_2$ is a face of $\nem (X)$.
It follows that $\nem (X,D(R_3))=S+R_3$ where $S$ is a ray of
$R_1+R_2$. If follows that $C\cdot H^\prime \ge 0$ for any curve
$C\in D(R_3)$. This implies that $H^\prime$ is $nef$ and only
extremal rays $R_2,R_3$ are orthogonal to $H^\prime$. Besides,
like above, $(H^\prime )^3\ge H^3>0$. Thus, $R_2+R_3$ is a
2-dimensional face
of $\nem (X)$ of Kodaira dimension $3$.

As we have mentioned above, $\nem (X,D(R_3))=S+R_3$, where
$S\in R_1+R_2$. Here $S\not=R_1$ since otherwise,
$R_3R_1$ is a single arrow. We have $D(R_3)\cap D(R_1)\not=\emptyset$.
Since $R_2+R_3$ is a face of $\nem (X)$, it follows that
$\nem (X,D(R_1))\subset R_1+R_2+R_3$. Then, like above, we
can see that $R_1+R_3$ is a $2$-dimensional face of $\nem (X)$ of
Kodaira dimension $3$. By Theorems 2.1.7 and 2.1.8, each extremal
ray $R_1,R_2,R_3$ does not belong to a pair of the type $\bb_2$.

If $T=\{ R_1,R_2,R_3\}$, we get the case (C).

Suppose that $T$
contains an extremal ray $Q$ different from $R_1,R_2,R_3$.
By condition, then $\{ R_1,R_2,R_3\}$ is elliptic. As we have seen,
$\nem (X,D(R_1)),\nem (X,D(R_2)),\nem (X,D(R_3))$ are contained in
$R_1+R_2+R_3$. Then like above (using Proposition 1.2.1 and Lemmas
1.2.8 and 1.2.9), we can prove that
$R_1+R_2+R_3$ is a $3$-dimensional face of Kodaira dimension $3$
of $\nem (X)$.

Let $Q\in T$ is different from $R_1,R_2,R_3$.
Since $R_1+R_2+R_3$ is a face of $\nem (X)$,
by statements (1) and (2) of Lemma 2.1.9,
the $Q$ is divisorially disjoint with $R_1$ and $R_2$.
Then, like above, using that $\nem (X,D(R_1))\subset R_1+R_2+R_3$,
we can prove that $R_1+Q$ is a $2$-dimensional face of Kodaira dimension
$3$ of $\nem (X)$. For $Q=R_2,R_3$ we have proven the same above.
Thus, this statement valid for any $Q\in T$ different from $R_1$.
Let us consider the contraction $f:X\to X^\prime$ of
$R_1$. By the statement above, $f(Q)$ is a divisorial extremal ray
for any $Q\in T,\ Q\not= R_1$. Evidently, $f(R_2)$ has the type (I).
Thus, considering $T^\prime =f(T)$ we get the case we had
considered above which gives rise the type (D) of Theorem.
Thus, $T^\prime$ has the type (D). One can easily see that then
$T$ has the type (C) and each extremal ray of $T$ does not belong to
a pair of the type $\bb_2$.

Let us prove the last statement.
Thus, we consider $T$ which has the type (B).
We use induction on $\# T$. For $\# T=2$, the statement
follows from Lemma 2.1.2.
If $\# T>2$, let us consider the contraction $f:X\to X^\prime$ of the
extremal ray $R_{11}$  which
we had used above when we considered this case.
Then $T^\prime =f(T)$ has the type (B) again and has the same
properties as $T$ but $\# T^\prime =n-1$.
By induction, $T^\prime = f(T)$
is contained in a face of Kodaira dimension $3$ of
$\nem (X^\prime)$. Then, evidently, $T$ is also contained in a face of
Kodaira dimension $3$ of $\nem (X)$.
This finishes the proof.
\enddemo

In Sect. 5, we will need another property of divisorial extremal rays
which shows importance of quasi-Lanner sets of divisorial extremal
rays.

\definition{Definition 2.2.7} A set $S=\{Q_1,...,Q_k\}$ of
divisorial extremal rays is called {\it semi-elliptic} if
for any non-negative $a_1,...,a_k$ there exists $j,\ 1\le j \le k,$
such that
$$
Q_j \cdot (a_1D(Q_1)+ \cdots + a_kD(Q_k))\le 0.
$$
\enddefinition

In particular, by Lemma 1.2.8, an elliptic set of
divisorial extremal rays is semi-elliptic.

\proclaim{Proposition 2.2.8} Let $\La=\{R_1,...,R_n\}$ be a set of
extremal rays of the type (I) or (II) with different divisors
$D(R_1),...,D(R_n)$. Assume that $\La$ is not semi-elliptic, i. e.
there are $a_1,...,a_n$ such that
$R_i\cdot (a_1D(R_1)+\cdots +a_nD(R_n))\ge 0$ for all $1\le i \le n$,
and one of these inequalities is strict.
Assume that any proper subset of $\La$ is semi-elliptic (i. e.
$\La$ is minimal non semi-elliptic).
Then $\La$ is quasi-Lanner (see Definition 2.2.3). Besides, any
proper subset of
$\La$ either is elliptic or is
connected parabolic with $\sharp \La -1$ elements.
\endproclaim

\demo{Proof} We use Proposition 2.2.4 as the equivalent definition of
an elliptic set of divisorial extremal rays.

If any proper subset  $\La^\prime\subset \La$ is elliptic, then
the set $\La$
is Lanner and is then quasi-Lanner.

Assume that $\La$ contains a
semi-elliptic subset  $\La^\prime$ which is not elliptic.
We should prove
that then $\La^\prime$ is connected parabolic and
$\sharp \La^\prime=\sharp \La -1$.
Let $\Cal P$ be a minimal semi-elliptic and not elliptic subset
of $\La^\prime$. Then any
proper subset of $\Cal P$ is elliptic (use Proposition 2.2.4).
It follows that $\Cal P$ is connected parabolic.

Let $R\in \La - {\Cal P}$.
Assume that there exists
an arrow from $R$ to (an element of)  $\Cal P$. Let
${\Cal P}=\{Q_1,...,Q_r\}$ and for some positive
$c_1,...,c_r$ we have ${\Cal P}\cdot (c_1D(Q_1)+\cdots +c_rD(Q_r))=0$.
Then for sufficiently big $\lambda >0$ we evidently have
$R\cdot (\lambda (c_1D(Q_1)+\cdots +c_rD(Q_r))+D(R))>0$
and ${\Cal P}\cdot (c_1D(Q_1)+\cdots +c_rD(Q_r))+D(R))\ge 0$.
Thus $\{ R \}\cup {\Cal P}$ is not semi-elliptic. Since $\La$ is
minimal non semi-elliptic, we get $\La =\{R \}\cup {\Cal P}$ and
$\La^\prime = {\Cal P}$ is connected parabolic with
$\sharp \La -1$ elements.

Thus, we can suppose that there is not  an arrow
from $\La - {\Cal P}$ to
${\Cal P}$. If $\La -{\Cal P}$ is not elliptic, arguing similarly,
we can find a connected parabolic subset
${\Cal P}^\prime \subset \La-{\Cal P}$. Moreover,
$\Cal P$, ${\Cal P}^\prime$ are divisorially disjoint to one another
and there is not an arrow from
$\La-({\Cal P}\cup {\Cal P}^\prime )$ to
${\Cal P}\cup {\Cal P}^\prime$. Continuing these considerations, we
prove that $\La$ is a disjoint union of connected
parabolic subsets ${\Cal P}_1,...,
{\Cal P}_m$ and an elliptic subset $\E$ such that connected
parabolic subsets ${\Cal P}_1,...,{\Cal P}_m$ are divisorially
disjoint to one another and there does not exist
an arrow from $\E$ to connected parabolic
subsets ${\Cal P}_1,...,{\Cal P}_m$.
If follows that $\La$ is semi-elliptic. We get a contradiction.
This finishes the proof.
\enddemo

 From our proof we also get

\proclaim{Proposition 2.2.9} Let $S$ be a set of
extremal rays of the type (I) or (II) with different divisors. Then
$S$ is semi-elliptic if and only if
$S$ is a disjoint union of connected
parabolic subsets ${\Cal P}_1,...,
{\Cal P}_m$ and an elliptic subset $\E$ such that
parabolic subsets ${\Cal P}_1,...,{\Cal P}_m$ are divisorially
disjoint to one another and there does not exist
an arrow from any element of $\E$ to any parabolic
subset ${\Cal P}_1,...,{\Cal P}_m$.
\endproclaim

\subhead
3. A refined variant of the Diagram Method
\endsubhead

Here we want to prove a more strong variant of
the Diagram Method Theorem which uses results of
both Sections 1 and 2. This variant will be useful below
for Calabi-Yau 3-folds.

We study 3-folds $X$ from the class $\LT$ with conditions
(a), (b) and (c) below:

(a) $X$ has a finite polyhedral Mori cone $\nem (X)$;

(b) $X$ does not have a small extremal ray;

(c) $\nem (X)$ does not have a face of numerical
Kodaira dimension $1$ or $2$.

 From (b) and (c), it follows that all extremal rays of $X$ are
divisorial of the type (I) or (II).

Using Theorem 2.2.6, we define analogous constants to
the constants $d(X)$, $C_1(X)$, $C_2(X)$ of Definition 1.3.1.

\definition{Definition 3.1}
We introduce some invariants of $X$.

Invariants $k$, $l$, $l_2$ are numbers
of divisorially connected components
of the types $\Aa_1$, $\cc_k,\ k\ge 2$, $\cc_2$ respectively
of the set of all special divisorial extremal rays on $X$.
See Theorem 2.1.7.

The invariants
$$
\split
&n(X)_D=\max_{F\in (D)}{\sharp F} - 1;\ \ \ \
n(X)_C=\max_{F\in (C)}{\sharp F}-1;\\
&n(X)_A=\max_{F\in (A)}{\sharp F}-1;\ \ \ \
d(X)_A=\max_{F\in (A)}{\text{diam\ } G(F)}.
\endsplit
$$
Here $F\in (A)$ (respectively $F\in (C)$ or $F\in (D)$)
means that $F$ runs through all $E$-sets of the type
(A) (respectively (C) or (D)) of Theorem 2.2.6
such that any proper subset of $F$
is extremal of Kodaira dimension $3$ and $F$ satisfies the
condition (iii) of Sect. 1.1 (in particular, $F$ is Lanner).
These invariants are
analogous to the invariant $d(X)$ of Definition 1.3.1.

Let $S$ be a set of divisorial extremal rays.
We denote as $S^\prime$
the subset of $S$ which one gets by removing all extremal rays
of the type (I) and all extremal rays of the type (II) which are
ends of single arrows
of $X$ (i.e. this extremal ray belongs to a component of the type
$\cc_k$ of the set of all special extremal rays on $X$
and is the end of a single arrow of this component). We define
a symmetric distance $\rho_A(R_1,R_2)$ in $S$ by the formula
$$
\rho _A(R_1, R_2)=
\cases
0,                         &\text{if $R_1=R_2$},\\
\rho(R_1,R_2)_{S^\prime},  &\text{if $\{R_1,R_2\} \subset S^\prime$,}\\
+\infty,                   &\text{otherwise}.
\endcases
$$
Here $\rho(R_1,R_2)_{S^\prime}$ denotes the ordinary distance
in the graph $G(S^\prime)$ using oriented paths.
The set $S^\prime$ does not have single arrows.
Thus, $\rho_A(R_1,R_2)$ is symmetric.

The invariant $C(X)_A$ is defined by the condition:
$$
\sharp \{ \{R_1, R_2\}\subset \E-\E_0 \mid
1 \le \rho_A (R_1,R_2)\le 2d(X)_A+1 \} \le C(X)_A \sharp (\E-\E_0)
$$
for any extremal set $\E$ of the Kodaira dimension $3$ of
divisorial extremal rays on $X$
and any its divisorially connected
subset $\E_0\subset \E$ of extremal rays of the type (II).
\enddefinition

We have the following refinement of Basic Theorem 1.3.2.

\proclaim{Basic Theorem 3.2} Let $X$ be a $3$-fold from the class $\LT$
and $X$ satisfies conditions (a), (b) and (c) above.
Then we have assertions (1), (2) and (3) below:

(1) $X$ does not have a pair of extremal rays of the type $\bb_2$ and
the Mori cone $\nem (X)$ is simplicial.
Thus, by Theorem
2.1.7, each divisorially connected component of the set of all
special divisorial extremal rays has
the type $\Aa_1$, $\cc_n$, $\Tt_3$ or $\Tt_3^\prime$.

(2) If the set of all special divisorial extremal rays on $X$ has
$k$ divisorially connected components of the type $\Aa_1$
and $l$ connected components of types $\cc_{n_1},...,\cc_{n_l}$
respectively, then
$$
k+(n_1-1)+...+(n_l-1)\le q(X)
$$
(see Definition 1.3.1 for the invariant $q(X)$).
Besides, if there exists a connected component of the type $\Tt_3$ or
$\Tt_3^\prime$, then every extremal ray of $X$ belongs to this connected
component and the Picard number $\rho (X) =3$ (in particular, $k=l=0$).

(3) We have the inequality:
$$
\split
\rho (X)=\dim N_1(X) &\le kn(X)_D+l_2\max {\{n(X)_C ,n(X)_A\}}+
8C(X)_A + 6 \\
&\le q(X)\max{\{n(X)_D, n(X)_C, n(X)_A \}} + 8C(X)_A+6.
\endsplit
$$
\endproclaim

\demo{Proof} (1) follows from the statement (1) of
Basic Theorem 1.3.2 and Theorem 2.1.7.

To prove second part of (2), we use the following Lemmas which are
important as itself (in fact, they are contained in Section 1)
and follow from the statement (1), Lemma 1.2.8, Lemma 1.1.1
and Theorem 1.2.13.

\proclaim{Lemma 3.3} Let $X$ be a $3$-fold from the class $\LT$
and $X$ satisfies the conditions (a), (b) and (c) above.
Let $T=\{R_1,...,R_n\}$ be a set of extremal rays of $X$.
Then

(i) If $T$ is extremal, there are positive $a_1,...,a_n$
such that $R_i\cdot (a_1D(R_1)+...+a_nD(R_n))<0$ for all $1\le i\le n$,
and for any non-negative $b_1,...,b_n$ which are not all equal to zero
there exists $i,\ 1\le i\le n$, such that
$R_i\cdot (b_1D(R_1)+...+b_nD(R_n))<0$.

(ii) If $T$ is not extremal, then $T$ contains an $E$-subset $\La$ and
there are non-negative $c_1,...,c_n$ which are not all equal to zero
such that $c_1D(R_1)+...+c_nD(R_n)$ is $nef$.
\endproclaim

\proclaim{Lemma 3.4} Let $X$ be a $3$-fold from the class $\LT$ and
$X$ satisfies the conditions (a), (b) and (c) above.
Let $T_1=\{R_1,...,R_n\}$ and
$T_2=\{Q_1,...,Q_m\}$ are two sets of extremal rays on $X$ such
that $R_i\cdot D(Q_j)=0$ for all $1 \le i \le n$
and all $1\le j \le m$ (equivalently, there does not
exist an arrow from $T_1$ to $T_2$). Then either $T_1$ or $T_2$ is
extremal.

If $T_1\cup T_2$ is the set of all extremal
rays on $X$ and $T_2\not=\emptyset$,
then $T_1$ is extremal and $T_2$ is not extremal.
\endproclaim

\demo{Proof} We explain only the last statement. Let us suppose
that $T_2$ is extremal. By condition and Lemma 3.3,
there are positive $a_1,...a_m$ such
that $R\cdot (a_1D(Q_1)+...+a_mD(Q_m))\le 0$ for any extremal ray
$R\in T_1\cup T_2$. By condition, $T_1\cup T_2$ is
the set of all extremal rays on $X$. It follows that,
$C\cdot (a_1D(Q_1)+...+a_mD(Q_m))\le 0$ for any curve $C$ on $X$.
Evidently, this is absurd. Thus, $T_2$ is not extremal.
By Lemma 3.3, there are non-negative $c_1,...,c_m$ such that
$D(Q)=c_1D(Q_1)+...+c_mD(Q_m)$ is $nef$. By condition, $T_1 \cdot D(Q)=0$.
It follows that $T_1$ is extremal.
\enddemo

 From the last statement of Lemma 3.4, it follows

\proclaim{Lemma 3.5} Let $X$ be a $3$-fold from the class $\LT$ and
$X$ satisfies the conditions (a), (b) and (c) above.
Then the set of all extremal rays on $X$ is divisorially connected and
is not extremal.
\endproclaim

Now, let us prove the statement (2) of Theorem 3.2.
Assume that the set of
special extremal rays has a connected component $P$ of the type
$\Tt_3$ or $\Tt_3^\prime$. Let $R$ be an extremal ray and
$R\not\in P$. By Theorem 2.1.8, then $R\cdot D(Q)=0$ for any extremal
ray $Q\in P$. We get a contradiction with Lemma 3.5.
Thus, any extremal ray of $X$ belongs to $P$. Since $P$ contains
exactly $3$ extremal rays, $\rho (X)\le 3$. By definition,
$P$ contains extremal rays $Q_1,Q_2, Q_3$ which define single arrows
$Q_1Q_2$ and $Q_2Q_3$. By Lemma 2.1.2, $Q_1+Q_2$ and $Q_2+Q_3$ are
$2$-dimensional faces of $\nem (X)$ which are evidently different.
It follows that $\rho (X)>2$. Thus, $\rho (X)=3$.

Let us prove
first part of (2). This is similar to the proof of the statement (2)
of Basic Theorem 1.3.2. Let $P_1=\{R_1\} ,...,P_k=\{R_k\}$ are connected
components of the type $\Aa_1$ and $S_1,...,S_l$ are connected
components of the types $\cc_{n_1},...,\cc_{n_l}$ of the set of
special extremal rays on $X$ where
$S_t=\{Q_{t1}, Q_{t2},...,Q_{tn_t}\}$ and
$Q_{t2}Q_{t1},...,Q_{tn_t}Q_{t1}$ are single arrows.

The set of all special extremal rays on $X$ evidently satisfies
the condition (i) of Lemma 3.3 and is then extremal.
Let $\E$ be a maximal extremal set of extremal rays on $X$
containing the set of all special extremal rays and such that
each divisorially connected component of $\E$ contains at least one special
extremal ray.

By Theorem 2.2.6, each connected component of $\E$ contains exactly
one connected component of the set of special extremal rays. Thus,
$\E$ contains
$k+l$ connected components $\E_1,...,\E_k, \E_{k+1},..., \E_{k+l}$.
Here $R_i\in \E_i$ and $\E_i$ has the type (D) for $1\le i \le k$.
And $S_j\subset \E_{k+j}$ and $\E_{k+j}$ either has the type (B) or (C)
for $1\le j \le l$. Changing numeration, we can suppose that
$\E_{k+j}$ has the type (C) for $1\le j \le r_1$ (in particular,
$n_1=...=n_{r_1}=2$) and $\E_{k+j}$ has the type (B) for $r_1<j\le l$.
For $r_1<j\le l$, we additionally denote as
$\E_{k+j,2},...,\E_{k+j,n_j}$ the corresponding chains containing
$Q_{j2},...,Q_{jn_j}$ respectively. These chains give
divisorially connected components of $\E_{k+j}-\{ Q_{j1}\}$.

Since $\E$ is extremal, by Lemma 3.3,
there are positive $a(R)$ for $R\in \E$
such that
$$
R\cdot \sum_{R\in \E}{a(R)D(R)}<0
$$
for any $R\in \E$.

Let us consider divisorially disjoint to one another
sets $F_v,\ v\in V$, of extremal rays which are
equal to one of sets $\E_1,...,\E_{k+r_1}$ or
$\E_{k+j,u}$, $r_1<j\le l,\ 2\le u\le n_j$, and the corresponding divisors
$$
D(F_v)=\sum_{R \in F_v}{a(R)D(R)}
$$
(here $\sharp V=k+(n_1-1)+...+(n_l-1)$).
One can easily see that
$C\cdot D(F_v)\le 0$ for any $R\in F_v$ and any curve $C\subset D(R)$.

Similarly to Lemma 1.3.4 and the proof of the statement
(2) of Basic Theorem 1.3.2, we can find an extremal set
$A=\{U_v \mid v\in V\}$ containing extremal rays
$U_v$ of the type (II) such that we have the
property: $U_v\cdot D(F_v)>0$ but $U_v\cdot D(F_{v^\prime})=0$ for
$v^\prime \not=v$. Similarly to the proof of the statement (2) of
Theorem 1.3.2, one proves that the graph $G(A)$ is full. Thus,
$\sharp A=\sharp V=k+(n_1-1)+...+(n_l-1)\le q(X)$.
This finishes the proof of the statement (2) of Theorem 3.2.

Proof of (3). It is based on the following

\proclaim{Lemma 3.6}
Let $X$ be a $3$-fold from the class $\LT$
and $X$ satisfies conditions (a), (b) and (c) above. Assume that
$\rho (X)>3$.
Then there exists
an extremal divisorially connected set $\E_0$
containing only extremal rays of the type
(II) and such
that any $E$-set $\La$ does not contain extremal rays of the type (I)
and extremal rays of the type (II) which are
terminal vertices of single arrows on $X$
(in particular, $\La$ has the type (A) of Theorem 2.2.6)
if $\La$ has at least two elements which do not belong to $\E_0$ and
$\La^\prime \cup \E_0$ is extremal for any proper subset
$\La^\prime \subset \La$.

Besides,
$$
\sharp \E_0 \le kn(X)_D+l_2\max{\{n(X)_C,n(X)_A \}}.
$$
\endproclaim

\demo{Proof}
We numerate as $R_1,...,R_k$ the whole set of
divisorial extremal rays of the type (I) and as
$R_{k+1},...,R_{k+l_2}$ the whole set of divisorial
extremal rays of the type (II) such that $R_{k+i}, 1\le i\le l_2,$
belongs to a connected component of the type $\cc_2$
of the set of all special extremal rays and $R_{k+i}$
is the terminal vertex of the single arrow of this component.

We construct $\E_0$ in $k+l_2$ steps as a sequence
$$
\emptyset = (\E_0)_0\subset ... \subset (\E_0)_{k+l_2}=\E_0.
$$
Here $(\E_0)_t, 1\le t \le k+l_2$, has the property that there does not
exist an $E$-set $\La$ which contains one of
extremal rays $R_1,...,R_t$, and $\La$ contains at least two elements which
do not belong to $(\E_0)_t$, and $(\E_0)_t\cup \La^\prime$ is extremal
for any proper subset $\La^\prime \subset \La$.

Suppose that we have constructed $(\E_0)_t$ with properties above
and $t<k+l_2$.
Assume that there exists an $E$-set $U$
which contains $R_{t+1}$, and $U$ contains at least two elements which
do not belong to $(\E_0)_t$, and $U^\prime \cup (\E_0)_t$ is
extremal for any proper subset $U^\prime \subset U$. Then we set
$$
(\E_0)_{t+1}=(\E_0)_t\cup (U-\{R_{t+1}\}).
$$
The $(\E_0)_{t+1}$ is extremal by the
conditions on $(\E_0)_t$ and $U$ above.
Let us suppose that there exists an $E$-set $\La$
which contains $R_i,\ 1\le i \le t+1$, and
$\La$ contains at least two elements which
do not belong to $(\E_0)_{t+1}$, and $\La^\prime \cup (\E_0)_{t+1}$ is
extremal for any proper subset $\La^\prime \subset \La$. If
$1\le i\le t$, we have similar properties for $(\E_0)_t$ and
get the contradiction. Thus, $i=t+1$. By our conditions,
there exists $Q\in \La-\{R_{t+1}\}$ such that
$Q\notin (\E_0)_{t+1}$ and $(\E_0)_{t+1}\cup (\La-\{Q\})$ is extremal. Since
$R_{t+1}\in \La-\{Q\}$ and
$U-\{R_{t+1}\}\subset (\E_0)_{t+1}$, it
follows that $U\subset (\E_0)_{t+1}\cup (\La-\{Q\})$ is extremal
because $(\E_0)_{t+1}\cup (\La-\{Q\})$ is extremal.
We get a contradiction since $U$ is an $E$-set. Thus,
$(\E_0)_{t+1}$ has desirable properties. Here, by construction and
Theorem 2.2.6, the set $U$ has the type (D)
if $1\le t+1\le k$, and the set $U$
has the type (C) or (A) if $k < t+1 \le k+l_2$.
If there does not exist the $E$-set $U$ above, we just put
$(\E_0)_{t+1}=(\E_0)_t$.

By our construction,
$\sharp \E_0 \le kn(X)_D+l_2\max{\{n(X)_C, n(X)_A \}}$.
By our construction, $\E_0$ is divisorially connected
since any two different $E$-sets are connected by arrows
(by Lemma 3.4) and extremal rays $R_1,...,R_{k+l_2}$ are
divisorially disjoint.

Let $R_{k+l_2+1},...,R_{k+l}$ be terminal vertices of arrows of
 all components
of the types $\cc_{k}, k\ge 3$. By our construction,
any $E$-set $\La$ does not contain
any of extremal rays $R_1,...,R_{k+l_2}$ if
$\La$ has at least two elements which do not belong to $\E_0$ and
$\La^\prime\cup \E_0$ is extremal for any proper subset
$\La^\prime \subset \La$.
We claim that we even have more:
the $E$-set $\La$ above does not contain any of extremal rays
$R_1,...,R_{k+l}$. Actually, if the $\E$-set $\La$ contains $R_i$,
$k+l_2 < i \le k+l$,
we get a contradiction with Theorems 2.1.8 and 2.2.6: the set $\La$
should be of the type (B), but then it is extremal and cannot be
an $E$-set.

Since $\rho (X)>3$, by the statement (2), any terminal vertex of a
single arrow on $X$ is one of extremal rays $R_{t+1},...,R_{t+l}$.
This finishes the proof of Lemma.

\enddemo

We continue the proof of Theorem 3.2.
Let us consider the face $\gamma \subset \M(X)=NEF(X)/{\Bbb R}^+$.
which is orthogonal to $\E_0$. Let us consider
constants
$C_1(X)^\prime$ and $C_2(X)^\prime$ which are defined by the properties:
$$
\sharp \{ (R_1, R_2)\in (\E-\E_0) \times (\E-\E_0)
\mid 1 \le \rho (R_1,R_2)\le d(X)_A \} \le C_1(X)^\prime \sharp (\E-\E_0);
$$
and
$$
\sharp \{ (R_1, R_2)\in (\E-\E_0) \times (\E-\E_0)
\mid d(X)_A+1\le \rho (R_1,R_2) \le 2d(X)_A+1\} \le C_2(X)^\prime
\sharp (\E-\E_0) .
$$
for any extremal set $\E$ which contains $\E_0$ and distance $\rho$ in the
graph $G(\E)$.

Directly applying to $\gamma$ Theorem 1.2 from \cite{N8}, we get
the inequality
$$
\dim \gamma <(16/3)C_1(X)^\prime + 4C_2(X)^\prime +6
\tag3-1
$$
which is less than we want.

To get the desirable estimate
$$
\dim \gamma < 8C(X)_A+6,
\tag3-2
$$
one needs to change a little the proof of Theorem 1.2 from \cite{N8}
for our case
(evidently, $C(X)_A \le (C_1(X)^\prime +C_2(X)^\prime )/2$ and
(3-2) implies (3-1)).

Let $\angle$ be an
oriented (plane) angle of $\gamma $.
Let $\R(\angle)$ be
the set of all extremal rays of $\nem (X)$ which are orthogonal
to the vertex of $\angle$.
The set
$\R(\angle)$ is a disjoint union
$$
\R(\angle)=\R(\angle^\perp)\cup \{ R_1(\angle)\} \cup \{ R_2(\angle)\}
$$
where $\R(\angle^\perp)$ contains all rays orthogonal to the plane of
the angle $\angle$, the rays $R_1(\angle )$ and $R_2(\angle )$
are orthogonal to the
first and second side of the oriented angle $\angle$, respectively.
Evidently, the set $\R(\angle)$ and the ordered pair
of rays $(R_1(\angle ), R_2(\angle ))$ define
the oriented angle $\angle$ uniquely.
We define the weight $\sigma_A (\angle )$ by the formula (for the proof
of Theorem 1.2 from \cite{N8}, the definition of the weight
$\sigma (\angle )$ was different):
$$
\sigma_A (\angle )=
\cases
1/2, &\text{if $1\le \rho_A(R_1(\angle ), R_2(\angle ))\le 2d(X)_A+1$,}\\
0, &\text{if $2d(X)_A+1 < \rho_A(R_1(\angle ), R_2(\angle ))$.}
\endcases
\tag3-3
$$
Similarly to the proof of Theorem 1.2 from \cite{N8},
one can check conditions of Lemma 1.4 (Vinberg's Lemma)
from \cite{N8} for the polyhedron ${\Cal M} = \gamma$
with the constants $C=C(X)_A$ and $D=0$.
The proof even is simpler because
the weight $\sigma_A (\angle^\prime)$ of the angle with opposite
orientation is equal to $\sigma_A (\angle)$. It follows (3-2).
For convenience of a reader, we recall
Vinberg's Lemma.

\proclaim{Lemma 3.8} Let $\Cal M$ be a convex simple polyhedron of
a dimension $n$. Let $C$ and $D$ are some numbers.
Suppose that oriented
angles ($2$-dimensional, plane) of $\Cal M$ are supplied with
weights and the following
conditions (1) and (2) hold:

(1) The sum of weights of all oriented angles at any vertex of
$\Cal M$ is not greater than $Cn+D$.

(2) The sum of weights of all oriented angles of any $2$-dimensional
face of $\Cal M$ is at least $5-k$ where $k$ is the number of
vertices of the $2$-dimensional face.

Then
$$
n<8C+5+
\cases
1+8D/n &\text{if $n$ is even,}\\
(8C+8D)/(n-1) &\text{if $n$ is odd}
\endcases .
$$
In particular, for $C\ge 0$ and $D=0$, we have
$$
n<8C+6.
$$
\endproclaim

Since $\dim \gamma =\dim N_1(X)-\sharp \E_0-1$ and
$\sharp \E_0\le kn(X)_D+l_2\max {\{n(X)_C ,n(X)_A\}}$, we
get the estimate $\dim N_1(X)\le
kn(X)_D+l_2\max {\{n(X)_C ,n(X)_A\}}$.
By the statement (2), $k+l_2\le q(X)$. Thus,
$\dim N_1(X)\le q(X)\max{\{n(X)_D, n(X)_C, n(X)_A\}} + 8C(X)_A+6$.
This finishes the proof of Theorem 3.2.
\enddemo

\remark{Remark 3.8} Let us consider Fano $3$-folds $X$ with
$\Bbb Q$-factorial terminal singularities and which satisfy
the conditions (a), (b) and (c). Then the constants
$q(X)\le 1$, $n(X)_D\le 1$, $n(X)_C=-1$, $n(X)_A\le 1$,
$d(X)_A=1$ and $C(X)_A=0$
(see \cite{N8}).
Thus, by Theorem 3.2, we get the main result of \cite{N8}:
$k+(n_1-1)+...+(n_l-1) \le 1$ and $\rho (X)\le 7$.
\endremark

In Sects. 4 and 5, we apply Basic Theorems 1.3.2 and 3.2 to
Calabi-Yau $3$-folds.

\head
4. Application to Calabi-Yau $3$-folds
\endhead

\subhead
4.1. Reminding
\endsubhead

Here we recall general facts about Calabi-Yau $3$-folds we
shall use.

A projective algebraic $3$-dimensional manifold
$X$ over $\Bbb C$
is called Calabi-Yau if the canonical class $K_X=0$ and
$h^1(X, {\Cal O}_X)=h^2(X, {\Cal O}_X)=0$.
Then $N^1(X)=H^2(X;\r)$.

It was shown by Wilson \cite{W1} and \cite{W2}
that the $nef$-cone $NEF(X)$ is locally rational polyhedral away from
the cubic intersection hypersurface
$$
{\Cal W}_X=\{ h\in H^2(X;\r)\mid h^3=0\}.
$$
The main tool here is to use results of
Kawamata \cite{Ka1} and Shokurov \cite{Sh} (which generalize
results of Mori \cite{Mo1}) about "polyhedrality" of
Mori cone of algebraic varieties with log-terminal singularities.
Thus, we can speak about faces $\gamma \subset NEF(X)$ considering
only faces of the locally polyhedral cone $NEF(X)$
away from the cone $\Cal W$. Then the
corresponding face $\gamma ^\prime = \gamma ^\perp$
of Mori cone $\nem (X)$ is finite
polyhedral and $\dim \gamma^\prime =
\codim (\gamma^\prime )^\perp $ where
$(\gamma^\prime)^\perp =\gamma$.
Evidently, these faces $\gamma \prime $ are exactly finite
polyhedral faces of $\nem (X)$ of numerical Kodaira dimension $3$ and
with the property $\codim (\gamma ^\prime )^\perp =\dim \gamma ^\prime$.
It is known (see \cite{W1}, \cite{W2}, \cite{O}) that
these faces are contractible and have Kodaira dimension $3$ (one
should apply theory of log-terminal extremal contractions of \cite{Ka1} and
\cite{Sh}; in fact, one need these results to prove
the results of Wilson we mentioned above).
Thus, we can speak about extremal rays of Kodaira dimension $3$
(they are orthogonal to faces of $NEF(X)$ of highest dimension) and
extremal rays of the types (I), (II) and (III) (small).
Applying results of \cite{Ka1} and \cite{Sh}, we get that a sequence
$X \to X^\prime $ of contractions of divisorial
extremal rays gives rise a $3$-fold with $\Bbb Q$-factorial (canonical)
singularities and $X^\prime$ inherits all properties of
$X$ we mentioned above. In particular, we get that
$X$ and $X^\prime $ belong to the class $\LT$.

The following fact is very important for us (see \cite{W2} and also
appendix of Shokurov): For an extremal ray $R$ of the type (II) of a
$3$-dimensional Calabi-Yau manifold $X$, there
exists a curve $C\in R$ such that
$$
C\cdot D(R)=-2.
\tag4-1-1
$$
In fact, this curve
$C$ is the general fiber of the $f:D(R)\to f(D(R))$ for the contraction
$f:X \to X^\prime $ of the extremal ray $R$. But, we don't need the last
fact in this paper.

Let $Q$ be another divisorial extremal ray of $X$. Using this curve $C\in R$,
we correspond to the arrow $RQ$ the weight $C\cdot D(Q)$. Thus, from now on,
arrows of the graph $G(T)$ of a set $T$ of extremal rays of the type (II)
have the corresponding weights. We don't show the corresponding
weight only if this weight is equal to $1$.

For an extremal ray $R$ of the type (I) we choose some curve $C\in R$.
This gives rise the corresponding weights of arrows too.
As was shown by Shokurov (see Appendix),
there exists a curve $C\in R$ such that $C\cdot D(R)=-1,-2$ or $-3$.
We shall use this result later. Now,
we fix the choice of the $C$ showing the weight $-k=C\cdot D(R)$
of the vertex corresponding to $R$. Thus, from now on, arrows and
black vertices of the graph $G(T)$ of any set
$T$ of extremal rays
of the type (I) or (II) are equipped by
weights. For a manifold $X$ these weights are integers.

Here we want to apply the theory we developed above for studying of
3-dimensional Calabi-Yau manifolds.

\subhead
4.2. Sets of divisorial extremal rays on Calabi-Yau manifolds
\endsubhead

We describe here elliptic, parabolic and Lanner (some
quasi-Lanner too) sets of extremal rays of the type (I) or (II) on
Calabi--Yau $3$-dimensional manifolds.
By Theorem 2.2.6, they have types (A), (B), (C), (D) or (E), (E').
By results of Section 1.2 and Theorem 2.2.6, this is
enough for the description of sets
of extremal rays of
the type (I) or (II) which are either
extremal of Kodaira dimension $3$ or are $E$-sets $\La$ such that
each proper subset of $\La$ is extremal of Kodaira dimension $3$
and $\La$ satisfies the condition (iii) of Sect. 1.1.

\subsubhead
4.2.1. Sets ot the type (A)
\endsubsubhead

We have the following results:

\proclaim{Theorem 4.2.1.1} Let $X$ be a $3$-dimensional
Calabi-Yau manifold. Let
$T=\{R_1,...,R_n\}$ be a set of
extremal rays of the type (II) on $X$ such that all divisors $D(R_1),...,
D(R_n)$ are different and the graph $G(T)$ does not have a single arrow
(thus, it has the type (A)) and is connected.

Then:

(a) The $T$ is elliptic if and only if $G(T)$ is a Dynkin
diagram of the root system type
$$
{\bold A}_n, {\bold B}_n, {\bold C}_n, {\bold D}_n, {\bold E}_6,
{\bold E}_7, {\bold E}_8, {\bold F}_4, {\bold G}_2
$$
of the Table 4.1.

(b) The $T$ is parabolic if and only if $G(T)$ is an extended
Dynkin diagram of the root system type
$$
\tilde{\bold  A}_n,
\tilde {\bold B}_n,
\widetilde{{\bold B}{\bold C}}_n,
\tilde {\bold C}_n,
\widetilde{{\bold B}{\bold D}}_n,
\widetilde{{\bold C}{\bold D}}_n,
\tilde{\bold D}_n,
\tilde{\bold E}_6,
\tilde{\bold E}_7,
\tilde{\bold E}_8,
\widetilde{{\bold B}{\bold F}}_4,
\widetilde{{\bold C}{\bold F}}_4,
\widetilde{{\bold A}{\bold G}}_2,
\widetilde{{\bold G}{\bold A}}_2,
$$
of the Table 4.2 where weights of vertices
show the coefficients $a_i$ of Definition 2.2.2.

(c) The $T$ is Lanner if and only if $G(T)$ is one of diagrams of
the Table 4.3.
\endproclaim

\demo{Proof} Applying (4-1-1), this is standard and
well-known in fact.
\enddemo

\subsubhead
4.2.2. Sets ot the type (B)
\endsubsubhead

\proclaim{Theorem 4.2.2.1} Let $X$ be a $3$-dimensional
Calabi-Yau manifold. Let
$T=\{R_1,...,R_n\}$ be a set of
extremal rays of the type (II) on $X$ such that all divisors $D(R_1),...,
D(R_n)$ are different and the graph $G(T)$ has the type (B) of
Theorem 2.2.2. Thus,
$T=\{ R\}\cup C_1\cup ... \cup C_k$, $k\ge 1$, where
$$
C_1=\{ R_{11},...,R_{1n_1}\},\  C_2=\{ R_{21},...,R_{2n_2}\},...,
C_k=\{ R_{k1},...,R_{kn_k}\}
$$
are divisorially disjoint to one another chains, and all arrows between
$R$ and extremal rays of these
chains are single arrows $R_{j1}R,\ j=1,...,k$.

Then the $T$ cannot be parabolic or Lanner. The $T$ is elliptic if
and only if
each chain $C_1,...,C_k$ is one of the chains
$$
{\bold A}_n, {\bold B}_n, {\bold C}_n, {\bold F}_4, {\bold G}_2
$$
of Theorem 4.2.1.1, weights of single arrows $R_{j1}R,\ j=1,...,k$,
are arbitrary natural numbers.
\endproclaim

\demo{Proof} If $T$ is either elliptic, or parabolic or Lanner,
the chains $C_1,...,C_k$ are elliptic because they give proper subsets of $T$.
By Theorem 4.2.1.1, they have the types
${\bold A}_n, {\bold B}_n, {\bold C}_n, {\bold F}_4$ or ${\bold G}_2$.
One can easily show (see the proof of Theorem 2.2.6)
that then $T$ is elliptic for arbitrary weights of arrows
$R_{j1}R$ which should be natural numbers. This finishes the proof.
\enddemo

\subsubhead
4.2.3. Sets of the type (C)
\endsubsubhead

\proclaim{Theorem 4.2.3.1} Let $X$ be a $3$-dimensional
Calabi-Yau manifold. Let
$T=\{R_1,...,R_n\}$ be a set of
extremal rays of the type (II) on $X$ such that all divisors $D(R_1),...,
D(R_n)$ are different and the graph $G(T)$ has the type (C) of
Theorem 2.2.6.

Then $T$ is elliptic, parabolic, Lanner or quasi-Lanner if and only if
$G(T)$ is elliptic, parabolic, Lanner or quasi-Lanner diagram
respectively of the Table 4.4 below.
\endproclaim

\demo{Proof} The corresponding calculations are very simple using
(4-1-1).
\enddemo

\subsubhead
4.2.4. Sets of the type (D)
\endsubsubhead

\proclaim{Theorem 4.2.4.1} Let $X$ be a $3$-dimensional
Calabi-Yau manifold. Let $T=\{R_1,...,R_n\}$ be a set of
extremal rays of the type (I) or (II) on $X$ such that
the graph $G(T)$ has the type (D) of
Theorem 2.2.6.

Then $T$ is elliptic, parabolic, Lanner or quasi-Lanner if and only if
$G(T)$ is elliptic, parabolic, Lanner or quasi-Lanner diagram
respectively of the Table 4.5 below (we recall that the black vertex
corresponds to an extremal ray of the type (I)).
\endproclaim

\demo{Proof} The corresponding calculations are very simple using
(4-1-1). We remark that $G(T)$ without the black vertex should have
one of the types
${\bold A}_n$, ${\bold B}_n$, ${\bold C}_n$, ${\bold F}_4$ or ${\bold G}_2$
by Theorem 3.2.1.1 (a).
\enddemo

\subsubhead
4.2.5. Sets of the type (E)
\endsubsubhead

\proclaim{Theorem 4.2.5.1} Let $X$ be a $3$-dimensional
Calabi-Yau manifold. Let $T=\{R_1,R_2,R_3\}$ be a set of
extremal rays of the type (II) on $X$ such that
the graph $G(T)$ has the type (E) (triangle or special triangle)
of Theorem 2.2.6.

Then $T$ is elliptic, parabolic, Lanner or quasi-Lanner if and only if
$G(T)$ is elliptic, parabolic, Lanner or quasi-Lanner diagram
respectively of the Table 4.6 below.
\endproclaim

\demo{Proof} The corresponding calculations are very simple using
(4-1-1).
\enddemo

\newpage

\centerline{\hfil
\hbox to 245pt{
\hfil
\hbox to15pt{\hfill \vbox to30pt{\vfill ${\bold A}_n\ :$ \vfill}\hfill}
\hskip10pt plus2pt minus2pt
\hbox to15pt{\hfill \vbox to30pt{\vfill $\bigcirc$ \vfill}\hfill}
\hskip8pt plus2pt minus2pt
\hbox to30pt{\hfill
\vbox to30pt{\vfill
\nointerlineskip
\hbox to25pt{\rightarrowfill} \nointerlineskip
\hbox to25pt{\leftarrowfill} \nointerlineskip
\vfill}
\hfill}
\hskip-14pt plus2pt minus2pt
\hbox to15pt{\hfill \vbox to30pt{\vfill $\bigcirc$ \vfill}\hfill}
\hskip8pt plus2pt minus2pt
\hbox to30pt{\hfill
\vbox to30pt{\vfill \nointerlineskip
\hbox to25pt{\rightarrowfill} \nointerlineskip
\hbox to25pt{\leftarrowfill} \nointerlineskip
\vfill}
\hfill}
\hskip-14pt plus2pt minus2pt
\hbox to15pt{\hfill \vbox to30pt{\vfill $\bigcirc$ \vfill}\hfill}
\hskip5pt plus2pt minus2pt
\hbox to15pt{\hfill \vbox to30pt{\vfill $\cdots$ \vfill}\hfill}
\hskip5pt plus2pt minus2pt
\hbox to15pt{\hfill \vbox to30pt{\vfill $\bigcirc$ \vfill}\hfill}
\hskip8pt plus2pt minus2pt
\hbox to30pt{\hfill
\vbox to30pt{\vfill \nointerlineskip
\hbox to25pt{\rightarrowfill} \nointerlineskip
\hbox to25pt{\leftarrowfill} \nointerlineskip
\vfill}
\hfill}
\hskip-14pt plus2pt minus2pt
\hbox to15pt{\hfill \vbox to30pt{\vfill $\bigcirc$ \vfill}\hfill}
\hskip8pt plus2pt minus2pt
\hbox to30pt{\hfill
\vbox to30pt{\vfill \nointerlineskip
\hbox to25pt{\rightarrowfill} \nointerlineskip
\hbox to25pt{\leftarrowfill} \nointerlineskip
\vfill}
\hfill}
\hskip-14pt plus2pt minus2pt
\hbox to15pt{\hfill \vbox to30pt{\vfill $\bigcirc$ \vfill}\hfill}
\hfil
}
\hskip55pt}


\centerline{\hfil
\hbox to 245pt{
\hfil
\hbox to15pt{\hfill \vbox to30pt{\vfill ${\bold B}_n\ :$ \vfill}\hfill}
\hskip10pt plus2pt minus2pt
\hbox to15pt{\hfill \vbox to30pt{\vfill $\bigcirc$ \vfill}\hfill}
\hskip8pt plus2pt minus2pt
\hbox to30pt{\hfill
\vbox to30pt{\vskip8pt \nointerlineskip
\hbox to25pt{\rightarrowfill} \nointerlineskip
\hbox to25pt{\leftarrowfill} \nointerlineskip
\hbox to25pt{\hfill $2$ \hfill}
\vfill}
\hfill}
\hskip-14pt plus2pt minus2pt
\hbox to15pt{\hfill \vbox to30pt{\vfill $\bigcirc$ \vfill}\hfill}
\hskip8pt plus2pt minus2pt
\hbox to30pt{\hfill
\vbox to30pt{\vfill \nointerlineskip
\hbox to25pt{\rightarrowfill} \nointerlineskip
\hbox to25pt{\leftarrowfill} \nointerlineskip
\vfill}
\hfill}
\hskip-14pt plus2pt minus2pt
\hbox to15pt{\hfill \vbox to30pt{\vfill $\bigcirc$ \vfill}\hfill}
\hskip5pt plus2pt minus2pt
\hbox to15pt{\hfill \vbox to30pt{\vfill $\cdots$ \vfill}\hfill}
\hskip5pt plus2pt minus2pt
\hbox to15pt{\hfill \vbox to30pt{\vfill $\bigcirc$ \vfill}\hfill}
\hskip8pt plus2pt minus2pt
\hbox to30pt{\hfill
\vbox to30pt{\vfill \nointerlineskip
\hbox to25pt{\rightarrowfill} \nointerlineskip
\hbox to25pt{\leftarrowfill} \nointerlineskip
\vfill}
\hfill}
\hskip-14pt plus2pt minus2pt
\hbox to15pt{\hfill \vbox to30pt{\vfill $\bigcirc$ \vfill}\hfill}
\hskip8pt plus2pt minus2pt
\hbox to30pt{\hfill
\vbox to30pt{\vfill \nointerlineskip
\hbox to25pt{\rightarrowfill} \nointerlineskip
\hbox to25pt{\leftarrowfill} \nointerlineskip
\vfill}
\hfill}
\hskip-14pt plus2pt minus2pt
\hbox to15pt{\hfill \vbox to30pt{\vfill $\bigcirc$ \vfill}\hfill}
\hfil
}
\hskip55pt}


\centerline{\hfil
\hbox to 245pt{
\hfil
\hbox to15pt{\hfill \vbox to30pt{\vfill ${\bold C}_n\ :$ \vfill}\hfill}
\hskip10pt plus2pt minus2pt
\hbox to15pt{\hfill \vbox to30pt{\vfill $\bigcirc$ \vfill}\hfill}
\hskip8pt plus2pt minus2pt
\hbox to30pt{\hfill
\vbox to30pt{\vfill
\hbox to25pt{\hfill $2$ \hfill}
\nointerlineskip
\hbox to25pt{\rightarrowfill} \nointerlineskip
\hbox to25pt{\leftarrowfill} \nointerlineskip
\vskip8pt}
\hfill}
\hskip-14pt plus2pt minus2pt
\hbox to15pt{\hfill \vbox to30pt{\vfill $\bigcirc$ \vfill}\hfill}
\hskip8pt plus2pt minus2pt
\hbox to30pt{\hfill
\vbox to30pt{\vfill \nointerlineskip
\hbox to25pt{\rightarrowfill} \nointerlineskip
\hbox to25pt{\leftarrowfill} \nointerlineskip
\vfill}
\hfill}
\hskip-14pt plus2pt minus2pt
\hbox to15pt{\hfill \vbox to30pt{\vfill $\bigcirc$ \vfill}\hfill}
\hskip5pt plus2pt minus2pt
\hbox to15pt{\hfill \vbox to30pt{\vfill $\cdots$ \vfill}\hfill}
\hskip5pt plus2pt minus2pt
\hbox to15pt{\hfill \vbox to30pt{\vfill $\bigcirc$ \vfill}\hfill}
\hskip8pt plus2pt minus2pt
\hbox to30pt{\hfill
\vbox to30pt{\vfill \nointerlineskip
\hbox to25pt{\rightarrowfill} \nointerlineskip
\hbox to25pt{\leftarrowfill} \nointerlineskip
\vfill}
\hfill}
\hskip-14pt plus2pt minus2pt
\hbox to15pt{\hfill \vbox to30pt{\vfill $\bigcirc$ \vfill}\hfill}
\hskip8pt plus2pt minus2pt
\hbox to30pt{\hfill
\vbox to30pt{\vfill \nointerlineskip
\hbox to25pt{\rightarrowfill} \nointerlineskip
\hbox to25pt{\leftarrowfill} \nointerlineskip
\vfill}
\hfill}
\hskip-14pt plus2pt minus2pt
\hbox to15pt{\hfill \vbox to30pt{\vfill $\bigcirc$ \vfill}\hfill}
\hfil
}
\hskip55pt}


\centerline{\hfil
\hbox to 290pt{
\hfil
\hbox to15pt{\hfill \vbox to60pt{\vfill ${\bold D}_n\ :$ \vfill}\hfill}
\hskip25pt plus2pt minus2pt
\hbox to 20pt{\hfill
\vbox to 60pt {\vfill
\hbox to 15pt{\hfill $\bigcirc$ \hfill}
\nointerlineskip\vfill
\hbox to15pt{\hfill $\uparrow\downarrow$ \hfill}
\nointerlineskip\vfill
\hbox to15pt{\hfill $\bigcirc$ \hfill}
\nointerlineskip\vfill
\hbox to15pt{\hfill $\uparrow\downarrow$ \hfill}
\nointerlineskip\vfill
\hbox to15pt{\hfill $\bigcirc$ \hfill}
\vfill}
\hfill}
\hskip-10pt plus2pt minus2pt
\hbox to30pt{\hfill
\vbox to60pt{\vfill
\nointerlineskip
\hbox to25pt{\rightarrowfill} \nointerlineskip
\hbox to25pt{\leftarrowfill} \nointerlineskip
\vfill}
\hfill}
\hskip-14pt plus2pt minus2pt
\hbox to15pt{\hfill \vbox to60pt{\vfill $\bigcirc$ \vfill}\hfill}
\hskip5pt plus2pt minus2pt
\hbox to15pt{\hfill \vbox to60pt{\vfill $\cdots$ \vfill}\hfill}
\hskip5pt plus2pt minus2pt
\hbox to15pt{\hfill \vbox to60pt{\vfill $\bigcirc$ \vfill}\hfill}
\hskip8pt plus2pt minus2pt
\hbox to30pt{\hfill
\vbox to60pt{\vfill \nointerlineskip
\hbox to25pt{\rightarrowfill} \nointerlineskip
\hbox to25pt{\leftarrowfill} \nointerlineskip
\vfill}
\hfill}
\hskip-14pt plus2pt minus2pt
\hbox to15pt{\hfill \vbox to60pt{\vfill $\bigcirc$ \vfill}\hfill}
\hskip8pt plus2pt minus2pt
\hbox to30pt{\hfill
\vbox to60pt{\vfill \nointerlineskip
\hbox to25pt{\rightarrowfill} \nointerlineskip
\hbox to25pt{\leftarrowfill} \nointerlineskip
\vfill}
\hfill}
\hskip-14pt plus2pt minus2pt
\hbox to15pt{\hfill \vbox to60pt{\vfill $\bigcirc$ \vfill}\hfill}
\hfil
}
\hskip50pt}


\centerline{\hfil
\hbox to 220pt{
\hfil
\hbox to15pt{\hfill \vbox to30pt{\vskip2pt${\bold E}_6\ :$ \vfill}\hfill}
\hskip10pt plus2pt minus2pt
\hbox to15pt{\hfill\vbox to30pt{\vskip2pt $\bigcirc$ \vfill}\hfill}
\hskip15pt plus2pt minus2pt
\hbox to30pt{\hfill
\vbox to30pt{
\nointerlineskip
\hbox to25pt{\rightarrowfill} \nointerlineskip
\hbox to25pt{\leftarrowfill} \nointerlineskip
\vfill}
\hfill}
\hskip-14pt plus2pt minus2pt
\hbox to15pt{\hfill \vbox to30pt{\vskip2pt $\bigcirc$ \vfill}\hfill}
\hskip8pt plus2pt minus2pt
\hbox to30pt{\hfill
\vbox to30pt{\nointerlineskip
\hbox to25pt{\rightarrowfill} \nointerlineskip
\hbox to25pt{\leftarrowfill} \nointerlineskip
\vfill}
\hfill}
\hskip-5pt plus2pt minus2pt
\hbox to 20pt{\hfill
\vbox to 30pt {\vskip2pt
\hbox to15pt{\hfill $\bigcirc$ \hfill}
\nointerlineskip\vfill
\hbox to15pt{\hfill $\uparrow\downarrow$ \hfill}
\nointerlineskip\vfill
\hbox to15pt{\hfill $\bigcirc$ \hfill}
}
\hfill}
\hskip-10pt plus2pt minus2pt
\hbox to30pt{\hfill
\vbox to30pt{\nointerlineskip
\hbox to25pt{\rightarrowfill} \nointerlineskip
\hbox to25pt{\leftarrowfill} \nointerlineskip
\vfill}
\hfill}
\hskip-14pt plus2pt minus2pt
\hbox to15pt{\hfill \vbox to30pt{\vskip2pt $\bigcirc$ \vfill}\hfill}
\hskip8pt plus2pt minus2pt
\hbox to30pt{\hfill
\vbox to30pt{\nointerlineskip
\hbox to25pt{\rightarrowfill} \nointerlineskip
\hbox to25pt{\leftarrowfill} \nointerlineskip
\vfill}
\hfill}
\hskip-14pt plus2pt minus2pt
\hbox to15pt{\hfill \vbox to30pt{\vskip2pt $\bigcirc$ \vfill}\hfill}
\hfil
}
\hskip80pt}


\centerline{\hfil
\hbox to 250pt{
\hfil
\hbox to15pt{\hfill \vbox to30pt{\vskip2pt${\bold E}_7\ :$ \vfill}\hfill}
\hskip10pt plus2pt minus2pt
\hbox to15pt{\hfill\vbox to30pt{\vskip2pt $\bigcirc$ \vfill}\hfill}
\hskip15pt plus2pt minus2pt
\hbox to30pt{\hfill
\vbox to30pt{
\nointerlineskip
\hbox to25pt{\rightarrowfill} \nointerlineskip
\hbox to25pt{\leftarrowfill} \nointerlineskip
\vfill}
\hfill}
\hskip-14pt plus2pt minus2pt
\hbox to15pt{\hfill \vbox to30pt{\vskip2pt $\bigcirc$ \vfill}\hfill}
\hskip8pt plus2pt minus2pt
\hbox to30pt{\hfill
\vbox to30pt{\nointerlineskip
\hbox to25pt{\rightarrowfill} \nointerlineskip
\hbox to25pt{\leftarrowfill} \nointerlineskip
\vfill}
\hfill}
\hskip-5pt plus2pt minus2pt
\hbox to 20pt{\hfill
\vbox to 30pt {\vskip2pt
\hbox to15pt{\hfill $\bigcirc$ \hfill}
\nointerlineskip\vfill
\hbox to15pt{\hfill $\uparrow\downarrow$ \hfill}
\nointerlineskip\vfill
\hbox to15pt{\hfill $\bigcirc$ \hfill}
}
\hfill}
\hskip-10pt plus2pt minus2pt
\hbox to30pt{\hfill
\vbox to30pt{\nointerlineskip
\hbox to25pt{\rightarrowfill} \nointerlineskip
\hbox to25pt{\leftarrowfill} \nointerlineskip
\vfill}
\hfill}
\hskip-14pt plus2pt minus2pt
\hbox to15pt{\hfill \vbox to30pt{\vskip2pt $\bigcirc$ \vfill}\hfill}
\hskip8pt plus2pt minus2pt
\hbox to30pt{\hfill
\vbox to30pt{\nointerlineskip
\hbox to25pt{\rightarrowfill} \nointerlineskip
\hbox to25pt{\leftarrowfill} \nointerlineskip
\vfill}
\hfill}
\hskip-14pt plus2pt minus2pt
\hbox to15pt{\hfill \vbox to30pt{\vskip2pt $\bigcirc$ \vfill}\hfill}
\hskip8pt plus2pt minus2pt
\hbox to30pt{\hfill
\vbox to30pt{\nointerlineskip
\hbox to25pt{\rightarrowfill} \nointerlineskip
\hbox to25pt{\leftarrowfill} \nointerlineskip
\vfill}
\hfill}
\hskip-14pt plus2pt minus2pt
\hbox to15pt{\hfill \vbox to30pt{\vskip2pt $\bigcirc$ \vfill}\hfill}
\hfil
}
\hskip50pt}

\centerline{\hfil
\hbox to 290pt{
\hfil
\hbox to15pt{\hfill \vbox to30pt{\vskip2pt${\bold E}_8\ :$ \vfill}\hfill}
\hskip10pt plus2pt minus2pt
\hbox to15pt{\hfill\vbox to30pt{\vskip2pt $\bigcirc$ \vfill}\hfill}
\hskip15pt plus2pt minus2pt
\hbox to30pt{\hfill
\vbox to30pt{
\nointerlineskip
\hbox to25pt{\rightarrowfill} \nointerlineskip
\hbox to25pt{\leftarrowfill} \nointerlineskip
\vfill}
\hfill}
\hskip-14pt plus2pt minus2pt
\hbox to15pt{\hfill \vbox to30pt{\vskip2pt $\bigcirc$ \vfill}\hfill}
\hskip8pt plus2pt minus2pt
\hbox to30pt{\hfill
\vbox to30pt{\nointerlineskip
\hbox to25pt{\rightarrowfill} \nointerlineskip
\hbox to25pt{\leftarrowfill} \nointerlineskip
\vfill}
\hfill}
\hskip-5pt plus2pt minus2pt
\hbox to 20pt{\hfill
\vbox to 30pt {\vskip2pt
\hbox to15pt{\hfill $\bigcirc$ \hfill}
\nointerlineskip\vfill
\hbox to15pt{\hfill $\uparrow\downarrow$ \hfill}
\nointerlineskip\vfill
\hbox to15pt{\hfill $\bigcirc$ \hfill}
}
\hfill}
\hskip-10pt plus2pt minus2pt
\hbox to30pt{\hfill
\vbox to30pt{\nointerlineskip
\hbox to25pt{\rightarrowfill} \nointerlineskip
\hbox to25pt{\leftarrowfill} \nointerlineskip
\vfill}
\hfill}
\hskip-14pt plus2pt minus2pt
\hbox to15pt{\hfill \vbox to30pt{\vskip2pt $\bigcirc$ \vfill}\hfill}
\hskip8pt plus2pt minus2pt
\hbox to30pt{\hfill
\vbox to30pt{\nointerlineskip
\hbox to25pt{\rightarrowfill} \nointerlineskip
\hbox to25pt{\leftarrowfill} \nointerlineskip
\vfill}
\hfill}
\hskip-14pt plus2pt minus2pt
\hbox to15pt{\hfill \vbox to30pt{\vskip2pt $\bigcirc$ \vfill}\hfill}
\hskip8pt plus2pt minus2pt
\hbox to30pt{\hfill
\vbox to30pt{\nointerlineskip
\hbox to25pt{\rightarrowfill} \nointerlineskip
\hbox to25pt{\leftarrowfill} \nointerlineskip
\vfill}
\hfill}
\hskip-14pt plus2pt minus2pt
\hbox to15pt{\hfill \vbox to30pt{\vskip2pt $\bigcirc$ \vfill}\hfill}
\hskip8pt plus2pt minus2pt
\hbox to30pt{\hfill
\vbox to30pt{\nointerlineskip
\hbox to25pt{\rightarrowfill} \nointerlineskip
\hbox to25pt{\leftarrowfill} \nointerlineskip
\vfill}
\hfill}
\hskip-14pt plus2pt minus2pt
\hbox to15pt{\hfill \vbox to30pt{\vskip2pt $\bigcirc$ \vfill}\hfill}
\hfil
}
\hskip10pt}


\centerline{\hfil
\hbox to 165pt{
\hfil
\hbox to15pt{\hfill \vbox to30pt{\vfill ${\bold F}_4\ :$ \vfill}\hfill}
\hskip10pt plus2pt minus2pt
\hbox to15pt{\hfill \vbox to30pt{\vfill $\bigcirc$ \vfill}\hfill}
\hskip8pt plus2pt minus2pt
\hbox to30pt{\hfill
\vbox to30pt{\vfill
\hbox to25pt{\rightarrowfill} \nointerlineskip
\hbox to25pt{\leftarrowfill}
\vfill}
\hfill}
\hskip-14pt plus2pt minus2pt
\hbox to15pt{\hfill \vbox to30pt{\vfill $\bigcirc$ \vfill}\hfill}
\hskip8pt plus2pt minus2pt
\hbox to30pt{\hfill
\vbox to30pt{\vfill \hbox to25pt{\hfill $2$ \hfill} \nointerlineskip
\hbox to25pt{\rightarrowfill} \nointerlineskip
\hbox to25pt{\leftarrowfill}
\vskip8pt}
\hfill}
\hskip-14pt plus2pt minus2pt
\hbox to15pt{\hfill \vbox to30pt{\vfill $\bigcirc$ \vfill}\hfill}
\hskip8pt plus2pt minus2pt
\hbox to30pt{\hfill
\vbox to30pt{\vfill
\hbox to25pt{\rightarrowfill} \nointerlineskip
\hbox to25pt{\leftarrowfill}
\vfill}
\hfill}
\hskip-14pt plus2pt minus2pt
\hbox to15pt{\hfill \vbox to30pt{\vfill $\bigcirc$ \vfill}\hfill}
}
\hskip135pt}


\centerline{\hfil
\hbox to 95pt{\hfil
\hbox to15pt{\hfill \vbox to30pt{\vfill ${\bold G}_2\ :$ \vfill}\hfill}
\hskip10pt plus2pt minus2pt
\hbox to15pt{\hfill \vbox to30pt{\vfill $\bigcirc$ \vfill}\hfill}
\hskip8pt plus2pt minus2pt
\hbox to30pt{\hfill
\vbox to30pt{\vfill \hbox to25pt{\hfill $3$ \hfill} \nointerlineskip
\hbox to25pt{\rightarrowfill} \nointerlineskip
\hbox to25pt{\leftarrowfill}
\vskip8pt}
\hfill}
\hskip-14pt plus2pt minus2pt
\hbox to15pt{\hfill \vbox to30pt{\vfill $\bigcirc$ \vfill}\hfill}
}
\hskip205pt}

\vskip20pt

\centerline{{\bf Table 4.1.} Calabi--Yau elliptic diagrams without
single arrows}
\centerline{(classical Dynkin diagrams).}

\newpage

\line{
\hbox{\hfill
\hbox {\hfill \vbox to30pt
{\vfill \hbox{${\bold A}_1(a, b),\  ab=4: $}
\vfill}\hfill}
\hskip10pt plus 2pt minus 2pt
\vertex{2}
\hskip0pt plus 2pt minus 2pt
\darrupdown{a}{b}
\hskip0pt plus 2pt minus 2pt
\vertex{b}}
\hskip20pt
\hbox{
\hbox to 30pt{\hfill \vbox to50pt{\vfill $\tilde {\bold A}_n,\ n>1:$
\vfill}\hfill}
\hskip15pt plus2pt minus2pt
\hbox to 3pt{\vbox to50pt{$1$ \vfill $1$}}
\hskip15pt plus2pt minus 2pt
\vbox to 50pt{
\hbox to15pt{\hfill $\bigcirc$ \hfill}
\nointerlineskip\vfill
\hbox to15pt{\hfill $\uparrow\downarrow$ \hfill}
\nointerlineskip\vfill
\hbox to15pt{\hfill $\cdot$ \hfill}
\nointerlineskip\vskip1pt
\hbox to15pt{\hfill $\cdot$ \hfill}
\nointerlineskip\vskip1pt
\hbox to15pt{\hfill $\cdot$ \hfill}
\nointerlineskip\vfill
\hbox to15pt{\hfill $\uparrow\downarrow$ \hfill}
\nointerlineskip\vfill
\hbox to15pt{\hfill $\bigcirc$ \hfill}
}
\vbox to50pt{
\hbox to45pt{\rightarrowfill} \nointerlineskip
\hbox to45pt{\leftarrowfill} \nointerlineskip\vfill
\hbox to45pt{\rightarrowfill} \nointerlineskip
\hbox to45pt{\leftarrowfill}}
\vbox to 50pt{
\hbox to15pt{\hfill $\bigcirc$ \hfill}
\nointerlineskip\vfill
\hbox to15pt{\hfill $\uparrow\downarrow$ \hfill}
\nointerlineskip\vfill
\hbox to15pt{\hfill $\cdot$ \hfill}
\nointerlineskip\vskip1pt
\hbox to15pt{\hfill $\cdot$ \hfill}
\nointerlineskip\vskip1pt
\hbox to15pt{\hfill $\cdot$ \hfill}
\nointerlineskip\vfill
\hbox to15pt{\hfill $\uparrow\downarrow$ \hfill}
\nointerlineskip\vfill
\hbox to15pt{\hfill $\bigcirc$ \hfill}
}
\hskip-15pt plus2pt minus2pt
\hbox to3pt{\vbox to50pt{$1$ \vfill $1$}}
}
\hfill}


\centerline{\hfil
\hbox to300pt{
\hfil
\hbox to15pt{\hfill \vbox
to30pt{\vfill $\tilde {\bold B}_n\ :$ \vfill}\hfill}
\hskip20pt plus2pt minus2pt
\hbox to15pt{\hfill \vbox to30pt{
\hbox to15pt{\hfill $1$ \hfill}
\nointerlineskip\vfill
\hbox to15pt{\hfill $\bigcirc$ \hfill}
\nointerlineskip
\vskip10pt}\hfill}
\hskip5pt plus2pt minus2pt
\hbox to30pt{\hfill
\vbox to30pt{\vskip8pt \nointerlineskip
\hbox to25pt{\rightarrowfill} \nointerlineskip
\hbox to25pt{\leftarrowfill} \nointerlineskip
\hbox to25pt{\hfill $2$ \hfill}
\vfill}
\hfill}
\hskip0pt plus2pt minus2pt
\hbox to15pt{\hfill \vbox to30pt{
\hbox to15pt{\hfill $2$ \hfill}
\nointerlineskip\vfill
\hbox to15pt{\hfill $\bigcirc$ \hfill}
\nointerlineskip
\vskip10pt}\hfill}
\hskip0pt plus2pt minus2pt
\hbox to30pt{\hfill
\vbox to30pt{\vfill \nointerlineskip
\hbox to25pt{\rightarrowfill} \nointerlineskip
\hbox to25pt{\leftarrowfill} \nointerlineskip
\vfill}
\hfill}
\hskip0pt plus2pt minus2pt
\hbox to15pt{\hfill \vbox to30pt{
\hbox to15pt{\hfill $2$ \hfill}
\nointerlineskip\vfill
\hbox to15pt{\hfill $\bigcirc$ \hfill}
\nointerlineskip
\vskip10pt}\hfill}
\hskip5pt plus2pt minus2pt
\hbox to30pt{\hfill \vbox to30pt{\vfill $\cdots$ \vfill}\hfill}
\hskip5pt plus2pt minus2pt
\hbox to15pt{\hfill \vbox to30pt{
\hbox to15pt{\hfill $2$ \hfill}
\nointerlineskip\vfill
\hbox to15pt{\hfill $\bigcirc$ \hfill}
\nointerlineskip
\vskip10pt}\hfill}
\hskip0pt plus2pt minus2pt
\hbox to30pt{\hfill
\vbox to30pt{\vfill \nointerlineskip
\hbox to25pt{\rightarrowfill} \nointerlineskip
\hbox to25pt{\leftarrowfill} \nointerlineskip
\vfill}
\hfill}
\hskip0pt plus2pt minus2pt
\hbox to15pt{\hfill \vbox to30pt{
\hbox to15pt{\hfill $2$ \hfill}
\nointerlineskip\vfill
\hbox to15pt{\hfill $\bigcirc$ \hfill}
\nointerlineskip
\vskip10pt}\hfill}
\hskip0pt plus2pt minus2pt
\hbox to30pt{\hfill
\vbox to30pt{\vfill
\hbox to25pt{\hfill $2$ \hfill}
\nointerlineskip
\hbox to25pt{\rightarrowfill} \nointerlineskip
\hbox to25pt{\leftarrowfill} \nointerlineskip
\vskip8pt}
\hfill}
\hskip0pt plus2pt minus2pt
\hbox to15pt{\hfill \vbox to30pt{
\hbox to15pt{\hfill $1$ \hfill}
\nointerlineskip\vfill
\hbox to15pt{\hfill $\bigcirc$ \hfill}
\nointerlineskip
\vskip10pt}\hfill}
\hfil
}
\hskip0pt}


\centerline{\hfil
\hbox to 300pt{
\hfil
\hbox to15pt{\hfill \vbox
to30pt{\vfill $\widetilde{{\bold B}{\bold C}}_n:$ \vfill}\hfill}
\hskip25pt plus2pt minus2pt
\hbox to15pt{\hfill \vbox to30pt{
\hbox to15pt{\hfill $1$ \hfill}
\nointerlineskip\vfill
\hbox to15pt{\hfill $\bigcirc$ \hfill}
\nointerlineskip
\vskip10pt}\hfill}
\hskip5pt plus2pt minus2pt
\hbox to30pt{\hfill
\vbox to30pt{\vskip8pt \nointerlineskip
\hbox to25pt{\rightarrowfill} \nointerlineskip
\hbox to25pt{\leftarrowfill} \nointerlineskip
\hbox to25pt{\hfill $2$ \hfill}
\vfill}
\hfill}
\hskip0pt plus2pt minus2pt
\hbox to15pt{\hfill \vbox to30pt{
\hbox to15pt{\hfill $2$ \hfill}
\nointerlineskip\vfill
\hbox to15pt{\hfill $\bigcirc$ \hfill}
\nointerlineskip
\vskip10pt}\hfill}
\hskip0pt plus2pt minus2pt
\hbox to30pt{\hfill
\vbox to30pt{\vfill \nointerlineskip
\hbox to25pt{\rightarrowfill} \nointerlineskip
\hbox to25pt{\leftarrowfill} \nointerlineskip
\vfill}
\hfill}
\hskip0pt plus2pt minus2pt
\hbox to15pt{\hfill \vbox to30pt{
\hbox to15pt{\hfill $2$ \hfill}
\nointerlineskip\vfill
\hbox to15pt{\hfill $\bigcirc$ \hfill}
\nointerlineskip
\vskip10pt}\hfill}
\hskip0pt plus2pt minus2pt
\hbox to30pt{\hfill \vbox to30pt{\vfill $\cdots$ \vfill}\hfill}
\hskip0pt plus2pt minus2pt
\hbox to15pt{\hfill \vbox to30pt{
\hbox to15pt{\hfill $2$ \hfill}
\nointerlineskip\vfill
\hbox to15pt{\hfill $\bigcirc$ \hfill}
\nointerlineskip
\vskip10pt}\hfill}
\hskip0pt plus2pt minus2pt
\hbox to30pt{\hfill
\vbox to30pt{\vfill \nointerlineskip
\hbox to25pt{\rightarrowfill} \nointerlineskip
\hbox to25pt{\leftarrowfill} \nointerlineskip
\vfill}
\hfill}
\hskip0pt plus2pt minus2pt
\hbox to15pt{\hfill \vbox to30pt{
\hbox to15pt{\hfill $2$ \hfill}
\nointerlineskip\vfill
\hbox to15pt{\hfill $\bigcirc$ \hfill}
\nointerlineskip
\vskip10pt}\hfill}
\hskip0pt plus2pt minus2pt
\hbox to30pt{\hfill
\vbox to30pt{\vskip8pt \nointerlineskip
\hbox to25pt{\rightarrowfill} \nointerlineskip
\hbox to25pt{\leftarrowfill} \nointerlineskip
\hbox to25pt{\hfill $2$ \hfill}
\vfill}
\hfill}
\hskip0pt plus2pt minus2pt
\hbox to15pt{\hfill \vbox to30pt{
\hbox to15pt{\hfill $2$ \hfill}
\nointerlineskip\vfill
\hbox to15pt{\hfill $\bigcirc$ \hfill}
\nointerlineskip
\vskip10pt}\hfill}
\hfil
}
\hskip0pt}


\centerline{\hfil
\hbox to 300pt{
\hfil
\hbox to15pt{\hfill \vbox
to30pt{\vfill $\tilde {\bold C}_n\ :$ \vfill}\hfill}
\hskip20pt plus2pt minus2pt
\hbox to15pt{\hfill \vbox to30pt{
\hbox to15pt{\hfill $1$ \hfill}
\nointerlineskip\vfill
\hbox to15pt{\hfill $\bigcirc$ \hfill}
\nointerlineskip
\vskip10pt}\hfill}
\hskip5pt plus2pt minus2pt
\hbox to30pt{\hfill
\vbox to30pt{\nointerlineskip
\hbox to25pt{\hfill $2$ \hfill}
\nointerlineskip\vfill
\hbox to25pt{\rightarrowfill} \nointerlineskip
\hbox to25pt{\leftarrowfill} \nointerlineskip
\vskip8pt}
\hfill}
\hskip0pt plus2pt minus2pt
\hbox to15pt{\hfill \vbox to30pt{
\hbox to15pt{\hfill $1$ \hfill}
\nointerlineskip\vfill
\hbox to15pt{\hfill $\bigcirc$ \hfill}
\nointerlineskip
\vskip10pt}\hfill}
\hskip0pt plus2pt minus2pt
\hbox to30pt{\hfill
\vbox to30pt{\vfill \nointerlineskip
\hbox to25pt{\rightarrowfill} \nointerlineskip
\hbox to25pt{\leftarrowfill} \nointerlineskip
\vfill}
\hfill}
\hskip0pt plus2pt minus2pt
\hbox to15pt{\hfill \vbox to30pt{
\hbox to15pt{\hfill $1$ \hfill}
\nointerlineskip\vfill
\hbox to15pt{\hfill $\bigcirc$ \hfill}
\nointerlineskip
\vskip10pt}\hfill}
\hskip5pt plus2pt minus2pt
\hbox to30pt{\hfill \vbox to30pt{\vfill $\cdots$ \vfill}\hfill}
\hskip5pt plus2pt minus2pt
\hbox to15pt{\hfill \vbox to30pt{
\hbox to15pt{\hfill $1$ \hfill}
\nointerlineskip\vfill
\hbox to15pt{\hfill $\bigcirc$ \hfill}
\nointerlineskip
\vskip10pt}\hfill}
\hskip0pt plus2pt minus2pt
\hbox to30pt{\hfill
\vbox to30pt{\vfill \nointerlineskip
\hbox to25pt{\rightarrowfill} \nointerlineskip
\hbox to25pt{\leftarrowfill} \nointerlineskip
\vfill}
\hfill}
\hskip0pt plus2pt minus2pt
\hbox to15pt{\hfill \vbox to30pt{
\hbox to15pt{\hfill $1$ \hfill}
\nointerlineskip\vfill
\hbox to15pt{\hfill $\bigcirc$ \hfill}
\nointerlineskip
\vskip10pt}\hfill}
\hskip0pt plus2pt minus2pt
\hbox to30pt{\hfill
\vbox to30pt{\vskip8pt
\nointerlineskip
\hbox to25pt{\rightarrowfill} \nointerlineskip
\hbox to25pt{\leftarrowfill} \nointerlineskip\vfill
\hbox to25pt{\hfill $2$ \hfill}
}
\hfill}
\hskip0pt plus2pt minus2pt
\hbox to15pt{\hfill \vbox to30pt{
\hbox to15pt{\hfill $1$ \hfill}
\nointerlineskip\vfill
\hbox to15pt{\hfill $\bigcirc$ \hfill}
\nointerlineskip
\vskip10pt}\hfill}
\hfil
}
\hskip0pt}

\centerline{
\hbox to15pt{\hfill \vbox
to60pt{\vfill $\widetilde{{\bold B}{\bold D}}_n:$ \vfill}\hfill}
\hskip10pt
\hbox to10pt{\hfill
\vbox to60pt{$1$\vfill$2$\vfill$1$}\hfill}
\hskip5pt
\hbox to20pt{\hfill
\vbox to 60pt {\vfill
\hbox to 15pt{\hfill $\bigcirc$ \hfill}
\nointerlineskip\vfill
\hbox to15pt{\hfill $\uparrow\downarrow$ \hfill}
\nointerlineskip\vfill
\hbox to15pt{\hfill $\bigcirc$ \hfill}
\nointerlineskip\vfill
\hbox to15pt{\hfill $\uparrow\downarrow$ \hfill}
\nointerlineskip\vfill
\hbox to15pt{\hfill $\bigcirc$ \hfill}
\vfill}
\hfill}
\hskip-10pt
\raise15pt\hbox{\darr
\hskip-5pt
\vertex{2}
\hskip-20pt plus2pt minus2pt
\hbox to20pt{\hfill \vbox to30pt{\vfill $\cdots$ \vfill}\hfill}
\hskip0pt plus2pt minus2pt
\vertex{2}
\hskip-5pt
\darrup{2}
\hskip-5pt
\vertex{1}}
\hskip-10pt
\hbox to15pt{\hfill \vbox
to60pt{\vfill $\widetilde{{\bold C}{\bold D}}_n:$ \vfill}\hfill}
\hskip10pt
\hbox to10pt{\hfill
\vbox to60pt{$1$\vfill$2$\vfill$1$}\hfill}
\hskip5pt
\hbox to20pt{\hfill
\vbox to 60pt {\vfill
\hbox to 15pt{\hfill $\bigcirc$ \hfill}
\nointerlineskip\vfill
\hbox to15pt{\hfill $\uparrow\downarrow$ \hfill}
\nointerlineskip\vfill
\hbox to15pt{\hfill $\bigcirc$ \hfill}
\nointerlineskip\vfill
\hbox to15pt{\hfill $\uparrow\downarrow$ \hfill}
\nointerlineskip\vfill
\hbox to15pt{\hfill $\bigcirc$ \hfill}
\vfill}
\hfill}
\hskip-10pt
\raise15pt\hbox{\darr
\hskip-5pt
\vertex{2}
\hskip-20pt plus2pt minus2pt
\hbox to20pt{\hfill \vbox to30pt{\vfill $\cdots$ \vfill}\hfill}
\hskip0pt plus2pt minus2pt
\vertex{2}
\hskip-5pt
\darrdown{2}
\hskip-5pt
\vertex{2}}
}

\centerline{
\hbox to15pt{\hfill \vbox
to60pt{\vfill $\tilde{\bold D}_n:$ \vfill}\hfill}
\hskip5pt
\hbox to10pt{\hfill
\vbox to60pt{$1$\vfill$2$\vfill$1$}\hfill}
\hskip5pt
\hbox to20pt{\hfill
\vbox to 60pt {\vfill
\hbox to 15pt{\hfill $\bigcirc$ \hfill}
\nointerlineskip\vfill
\hbox to15pt{\hfill $\uparrow\downarrow$ \hfill}
\nointerlineskip\vfill
\hbox to15pt{\hfill $\bigcirc$ \hfill}
\nointerlineskip\vfill
\hbox to15pt{\hfill $\uparrow\downarrow$ \hfill}
\nointerlineskip\vfill
\hbox to15pt{\hfill $\bigcirc$ \hfill}
\vfill}
\hfill}
\hskip-10pt
\raise15pt\hbox{\darr
\hskip-5pt
\vertex{2}
\hskip-20pt plus2pt minus2pt
\hbox to20pt{\hfill \vbox to30pt{\vfill $\cdots$ \vfill}\hfill}
\hskip0pt plus2pt minus2pt
\vertex{2}
\hskip-5pt
\darr}
\hskip-5pt
\hbox to20pt{\hfill
\vbox to 60pt {\vfill
\hbox to 15pt{\hfill $\bigcirc$ \hfill}
\nointerlineskip\vfill
\hbox to15pt{\hfill $\uparrow\downarrow$ \hfill}
\nointerlineskip\vfill
\hbox to15pt{\hfill $\bigcirc$ \hfill}
\nointerlineskip\vfill
\hbox to15pt{\hfill $\uparrow\downarrow$ \hfill}
\nointerlineskip\vfill
\hbox to15pt{\hfill $\bigcirc$ \hfill}
\vfill}
\hfill}
\hskip-17pt
\hbox to10pt{\hfill
\vbox to60pt{$1$\vfill$2$\vfill$1$}\hfill}
\hskip10pt
\hbox to15pt{\hfill \vbox
to60pt{\vfill $\tilde{\bold E}_6:$ \vfill}\hfill}
\raise30pt\hbox{
\hskip5pt
\vertex{1}
\hskip-5pt
\darr
\hskip-5pt
\vertex{2}
\hskip-5pt
\darr
}
\hskip-10pt
\hbox to20pt{\hfill
\vbox to 60pt {\vfill
\hbox to15pt{\hfill $3$ \hfill}
\nointerlineskip\vfill
\hbox to 15pt{\hfill $\bigcirc$ \hfill}
\nointerlineskip\vfill
\hbox to15pt{\hfill $\uparrow\downarrow$ \hfill}
\nointerlineskip\vfill
\hbox to15pt{\hfill $\bigcirc$ \hfill}
\nointerlineskip\vfill
\hbox to15pt{\hfill $\uparrow\downarrow$ \hfill}
\nointerlineskip\vfill
\hbox to15pt{\hfill $\bigcirc$ \hfill}
\vfill}
\hfill}
\hskip-20pt
\hbox to10pt{\hfill
\vbox to60pt{\vskip30pt $2$\vfill$1$\vskip3pt}\hfill}
\hskip-10pt
\raise30pt\hbox{
\darr
\hskip-5pt
\vertex{2}
\hskip-5pt
\darr
\hskip-5pt
\vertex{1}}
}

\vskip-20pt

\centerline{
\lower10pt\hbox to15pt{\hfill \vbox
to60pt{\vfill $\tilde{\bold E}_7:$ \vfill}\hfill}
\hskip10pt
\vertex{1}
\hskip-5pt
\darr
\hskip-5pt
\vertex{2}
\hskip-5pt
\darr
\hskip-5pt
\vertex{3}
\hskip-5pt
\darr
\hskip-5pt
\hbox to20pt{\hfill
\vbox to 30pt {\vskip3pt
\hbox to15pt{\hfill $4$ \hfill}
\nointerlineskip\vfill
\hbox to 15pt{\hfill $\bigcirc$ \hfill}
\nointerlineskip\vfill
\hbox to15pt{\hfill $\uparrow\downarrow$ \hfill}
\nointerlineskip\vfill
\hbox to15pt{\hfill $\bigcirc$ \hfill}
\nointerlineskip\vfill
}
\hfill}
\hskip-20pt
\lower5pt\hbox to10pt{\hfill
\vbox to30pt{\vfill $2$}\hfill}
\hskip0pt
\darr
\hskip-5pt
\vertex{3}
\hskip-5pt
\darr
\hskip-5pt
\vertex{2}
\hskip-5pt
\darr
\hskip-5pt
\vertex{1}
}

\centerline{
\lower10pt\hbox to15pt{\hfill \vbox
to60pt{\vfill $\tilde{\bold E}_8:$ \vfill}\hfill}
\hskip10pt
\vertex{2}
\hskip-5pt
\darr
\hskip-5pt
\vertex{4}
\hskip-5pt
\darr
\hskip-5pt
\hbox to20pt{\hfill
\vbox to 30pt {\vskip3pt
\hbox to15pt{\hfill $6$ \hfill}
\nointerlineskip\vfill
\hbox to 15pt{\hfill $\bigcirc$ \hfill}
\nointerlineskip\vfill
\hbox to15pt{\hfill $\uparrow\downarrow$ \hfill}
\nointerlineskip\vfill
\hbox to15pt{\hfill $\bigcirc$ \hfill}
\nointerlineskip\vfill
}
\hfill}
\hskip-20pt
\lower5pt\hbox to10pt{\hfill
\vbox to30pt{\vfill $3$}\hfill}
\hskip0pt
\darr
\hskip-5pt
\vertex{5}
\hskip-5pt
\darr
\hskip-5pt
\vertex{4}
\hskip-5pt
\darr
\hskip-5pt
\vertex{3}
\hskip-5pt
\darr
\hskip-5pt
\vertex{2}
\hskip-5pt
\darr
\hskip-5pt
\vertex{1}
}

\centerline{\hfill
\hbox to15pt{\hfill \vbox
to30pt{\vfill $\widetilde{{\bold B}{\bold F}}_4:$ \vfill}\hfill}
\hskip20pt plus 2pt minus 2pt
\vertex{1}
\hskip0pt plus 2pt minus 2pt
\darr
\hskip0pt plus 2pt minus 2pt
\vertex{2}
\darr
\hskip0pt plus 2pt minus 2pt
\vertex{3}
\hskip0pt plus 2pt minus 2pt
\darrup{2}
\hskip0pt plus 2pt minus 2pt
\vertex{2}
\hskip0pt plus 2pt minus 2pt
\darr
\hskip0pt plus 2pt minus 2pt
\vertex{1}
\hskip70pt}


\centerline{\hfill
\hbox to15pt{\hfill \vbox
to30pt{\vfill $\widetilde{{\bold C}{\bold F}}_4:$ \vfill}\hfill}
\hskip20pt plus 2pt minus 2pt
\vertex{2}
\hskip0pt plus 2pt minus 2pt
\darr
\hskip0pt plus 2pt minus 2pt
\vertex{4}
\darrup{2}
\hskip0pt plus 2pt minus 2pt
\vertex{3}
\hskip0pt plus 2pt minus 2pt
\darr
\hskip0pt plus 2pt minus 2pt
\vertex{2}
\hskip0pt plus 2pt minus 2pt
\darr
\hskip0pt plus 2pt minus 2pt
\vertex{1}
\hskip70pt}


\centerline{\hfill
\hbox to15pt{\hfill \vbox
to30pt{\vfill $\widetilde{{\bold A}{\bold G}}_2:$ \vfill}\hfill}
\hskip20pt plus 2pt minus 2pt
\vertex{1}
\hskip-5pt plus 2pt minus 2pt
\darr
\hskip-5pt plus 2pt minus 2pt
\vertex{2}
\hskip-5pt plus 2pt minus 2pt
\darrup{3}
\hskip-5pt plus 2pt minus 2pt
\vertex{1}
\hskip0pt plus 2pt minus 2pt
\hbox to15pt{\hfill \vbox
to30pt{\vfill $\widetilde{{\bold G}{\bold A}}_2:$ \vfill}\hfill}
\hskip20pt plus 2pt minus 2pt
\vertex{1}
\hskip-5pt plus 2pt minus 2pt
\darr
\hskip-5pt plus 2pt minus 2pt
\vertex{2}
\hskip-5pt plus2pt minus2pt
\darrdown{3}
\hskip-5pt plus 2pt minus 2pt
\vertex{3}
\hskip0pt plus 2pt minus 2pt
\hfill}


\centerline{{\bf Table 4.2.} Calabi--Yau parabolic diagrams without single
arrows}
\centerline{(classical extended Dynkin diagrams)}

\newpage

\centerline{\hbox{
$\bigcirc$
\hskip-5pt
\lower13pt\darrupdown{t_{12}}{t_{21}}
\hskip-5pt
$\bigcirc$} \hskip20pt where $t_{12}t_{21}>4$.
\hfil}

\vskip15pt


\centerline{\hbox{
$\bigcirc$
\hskip0pt
\lower13pt\darrupdown{t_{12}}{t_{21}}
\hskip-5pt
$\bigcirc$
\hskip0pt
\lower13pt\darrupdown{t_{23}}{t_{32}}
\hskip-5pt
$\bigcirc$}
\hskip20pt where $0<t_{12}t_{21}<4,\ 0<t_{23}t_{32}<4,\
t_{12}t_{21}+t_{23}t_{32}>4$ \hfil}

\vskip20pt

\centerline{\hbox{
\vbox to70pt{
\hbox to100pt{\hfill$\bigcirc$\hfill}
\vskip10pt\nointerlineskip\vfil
\hbox to100pt{\hfil
$t_{12}$ $\nearrow$\hskip-5pt$\swarrow$ $t_{21}$ \hfill
$t_{32}$ $\nwarrow$\hskip-4pt$\searrow$ $t_{23}$\hfil}
\nointerlineskip\vfil
\hbox to100pt{\hskip10pt \raise12pt\hbox{$\bigcirc$}\hfil
\darrupdown{t_{13}}{t_{31}} \hfil \raise12pt\hbox{$\bigcirc$}\hskip10pt}
}
}
\hskip20pt \raise30pt\hbox{\vbox{
\hbox{where\hfil}
\hbox{\hfil $0<t_{12}t_{21}<4,\ 0<t_{23}t_{32}<4,\
0<t_{31}t_{13}<4$,\hfil}
\hbox{$t_{12}t_{21}+t_{23}t_{32}+t_{31}t_{13}>3$}
}}
}

\vskip20pt


\centerline{
\hbox to15pt{\hfill
\vbox to 40pt{
\hbox to15pt{\hfill $\bigcirc$ \hfill}
\nointerlineskip\vfill
\hbox to15pt{\hfill $\uparrow\downarrow$ \hfill}
\nointerlineskip\vfill
\hbox to15pt{\hfill $\bigcirc$ \hfill}
\nointerlineskip}\hfill}
\hskip-5pt
\lower1pt\hbox{\hfill
\vbox to52pt {\darrup{2}
\vskip2pt
\darr}
\hskip-5pt
\hbox to15pt{\hfill
\vbox to 40pt{
\hbox to15pt{\hfill $\bigcirc$ \hfill}
\nointerlineskip\vfill
\hbox to15pt{\hfill $\uparrow\downarrow$ \hfill}
\nointerlineskip\vfill
\hbox to15pt{\hfill $\bigcirc$ \hfill}
\nointerlineskip}\hfill}
}
\hskip5pt
\hbox to15pt{\hfill
\vbox to 40pt{
\hbox to15pt{\hfill $\bigcirc$ \hfill}
\nointerlineskip\vfill
\hbox to15pt{\hfill $\uparrow\downarrow$ \hfill}
\nointerlineskip\vfill
\hbox to15pt{\hfill $\bigcirc$ \hfill}
\nointerlineskip}\hfill}
\hskip-5pt
\lower1pt\hbox{\hfill
\vbox to52pt {\darrup{2}
\vskip2pt
\darrup{2}}
\hskip-5pt
\hbox to15pt{\hfill
\vbox to 40pt{
\hbox to15pt{\hfill $\bigcirc$ \hfill}
\nointerlineskip\vfill
\hbox to15pt{\hfill $\uparrow\downarrow$ \hfill}
\nointerlineskip\vfill
\hbox to15pt{\hfill $\bigcirc$ \hfill}
\nointerlineskip}\hfill}
}
\hskip5pt
\hbox to15pt{\hfill
\vbox to 40pt{
\hbox to15pt{\hfill $\bigcirc$ \hfill}
\nointerlineskip\vfill
\hbox to15pt{\hfill $\uparrow\downarrow$ \hfill}
\nointerlineskip\vfill
\hbox to15pt{\hfill $\bigcirc$ \hfill}
\nointerlineskip}\hfill}
\hskip-5pt
\lower1pt\hbox{\hfill
\vbox to52pt {\darrup{2}
\vskip2pt
\darrdown{2}}
\hskip-5pt
\hbox to15pt{\hfill
\vbox to 40pt{
\hbox to15pt{\hfill $\bigcirc$ \hfill}
\nointerlineskip\vfill
\hbox to15pt{\hfill $\uparrow\downarrow$ \hfill}
\nointerlineskip\vfill
\hbox to15pt{\hfill $\bigcirc$ \hfill}
\nointerlineskip}\hfill}
}
\hfil}


\vskip20pt

\centerline{
\hbox to15pt{\hfill
\vbox to 40pt{
\hbox to15pt{\hfill $\bigcirc$ \hfill}
\nointerlineskip\vfill
\hbox to15pt{\hfill $\uparrow\downarrow$ \hfill}
\nointerlineskip\vfill
\hbox to15pt{\hfill $\bigcirc$ \hfill}
\nointerlineskip}\hfill}
\hskip-5pt
\lower3pt\hbox to65pt{\hfill
\vbox to50pt {
\hbox to60pt{\hfill
\vbox to12pt{
\hbox to60pt{\hfill $2$ \hfill}
\nointerlineskip
\hbox to60pt{\rightarrowfill} \nointerlineskip\vskip-2pt
\hbox to60pt{\leftarrowfill}
\vfill}\hfill}
\vfill
\hbox to60pt{\hfill
\hbox to23pt{\hfill
\vbox to12pt{\vfill
\hbox to20pt{\rightarrowfill} \nointerlineskip\vskip-2pt
\hbox to20pt{\leftarrowfill}
\vfill}\hfill}
\hskip0pt \raise4pt\hbox{$\bigcirc$} \hskip0pt
\hbox to20pt{\hfill
\vbox to12pt{\vfill
\hbox to20pt{\rightarrowfill} \nointerlineskip\vskip-2pt
\hbox to20pt{\leftarrowfill}
\vfill}\hfill}
\hfill}
}
\hfill}
\hskip-5pt
\hbox to15pt{\hfill
\vbox to 40pt{
\hbox to15pt{\hfill $\bigcirc$ \hfill}
\nointerlineskip\vfill
\hbox to15pt{\hfill $\uparrow\downarrow$ \hfill}
\nointerlineskip\vfill
\hbox to15pt{\hfill $\bigcirc$ \hfill}
\nointerlineskip}\hfill}
\hfil}

\vskip20pt

\centerline{{\bf Table 4.3.} Calabi--Yau Lanner
diagrams without single arrows}
\centerline{(classical Lanner diagrams).}

\newpage

\centerline{\hbox to110pt{\hfil
\hbox to30pt{\hfil \vbox to70pt{\vfil $t_{21}$ \vfil}\hfil}
\hskip-13pt
\hbox to15pt{\hfil \vbox to70pt{\hbox{\hfil $\bigcirc$ \hfil}
\vfil
\hbox to10pt{\hskip2pt $\downarrow$ \hfil}
\vfil
\hbox{\hfil $\bigcirc$ \hfil}
}}
\hskip-10pt
\hbox to60pt{\hfil \vbox to60pt{
\vskip5pt
\hbox{\hfil \lower4pt\hbox{$t_{32}$} $\nwarrow$ \hskip-4pt $\searrow$
\raise4pt\hbox{$t_{23}$}\hfil}
\vfil
\hbox{\hfil \raise4pt\hbox{$t_{13}$} $\nearrow$ \hskip-5pt $\swarrow$
\lower4pt\hbox{$t_{31}$} \hfill}
\vskip5pt
}\hfil}
\hskip-30pt
\hbox to10pt{\hfil\vbox to60pt{\vfil $\bigcirc$ \vfil}\hfil}
\hfil}
\hskip10pt
\hbox to170pt{\hfil \vbox{
\hbox{Elliptic diagram where \hfil}
\hbox{$t_{21}t_{13}t_{32}+2t_{13}t_{31}+2t_{23}t_{32}<8$;}
\vskip3pt
\hbox{Parabolic diagram where \hfil}
\hbox{$t_{21}t_{13}t_{32}+2t_{13}t_{31}+2t_{23}t_{32}=8$;}
\vskip3pt
\hbox{Quasi-Lanner diagram where}
\hbox{$t_{13}t_{31}\le 4$, $t_{23}t_{32}\le 4$}
\hbox{and $t_{21}t_{13}t_{32}+2t_{13}t_{31}+2t_{23}t_{32}>8.$ }
}\hfil}
\hfil}

\vskip10pt

\centerline{\hskip40pt\hbox to90pt{\hfil
\hbox to15pt{\hfil \vbox to30pt{\hbox{\hfil $\bigcirc$ \hfil}
\vfil
\hbox to10pt{\hskip2pt $\downarrow$ \hfil}
\vfil
\hbox{\hfil $\bigcirc$ \hfil}
}}
\hskip-25pt
\hbox to60pt{\hfil \vbox to30pt{
\vskip3pt
\hbox{\hfil $\nwarrow$ \hskip-4pt $\searrow$ \hfil}
\vfil
\hbox{\hfil $\nearrow$ \hskip-5pt $\swarrow$ \hfill}
\vskip3pt
}\hfil}
\hskip-37pt
\hbox to10pt{\hfil\vbox to30pt{\vfil $\bigcirc$ \vfil}\hfil}
\hfil}
\hskip-20pt
\darr
\hskip-17pt
\hbox to10pt{\hfil\vbox to30pt{\vfil $\bigcirc$ \vfil}\hfil}
\hskip 20pt
\hbox to130pt{\hfil {\vbox to30pt{
\vfil
\hbox{Elliptic diagram.}
\vfil}
\hfil}
}
\hskip50pt}

\vskip10pt

\centerline{\hbox to90pt{\hfil
\hbox to30pt{\hfil \vbox to30pt{\vfil $2$ \vfil}\hfil}
\hskip-17pt
\hbox to15pt{\hfil \vbox to30pt{\hbox{\hfil $\bigcirc$ \hfil}
\vfil
\hbox to10pt{\hskip2pt $\downarrow$ \hfil}
\vfil
\hbox{\hfil $\bigcirc$ \hfil}
}}
\hskip-25pt
\hbox to60pt{\hfil \vbox to30pt{
\vskip3pt
\hbox{\hfil $\nwarrow$ \hskip-4pt $\searrow$ \hfil}
\vfil
\hbox{\hfil $\nearrow$ \hskip-5pt $\swarrow$ \hfill}
\vskip3pt
}\hfil}
\hskip-37pt
\hbox to10pt{\hfil\vbox to30pt{\vfil $\bigcirc$ \vfil}\hfil}
\hfil}
\hskip-17pt
\darr
\hskip-17pt
\hbox to10pt{\hfil\vbox to30pt{\vfil $\bigcirc$ \vfil}\hfil}
\hskip 20pt
\hbox to120pt{\hfil {\vbox to30pt{
\vfil
\hbox{Parabolic diagram.}
\vfil}
\hfil}
}
\hskip50pt}

\vskip5pt

\centerline{\hbox to130pt{\hfil
\hbox to30pt{\hfil \vbox to70pt{\vfil $t_{21}$ \vfil}\hfil}
\hskip-13pt
\hbox to15pt{\hfil \vbox to70pt{\hbox{\hfil $\bigcirc$ \hfil}
\vfil
\hbox to10pt{\hskip2pt $\downarrow$ \hfil}
\vfil
\hbox{\hfil $\bigcirc$ \hfil}
}}
\hskip-10pt
\hbox to60pt{\hfil \vbox to60pt{
\vskip5pt
\hbox{\hfil \lower4pt\hbox{$t_{32}$} $\nwarrow$ \hskip-4pt $\searrow$
\raise4pt\hbox{$t_{23}$}\hfil}
\vfil
\hbox{\hfil \raise4pt\hbox{$t_{13}$} $\nearrow$ \hskip-5pt $\swarrow$
\lower4pt\hbox{$t_{31}$} \hfill}
\vskip5pt
}\hfil}
\hskip-30pt
\hbox to10pt{\hfil\vbox to60pt{\vfil $\bigcirc$ \vfil}\hfil}
\hskip10pt
\raise15pt\darrupdown{t_{34}}{t_{43}}
\hskip-18pt
\hbox to10pt{\hfil\vbox to60pt{\vfil $\bigcirc$ \vfil}\hfil}
\hfil}
\hskip5pt
\raise5pt\hbox to190pt{\hfil\vbox{\vfil
\hbox{Quasi-Lanner diagram where}
\hbox{$t_{13}t_{31}+t_{34}t_{43}\le 4$, $t_{23}t_{32}+t_{34}t_{43}\le 4$,}
\hbox{$t_{21}t_{13}t_{32}+2t_{13}t_{31}+2t_{23}t_{32}\le 8$,}
\hbox{$t_{21}t_{13}t_{32}+2t_{13}t_{31}+2t_{23}t_{32}+2t_{34}t_{43}>8.$ }
\vfil}\hfil}
}

\vskip10pt

\centerline{\hbox to90pt{\hfil
\hbox to15pt{\hfil \vbox to30pt{\hbox{\hfil $\bigcirc$ \hfil}
\vfil
\hbox to10pt{\hskip2pt $\downarrow$ \hfil}
\vfil
\hbox{\hfil $\bigcirc$ \hfil}
}}
\hskip-25pt
\hbox to60pt{\hfil \vbox to30pt{
\vskip3pt
\hbox{\hfil $\nwarrow$ \hskip-4pt $\searrow$ \hfil}
\vfil
\hbox{\hfil $\nearrow$ \hskip-5pt $\swarrow$ \hfill}
\vskip3pt
}\hfil}
\hskip-37pt
\hbox to10pt{\hfil\vbox to30pt{\vfil $\bigcirc$ \vfil}\hfil}
\hfil}
\hskip-25pt
\darr
\hskip-17pt
\hbox to10pt{\hfil\vbox to30pt{\vfil $\bigcirc$ \vfil}\hfil}
\hskip7pt
\darr
\hskip-17pt
\hbox to10pt{\hfil\vbox to30pt{\vfil $\bigcirc$ \vfil}\hfil}
\hskip 20pt
\hbox to130pt{\hfil {\vbox to30pt{
\vfil
\hbox{Elliptic diagram}
\vfil}
\hfil}
}
\hskip50pt}

\vskip10pt

\centerline{\hbox to90pt{\hfil
\hbox to15pt{\hfil \vbox to30pt{\hbox{\hfil $\bigcirc$ \hfil}
\vfil
\hbox to10pt{\hskip2pt $\downarrow$ \hfil}
\vfil
\hbox{\hfil $\bigcirc$ \hfil}
}}
\hskip-25pt
\hbox to60pt{\hfil \vbox to30pt{
\vskip3pt
\hbox{\hfil $\nwarrow$ \hskip-4pt $\searrow$ \hfil}
\vfil
\hbox{\hfil $\nearrow$ \hskip-5pt $\swarrow$ \hfill}
\vskip3pt
}\hfil}
\hskip-37pt
\hbox to10pt{\hfil\vbox to30pt{\vfil $\bigcirc$ \vfil}\hfil}
\hfil}
\hskip-25pt
\darr
\hskip-17pt
\hbox to10pt{\hfil\vbox to30pt{\vfil $\bigcirc$ \vfil}\hfil}
\hskip7pt
\darrupdown{t_{45}}{t_{54}}
\hskip-17pt
\hbox to10pt{\hfil\vbox to30pt{\vfil $\bigcirc$ \vfil}\hfil}
\hskip 20pt
\hbox to140pt{\hfil {\vbox to30pt{
\vfil
\hbox{Lanner diagram}
\hbox{where $t_{45}t_{54}=2$.}
\vfil}
\hfil}
}
\hskip30pt}

\vskip10pt

\centerline{\hbox to90pt{\hfil
\hbox to30pt{\hfil \vbox to30pt{\vfil $2$ \vfil}\hfil}
\hskip-17pt
\hbox to15pt{\hfil \vbox to30pt{\hbox{\hfil $\bigcirc$ \hfil}
\vfil
\hbox to10pt{\hskip2pt $\downarrow$ \hfil}
\vfil
\hbox{\hfil $\bigcirc$ \hfil}
}}
\hskip-25pt
\hbox to60pt{\hfil \vbox to30pt{
\vskip3pt
\hbox{\hfil $\nwarrow$ \hskip-4pt $\searrow$ \hfil}
\vfil
\hbox{\hfil $\nearrow$ \hskip-5pt $\swarrow$ \hfill}
\vskip3pt
}\hfil}
\hskip-37pt
\hbox to10pt{\hfil\vbox to30pt{\vfil $\bigcirc$ \vfil}\hfil}
\hfil}
\hskip-17pt
\darr
\hskip-17pt
\hbox to10pt{\hfil\vbox to30pt{\vfil $\bigcirc$ \vfil}\hfil}
\hskip7pt
\darrupdown{t_{45}}{t_{54}}
\hskip-17pt
\hbox to10pt{\hfil\vbox to30pt{\vfil $\bigcirc$ \vfil}\hfil}
\hskip 20pt
\hbox to130pt{\hfil {\vbox to30pt{
\vfil
\hbox{Quasi-Lanner diagram}
\hbox{where $1\le t_{45}t_{54}\le 2$.}
\vfil}
\hfil}
}
\hskip50pt}

\vskip10pt

\centerline{\hbox to90pt{\hfil
\hbox to15pt{\hfil \vbox to30pt{\hbox{\hfil $\bigcirc$ \hfil}
\vfil
\hbox to10pt{\hskip2pt $\downarrow$ \hfil}
\vfil
\hbox{\hfil $\bigcirc$ \hfil}
}}
\hskip-25pt
\hbox to60pt{\hfil \vbox to30pt{
\vskip3pt
\hbox{\hfil $\nwarrow$ \hskip-4pt $\searrow$ \hfil}
\vfil
\hbox{\hfil $\nearrow$ \hskip-5pt $\swarrow$ \hfill}
\vskip3pt
}\hfil}
\hskip-37pt
\hbox to10pt{\hfil\vbox to30pt{\vfil $\bigcirc$ \vfil}\hfil}
\hfil}
\hskip-25pt
\darr
\hskip-17pt
\hbox to10pt{\hfil\vbox to30pt{\vfil $\bigcirc$ \vfil}\hfil}
\hskip7pt
\darr
\hskip-17pt
\hbox to10pt{\hfil\vbox to30pt{\vfil $\bigcirc$ \vfil}\hfil}
\hskip7pt
\darr
\hskip-17pt
\hbox to10pt{\hfil\vbox to30pt{\vfil $\bigcirc$ \vfil}\hfil}
\hskip 20pt
\hbox to130pt{\hfil {\vbox to30pt{
\vfil
\hbox{Elliptic diagram.}
\vfil}
\hfil}
}
\hskip50pt}

\vskip10pt

\vskip5pt

\centerline{\hbox to90pt{\hfil
\hbox to15pt{\hfil \vbox to30pt{\hbox{\hfil $\bigcirc$ \hfil}
\vfil
\hbox to10pt{\hskip2pt $\downarrow$ \hfil}
\vfil
\hbox{\hfil $\bigcirc$ \hfil}
}}
\hskip-25pt
\hbox to60pt{\hfil \vbox to30pt{
\vskip3pt
\hbox{\hfil $\nwarrow$ \hskip-4pt $\searrow$ \hfil}
\vfil
\hbox{\hfil $\nearrow$ \hskip-5pt $\swarrow$ \hfill}
\vskip3pt
}\hfil}
\hskip-37pt
\hbox to10pt{\hfil\vbox to30pt{\vfil $\bigcirc$ \vfil}\hfil}
\hfil}
\hskip-25pt
\darr
\hskip-17pt
\hbox to10pt{\hfil\vbox to30pt{\vfil $\bigcirc$ \vfil}\hfil}
\hskip7pt
\darr
\hskip-17pt
\hbox to10pt{\hfil\vbox to30pt{\vfil $\bigcirc$ \vfil}\hfil}
\hskip7pt
\darrupdown{t_{45}}{t_{54}}
\hskip-17pt
\hbox to10pt{\hfil\vbox to30pt{\vfil $\bigcirc$ \vfil}\hfil}
\hskip 20pt
\hbox to140pt{\hfil {\vbox to30pt{
\vfil
\hbox{Lanner diagram}
\hbox{where $t_{45}t_{54}=2$.}
\vfil}
\hfil}
}
\hskip50pt}

\vskip10pt

\centerline{\hbox to90pt{\hfil
\hbox to15pt{\hfil \vbox to30pt{\hbox{\hfil $\bigcirc$ \hfil}
\vfil
\hbox to10pt{\hskip2pt $\downarrow$ \hfil}
\vfil
\hbox{\hfil $\bigcirc$ \hfil}
}}
\hskip-25pt
\hbox to60pt{\hfil \vbox to30pt{
\vskip3pt
\hbox{\hfil $\nwarrow$ \hskip-4pt $\searrow$ \hfil}
\vfil
\hbox{\hfil $\nearrow$ \hskip-5pt $\swarrow$ \hfill}
\vskip3pt
}\hfil}
\hskip-37pt
\hbox to10pt{\hfil\vbox to30pt{\vfil $\bigcirc$ \vfil}\hfil}
\hfil}
\hskip-25pt
\darr
\hskip-17pt
\hbox to10pt{\hfil\vbox to30pt{\vfil $\bigcirc$ \vfil}\hfil}
\hskip7pt
\darr
\hskip-17pt
\hbox to10pt{\hfil\vbox to30pt{\vfil $\bigcirc$ \vfil}\hfil}
\hskip7pt
\darr
\hskip-17pt
\hbox to10pt{\hfil\vbox to30pt{\vfil $\bigcirc$ \vfil}\hfil}
\hskip7pt
\darrupdown{t_{45}}{t_{54}}
\hskip-17pt
\hbox to10pt{\hfil\vbox to30pt{\vfil $\bigcirc$ \vfil}\hfil}
\hskip 20pt
\hbox to130pt{\hfil {\vbox to30pt{
\vfil
\hbox{Quasi-Lanner diagram}
\hbox{where $t_{45}t_{54}=2$.}
\vfil}
\hfil}
}
\hskip50pt}

\vskip10pt

\centerline{{\bf Table 4.4.} Calabi-Yau diagrams of the type (C).}

\newpage

\hrule
\vskip2pt

\centerline{
\hbox to50pt{\hfill \vbox to30pt{\vfill
${\bold A}_n^\bullet (k;a,b)\ :$ \vfill}\hfill}
\hskip20pt
\vertexb{-k}
\hskip-5pt
\darrupdown{a}{b}
\hskip-5pt
\verte
\hskip-5pt
\darr
\hskip-5pt
\verte
\hskip0pt
\raise12pt\hbox{$\cdots$}
\hskip0pt
\verte
\hskip-5pt
\darr
\hskip-5pt
\verte
\hskip80pt}

Elliptic diagram iff $ab<k(n+1)/n$.
Parabolic diagram iff $ab=k(n+1)/n$.
Quasi-Lanner diagram iff $kn/(n-1)\ge ab> k(n+1)/n$;
in particular, $n\le k+1$.
Lanner diagram
iff $kn/(n-1)> ab> k(n+1)/n$; in particular, $n\le 2$ for $1\le k\le 3$.

\vskip2pt
\hrule width8cm
\vskip2pt

\centerline{
\hbox to50pt{\hfill \vbox to30pt{\vfill
${\bold B}_n^\bullet (k;a,b)_1\ :$ \vfill}\hfill}
\hskip20pt
\vertexb{-k}
\hskip-5pt
\darrupdown{a}{b}
\hskip-5pt
\verte
\hskip-5pt
\darr
\hskip-5pt
\verte
\hskip0pt
\raise12pt\hbox{$\cdots$}
\hskip0pt
\verte
\hskip-5pt
\darr
\hskip-5pt
\verte
\hskip-5pt
\darrup{2}
\hskip-5pt
\verte
\hskip40pt}

\centerline{
\hbox to50pt{\hfill \vbox to30pt{\vfill
${\bold C}_n^\bullet (k;a,b)_1\ :$ \vfill}\hfill}
\hskip20pt
\vertexb{-k}
\hskip-5pt
\darrupdown{a}{b}
\hskip-5pt
\verte
\hskip-5pt
\darr
\hskip-5pt
\verte
\hskip0pt
\raise12pt\hbox{$\cdots$}
\hskip0pt
\verte
\hskip-5pt
\darr
\hskip-5pt
\verte
\hskip-5pt
\darrdown{2}
\hskip-5pt
\verte
\hskip40pt}

Elliptic diagram iff $ab<k$.
Parabolic diagram iff $ab=k$.
Quasi-Lanner diagram iff $kn/(n-1)\ge ab > k$; in particular, $n\le k+1$.
Lanner diagram iff $kn/(n-1)>ab > k$; in particular, $n < k+1$.

\vskip2pt
\hrule width8cm
\vskip2pt

\centerline{
\hbox to50pt{\hfill \vbox to30pt{\vfill
${\bold B}_n^\bullet (k;a,b)_2\ :$ \vfill}\hfill}
\hskip20pt
\vertexb{-k}
\hskip-5pt
\darrupdown{a}{b}
\hskip-5pt
\verte
\hskip-5pt
\darrdown{2}
\hskip-5pt
\verte
\hskip-5pt
\darr
\hskip-5pt
\verte
\hskip0pt
\raise12pt\hbox{$\cdots$}
\hskip0pt
\verte
\hskip-5pt
\darr
\hskip-5pt
\verte
\hskip40pt}

\centerline{
\hbox to50pt{\hfill \vbox to30pt{\vfill
${\bold C}_n^\bullet (k;a,b)_2\ :$ \vfill}\hfill}
\hskip20pt
\vertexb{-k}
\hskip-5pt
\darrupdown{a}{b}
\hskip-5pt
\verte
\hskip-5pt
\darrup{2}
\hskip-5pt
\verte
\hskip-5pt
\darr
\hskip-5pt
\verte
\hskip0pt
\raise12pt\hbox{$\cdots$}
\hskip0pt
\verte
\hskip-5pt
\darr
\hskip-5pt
\verte
\hskip40pt}

Elliptic diagram iff $ab<2k/n$.
Parabolic diagram iff $ab=2k/n$.
Quasi-Lanner diagram iff $2k/(n-1)\ge ab > 2k/n$; in particular, $n\le 2k+1$.
Lanner diagram iff $2k/(n-1)> ab > 2k/n$;
in particular, $n\le 2$ for $1\le k\le 3$.

\vskip2pt
\hrule width8cm
\vskip2pt

\centerline{
\hbox to50pt{\hfill \vbox to30pt{\vfill
${\bold F}_4^\bullet (k;a,b)_1\ :$ \vfill}\hfill}
\hskip20pt
\vertexb{-k}
\hskip-5pt
\darrupdown{a}{b}
\hskip-5pt
\verte
\hskip-5pt
\darr
\hskip-5pt
\verte
\hskip-5pt
\darrup{2}
\hskip-5pt
\verte
\hskip-5pt
\darr
\hskip-5pt
\verte
\hskip80pt}

\centerline{
\hbox to50pt{\hfill \vbox to30pt{\vfill
${\bold F}_4^\bullet (k;a,b)_2\ :$ \vfill}\hfill}
\hskip20pt
\vertexb{-k}
\hskip-5pt
\darrupdown{a}{b}
\hskip-5pt
\verte
\hskip-5pt
\darr
\hskip-5pt
\verte
\hskip-5pt
\darrdown{2}
\hskip-5pt
\verte
\hskip-5pt
\darr
\hskip-5pt
\verte
\hskip80pt}

Elliptic diagram iff $ab<k/2$.
Parabolic diagram iff $ab=k/2$.
Quasi-Lanner diagram iff $k\ge ab > k/2$. Lanner diagram
iff $k> ab >k/2$.

\vskip2pt
\hrule width8cm
\vskip2pt

\centerline{
\hbox{\hbox to50pt{\hfill \vbox to30pt{\vfill
${\bold G}_2^\bullet (k;a,b)_1 :$ \vfill}\hfill}
\hskip15pt
\vertexb{-k}
\hskip-7pt
\darrupdown{a}{b}
\hskip-5pt
\verte
\hskip-5pt
\darrup{3}
\hskip-5pt
\verte}
\hskip-15pt
\hbox{
\hbox to50pt{\hfill \vbox to30pt{\vfill
${\bold G}_2^\bullet (k;a,b)_2 :$ \vfill}\hfill}
\hskip15pt
\vertexb{-k}
\hskip-7pt
\darrupdown{a}{b}
\hskip-5pt
\verte
\hskip-5pt
\darrdown{3}
\hskip-5pt
\verte}
}

Elliptic diagram iff $ab<k/2$.
Parabolic diagram iff $ab=k/2$.
Quasi-Lanner diagram iff $k\ge ab > k/2$.
Lanner diagram iff $k>ab>k/2$.

\vskip2pt
\hrule width8cm

\vskip10pt

\centerline{{\bf Table 4.5.}  Calabi--Yau diagrams of the type (D)}
\centerline{( $n$ is equal to the number of white vertices).}

\newpage


\centerline{\hbox{
\vbox to70pt{
\hbox to100pt{\hfill$\bigcirc$\hfill}
\vskip10pt\nointerlineskip\vfil
\hbox to100pt{\hfil
$t_{12}$ $\nearrow$    \hfill    $\searrow$ $t_{23}$\hfil}
\nointerlineskip\vfil
\hbox to100pt{\hskip10pt \raise12pt\hbox{$\bigcirc$}\hfil
\darrupdown{t_{13}}{t_{31}} \hfil \raise12pt\hbox{$\bigcirc$}\hskip10pt}
}
}
\hskip20pt \raise10pt\hbox{\vbox{
\hbox{Elliptic where \hfil}
\hbox{$t_{12}t_{23}t_{31}+2t_{13}t_{31}<8$}
\vskip3pt
\hbox{Parabolic where \hfil}
\hbox{$t_{12}t_{23}t_{31}+2t_{13}t_{31}=8$}
\vskip3pt
\hbox{Lanner where  \hfil}
\hbox{$t_{13}t_{31} < 4$ and $t_{12}t_{23}t_{31}+2t_{13}t_{31}>8$}
\vskip3pt
\hbox{Quasi-Lanner where  \hfil}
\hbox{$t_{13}t_{31}\le 4$ and $t_{12}t_{23}t_{31}+2t_{13}t_{31}>8$}
}}
}

\vskip25pt


\centerline{\hbox{
\vbox to70pt{
\hbox to100pt{\hfill$\bigcirc$\hfill}
\vfil\nointerlineskip\vfil
\hbox to100pt{\hfil
$t_{12}$ $\nearrow$    \hfill    $\searrow$ $t_{23}$\hfil}
\nointerlineskip\vfil
\hbox to100pt{\hskip10pt $\bigcirc$ \hfil
\hbox to40pt{\leftarrowfill} \hfil $\bigcirc$ \hskip10pt}
\vskip2pt
\hbox to100pt{\hfil $t_{31}$ \hfil}
}
}
\hskip20pt \hbox{\vbox{
\hbox{Elliptic where \hfil}
\hbox{$t_{12}t_{23}t_{31}<8$}
\vskip3pt
\hbox{Parabolic where \hfil}
\hbox{$t_{12}t_{23}t_{31}=8$}
\vskip3pt
\hbox{Lanner where \hfil $t_{12}t_{23}t_{31}>8$}}
}
}

\vskip30pt

\centerline{{\bf Table 4.6}: Calabi-Yau triangle and
special triangle diagrams. }

\subhead
4.3. Basic results on Calabi-Yau $3$-dimensional manifolds
\endsubhead

To prove our main results about
Calabi-Yau $3$-dimensional manifolds, we
use our results above (especially, Basic Theorems 1.3.2 and 3.2
and Theorem 2.2.6) and the result of V.V. Shokurov about
the length of extremal rays of the type (I) on Calabi-Yau
manifolds (see Appendix of V.V. Shokurov).

\proclaim{Theorem 4.3.1 (by the author and V.V. Shokurov)}
Let $X$ be a $3$-dimensional Calabi-Yau
manifold and $\rho (X) >40$.

Then one of two cases (i) or (ii) below hold:

(i) There exists a small extremal ray on $X$.

(ii) There exists a $nef$ element $h$ such that $h^3=0$
(thus, the $nef$ cone $NEF(X)$ and the cubic intersection hypersurface
${\Cal W}_X$ have a common point; here, we don't claim that $h$
is rational!).
\endproclaim

\demo{Proof} Assume that $X$ does not have a $nef$ element
$h$ with $h^3=0$. Since, by Wilson \cite{W1} and \cite{W2}, the
$NEF(X)$ is locally rational polyhedral
away ${\Cal W}_X$, using compactness arguments, we get that
$NEF(X)$ is rational finite polyhedral. It follows that
the $\nem (X)$ is rational finite polyhedral. Thus,
Theorem 4.3.1 follows from

\proclaim{Theorem 4.3.1'}
Let $X$ be a $3$-dimensional Calabi-Yau
manifold with the finite rational polyhedral
cone $NEF(X)$ (equivalently,
$\nem (X)$). Let us assume that $X$ does
not have a small extremal ray and
a $nef$ rational element $h$ such that $h^3=0$
(equivalently, all faces of
$\nem (X)$ have Kodaira dimension $3$).
Then we have the following statements (1), (2) and (3) about $X$:

(1) $X$ does not have a pair of extremal rays of the type $\bb_2$ and
$\nem (X)$ is simplicial;

(2) $X$ has $\le 2$ extremal rays of the type (I) (more generally,
$k+(n_1-1)+...+(n_l-1)\le 2$ with notations of Theorem 3.1);

(3) (by the author and V.V. Shokurov) The Picard number
$\rho (X)=\dim N_1(X) \le 4k+5l_2+29 \le 40$.
\endproclaim

\demo{Proof}
 From the statement (1) of Theorem 1.3.2, the statement (1) follows.

Let $\E$ be an extremal set of Kodaira dimension $3$ of
extremal rays of the type (I) or (II) on $X$.
By Proposition 1.2.1, the $\E$ is elliptic. By Theorem 2.2.6 and
calculations of Section 4.2, the graph $G(\E )$ is one of elliptic graphs
of the Tables 4.1, 4.4, 4.5, 4.6 and Theorem 4.2.2.1.

In particular, $\#\E \le 2$ if $\E$ contains extremal rays of
the type (II) and the graph $G(\E )$ is full (i.e. any
two vertices are joined by non-single arrows).
Thus, for the constant $q(X)$ of
Basic Theorems 1.3.2 and 3.2 (see Definition
1.3.1), we have $q(X)\le 2$. By Theorems 1.3.2, (2) and 3.2, (2)
we get the statement (2).

Now let us estimate $\rho (X)$. To demonstrate how
Basic Theorems 1.3.2, (3) and 3.2, (3) do work, we first
give worse estimates for $\rho (X)$ which give these
general Theorems.

First, let us apply Theorem 1.3.2, (3).

Considering Tables 4.1, 4.4, 4.5 and 4.6,
one can easily see that for $d\ge 2$, we have
(where we take the distance in the graph $G(\E )$):
$$
\sharp \{ (R_1, R_2)\in \E \times \E
\mid 1 \le \rho (R_1,R_2)\le d\} \le C_1(d) \sharp \E ;
$$
and
$$
\sharp \{ (R_1, R_2)\in \E \times \E
\mid d+1\le \rho (R_1,R_2) \le 2d+1\} \le C_2(d) \sharp \E .
$$
where
$$
C_1(d)\le 2d,\ \ \  C_2(d)\le 2(d+1).
\tag4-3-1
$$

Let $\La$ be an $E$-set on $X$. Since all extremal rays on $X$ are
of the type (I) or (II)
and $X$ does not have a pair of extremal rays of the type
$\bb_2$, the $\La$ is Lanner by Lemma 1.2.11. Using calculations of
Section 4.2, we then get that the diameter
$$
\text{diam\ }(G(\La ))\le 4.
$$
Here the maximum $\text{diam\ }(G(\La ))=4$ we get for
the Lanner diagram with $6$ vertices
of the Table 4.4 and the Lanner diagrams
${\bold F}_4^\bullet (k;a,b)_1$ and ${\bold F}_4^\bullet (k;a,b)_2$
of the type (D) of the Table 4.5. Thus, the constant $d(X)$ of
Basic Theorem 1.3.2 has the estimate $d(X)\le 4$.
Here, for Lanner diagrams of the type (D) (Table 4.5),
we use the result of V.V. Shokurov (see Appendix):
for an extremal ray $R$ of the type (I) on $X$, there exists a curve
$C\in R$ such that $C\cdot D(R)=-k$ where $1\le k\le 3$.
By (4-3-1), we can apply Theorem 1.3.2 with the constants
$C_1(X)=8$ and $C_2(X)=10$. By Theorem 1.3.2, we get
$\rho (X)=\dim N_1(X) \le (16/3)C_1(X) + 4C_2(X) + 6=
(16/3)\cdot 8+4\cdot 10+6=88+2/3<89$. Thus, Basic Theorem 1.3.2
gives the estimate
$$
\rho (X)\le 88.
$$

Now let us apply Basic Theorem 3.2, (3).
Similar considerations give constants:
$$
n(X)_D\le 4, n(X)_C\le 5, n(X)_A\le 4, d(X)_A \le 2,\ \ C(X)_A\le 5.
$$
Thus, by Basic Theorem 3.2, (3),
$$
\rho (X)\le 4k+5l_2+8\cdot 5+6=4k+5l_2+46,
$$
and
$$
\rho (X)\le 56.
$$
(since $k+l\le q(X)=2$).

To get the strong estimate
$$
\rho (X)\le 4k+5l_2+30\le 40,
$$
we should change for this case the proof of the statement (3)
of Basic Theorem 3.2.
We should change the definition (3-3) of the weight $\sigma_A (\angle)$
of oriented angles of $\gamma$ according
to \'E.B. Vinberg \cite{V1, \S 6, Sect. 1}. We denote
this new weight as $\sigma_{AV}(\angle )$.

This is the following: We consider the extremal set
$\R(\angle )^\prime$ (see Definition 3.1)
and connected components of the graph $G(R(\angle )^\prime)$.
This graph contains non-single arrows only and this connected
components have types
$$
{\bold A}_n, {\bold B}_n, {\bold C}_n, {\bold D}_n, {\bold E}_6,
{\bold E}_7, {\bold E}_8, {\bold F}_4, {\bold G}_2
$$
by Theorem 3.2.1.1, (a).

The weight $\sigma_{AV}(\angle )$ is non-negative. It is
not equal to zero only if
both extremal rays $R_1(\angle )$, $R_2(\angle )$ belong to one
connected component $S$ of the graph
$G(\R(\angle )^\prime)$ and we have one of five cases below:

(I) $\rho_A(R_1(\angle ), R_2(\angle ))=1$ (i.e.
$R_1(\angle)$ and $R_2(\angle)$ are adjoined in $G(S)$;

(II) the connected component $S$ contains $\le 7$ vertices;

(III) the connected component $S$ is classical (i.e. it has the type
${\bold A}_n$, ${\bold B}_n$, ${\bold C}_n$ or ${\bold D}_n$)
and $R_1(\angle), R_2(\angle )$ belong to the terminal interval of
the order $\le 6$ of $G(S)$ (for $S$ of the type ${\bold D}_n$,
$n\ge 5$, by definition, the terminal interval is
a connected subgraph of $G(S)$ which
contains a terminal vertex of $S$ and is invariant relative to
the involution automorphism of ${\bold D}_n$.

(IV) $S$ has the type ${\bold E}_8$ and the pair
$R_1(\angle )$, $R_2(\angle )$ is different from
pairs of vertices of ${\bold E}_8$ below marked by $\otimes$:

\vskip10pt

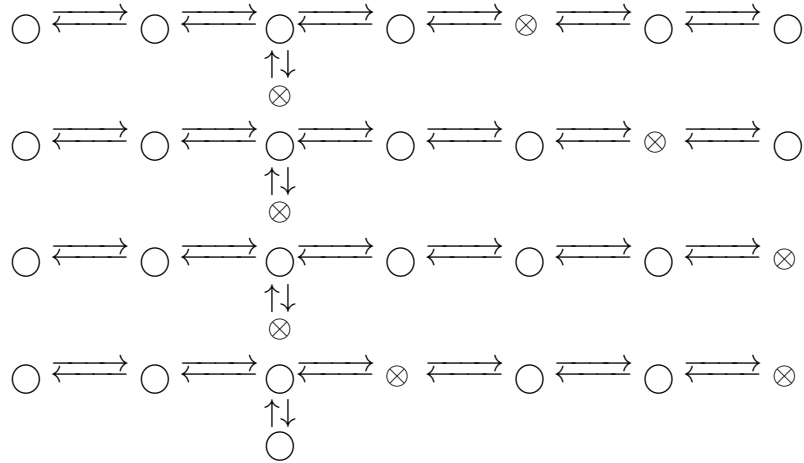
\centerline{\hfil
\hbox to 260pt{
\hfil
\hbox to15pt{\hfill\vbox to30pt{\vskip2pt $\bigcirc$ \vfill}\hfill}
\hskip8pt plus2pt minus2pt
\hbox to30pt{\hfill
\vbox to30pt{
\nointerlineskip
\hbox to25pt{\rightarrowfill} \nointerlineskip
\hbox to25pt{\leftarrowfill} \nointerlineskip
\vfill}
\hfill}
\hskip-14pt plus2pt minus2pt
\hbox to15pt{\hfill \vbox to30pt{\vskip2pt $\bigcirc$ \vfill}\hfill}
\hskip8pt plus2pt minus2pt
\hbox to30pt{\hfill
\vbox to30pt{\nointerlineskip
\hbox to25pt{\rightarrowfill} \nointerlineskip
\hbox to25pt{\leftarrowfill} \nointerlineskip
\vfill}
\hfill}
\hskip-5pt plus2pt minus2pt
\hbox to 20pt{\hfill
\vbox to 30pt {\vskip2pt
\hbox to15pt{\hfill $\bigcirc$ \hfill}
\nointerlineskip\vfill
\hbox to15pt{\hfill $\uparrow\downarrow$ \hfill}
\nointerlineskip\vfill
\hbox to15pt{\hfill $\otimes$ \hfill}
}
\hfill}
\hskip-10pt plus2pt minus2pt
\hbox to30pt{\hfill
\vbox to30pt{\nointerlineskip
\hbox to25pt{\rightarrowfill} \nointerlineskip
\hbox to25pt{\leftarrowfill} \nointerlineskip
\vfill}
\hfill}
\hskip-14pt plus2pt minus2pt
\hbox to15pt{\hfill \vbox to30pt{\vskip2pt $\bigcirc$ \vfill}\hfill}
\hskip8pt plus2pt minus2pt
\hbox to30pt{\hfill
\vbox to30pt{\nointerlineskip
\hbox to25pt{\rightarrowfill} \nointerlineskip
\hbox to25pt{\leftarrowfill} \nointerlineskip
\vfill}
\hfill}
\hskip-14pt plus2pt minus2pt
\hbox to15pt{\hfill \vbox to30pt{\vskip2pt $\otimes$ \vfill}\hfill}
\hskip8pt plus2pt minus2pt
\hbox to30pt{\hfill
\vbox to30pt{\nointerlineskip
\hbox to25pt{\rightarrowfill} \nointerlineskip
\hbox to25pt{\leftarrowfill} \nointerlineskip
\vfill}
\hfill}
\hskip-14pt plus2pt minus2pt
\hbox to15pt{\hfill \vbox to30pt{\vskip2pt $\bigcirc$ \vfill}\hfill}
\hskip8pt plus2pt minus2pt
\hbox to30pt{\hfill
\vbox to30pt{\nointerlineskip
\hbox to25pt{\rightarrowfill} \nointerlineskip
\hbox to25pt{\leftarrowfill} \nointerlineskip
\vfill}
\hfill}
\hskip-14pt plus2pt minus2pt
\hbox to15pt{\hfill \vbox to30pt{\vskip2pt $\bigcirc$ \vfill}\hfill}
\hfil
}
\hskip10pt}


\vskip5pt

\centerline{\hfil
\hbox to 260pt{
\hfil
\hbox to15pt{\hfill\vbox to30pt{\vskip2pt $\bigcirc$ \vfill}\hfill}
\hskip8pt plus2pt minus2pt
\hbox to30pt{\hfill
\vbox to30pt{
\nointerlineskip
\hbox to25pt{\rightarrowfill} \nointerlineskip
\hbox to25pt{\leftarrowfill} \nointerlineskip
\vfill}
\hfill}
\hskip-14pt plus2pt minus2pt
\hbox to15pt{\hfill \vbox to30pt{\vskip2pt $\bigcirc$ \vfill}\hfill}
\hskip8pt plus2pt minus2pt
\hbox to30pt{\hfill
\vbox to30pt{\nointerlineskip
\hbox to25pt{\rightarrowfill} \nointerlineskip
\hbox to25pt{\leftarrowfill} \nointerlineskip
\vfill}
\hfill}
\hskip-5pt plus2pt minus2pt
\hbox to 20pt{\hfill
\vbox to 30pt {\vskip2pt
\hbox to15pt{\hfill $\bigcirc$ \hfill}
\nointerlineskip\vfill
\hbox to15pt{\hfill $\uparrow\downarrow$ \hfill}
\nointerlineskip\vfill
\hbox to15pt{\hfill $\otimes$ \hfill}
}
\hfill}
\hskip-10pt plus2pt minus2pt
\hbox to30pt{\hfill
\vbox to30pt{\nointerlineskip
\hbox to25pt{\rightarrowfill} \nointerlineskip
\hbox to25pt{\leftarrowfill} \nointerlineskip
\vfill}
\hfill}
\hskip-14pt plus2pt minus2pt
\hbox to15pt{\hfill \vbox to30pt{\vskip2pt $\bigcirc$ \vfill}\hfill}
\hskip8pt plus2pt minus2pt
\hbox to30pt{\hfill
\vbox to30pt{\nointerlineskip
\hbox to25pt{\rightarrowfill} \nointerlineskip
\hbox to25pt{\leftarrowfill} \nointerlineskip
\vfill}
\hfill}
\hskip-14pt plus2pt minus2pt
\hbox to15pt{\hfill \vbox to30pt{\vskip2pt $\bigcirc$ \vfill}\hfill}
\hskip8pt plus2pt minus2pt
\hbox to30pt{\hfill
\vbox to30pt{\nointerlineskip
\hbox to25pt{\rightarrowfill} \nointerlineskip
\hbox to25pt{\leftarrowfill} \nointerlineskip
\vfill}
\hfill}
\hskip-14pt plus2pt minus2pt
\hbox to15pt{\hfill \vbox to30pt{\vskip2pt $\otimes$ \vfill}\hfill}
\hskip8pt plus2pt minus2pt
\hbox to30pt{\hfill
\vbox to30pt{\nointerlineskip
\hbox to25pt{\rightarrowfill} \nointerlineskip
\hbox to25pt{\leftarrowfill} \nointerlineskip
\vfill}
\hfill}
\hskip-14pt plus2pt minus2pt
\hbox to15pt{\hfill \vbox to30pt{\vskip2pt $\bigcirc$ \vfill}\hfill}
\hfil
}
\hskip10pt}


\vskip5pt

\centerline{\hfil
\hbox to 260pt{
\hfil
\hbox to15pt{\hfill\vbox to30pt{\vskip2pt $\bigcirc$ \vfill}\hfill}
\hskip8pt plus2pt minus2pt
\hbox to30pt{\hfill
\vbox to30pt{
\nointerlineskip
\hbox to25pt{\rightarrowfill} \nointerlineskip
\hbox to25pt{\leftarrowfill} \nointerlineskip
\vfill}
\hfill}
\hskip-14pt plus2pt minus2pt
\hbox to15pt{\hfill \vbox to30pt{\vskip2pt $\bigcirc$ \vfill}\hfill}
\hskip8pt plus2pt minus2pt
\hbox to30pt{\hfill
\vbox to30pt{\nointerlineskip
\hbox to25pt{\rightarrowfill} \nointerlineskip
\hbox to25pt{\leftarrowfill} \nointerlineskip
\vfill}
\hfill}
\hskip-5pt plus2pt minus2pt
\hbox to 20pt{\hfill
\vbox to 30pt {\vskip2pt
\hbox to15pt{\hfill $\bigcirc$ \hfill}
\nointerlineskip\vfill
\hbox to15pt{\hfill $\uparrow\downarrow$ \hfill}
\nointerlineskip\vfill
\hbox to15pt{\hfill $\otimes$ \hfill}
}
\hfill}
\hskip-10pt plus2pt minus2pt
\hbox to30pt{\hfill
\vbox to30pt{\nointerlineskip
\hbox to25pt{\rightarrowfill} \nointerlineskip
\hbox to25pt{\leftarrowfill} \nointerlineskip
\vfill}
\hfill}
\hskip-14pt plus2pt minus2pt
\hbox to15pt{\hfill \vbox to30pt{\vskip2pt $\bigcirc$ \vfill}\hfill}
\hskip8pt plus2pt minus2pt
\hbox to30pt{\hfill
\vbox to30pt{\nointerlineskip
\hbox to25pt{\rightarrowfill} \nointerlineskip
\hbox to25pt{\leftarrowfill} \nointerlineskip
\vfill}
\hfill}
\hskip-14pt plus2pt minus2pt
\hbox to15pt{\hfill \vbox to30pt{\vskip2pt $\bigcirc$ \vfill}\hfill}
\hskip8pt plus2pt minus2pt
\hbox to30pt{\hfill
\vbox to30pt{\nointerlineskip
\hbox to25pt{\rightarrowfill} \nointerlineskip
\hbox to25pt{\leftarrowfill} \nointerlineskip
\vfill}
\hfill}
\hskip-14pt plus2pt minus2pt
\hbox to15pt{\hfill \vbox to30pt{\vskip2pt $\bigcirc$ \vfill}\hfill}
\hskip8pt plus2pt minus2pt
\hbox to30pt{\hfill
\vbox to30pt{\nointerlineskip
\hbox to25pt{\rightarrowfill} \nointerlineskip
\hbox to25pt{\leftarrowfill} \nointerlineskip
\vfill}
\hfill}
\hskip-14pt plus2pt minus2pt
\hbox to15pt{\hfill \vbox to30pt{\vskip2pt $\otimes$ \vfill}\hfill}
\hfil
}
\hskip10pt}


\vskip5pt

\centerline{\hfil
\hbox to 260pt{
\hfil
\hbox to15pt{\hfill\vbox to30pt{\vskip2pt $\bigcirc$ \vfill}\hfill}
\hskip8pt plus2pt minus2pt
\hbox to30pt{\hfill
\vbox to30pt{
\nointerlineskip
\hbox to25pt{\rightarrowfill} \nointerlineskip
\hbox to25pt{\leftarrowfill} \nointerlineskip
\vfill}
\hfill}
\hskip-14pt plus2pt minus2pt
\hbox to15pt{\hfill \vbox to30pt{\vskip2pt $\bigcirc$ \vfill}\hfill}
\hskip8pt plus2pt minus2pt
\hbox to30pt{\hfill
\vbox to30pt{\nointerlineskip
\hbox to25pt{\rightarrowfill} \nointerlineskip
\hbox to25pt{\leftarrowfill} \nointerlineskip
\vfill}
\hfill}
\hskip-5pt plus2pt minus2pt
\hbox to 20pt{\hfill
\vbox to 30pt {\vskip2pt
\hbox to15pt{\hfill $\bigcirc$ \hfill}
\nointerlineskip\vfill
\hbox to15pt{\hfill $\uparrow\downarrow$ \hfill}
\nointerlineskip\vfill
\hbox to15pt{\hfill $\bigcirc$ \hfill}
}
\hfill}
\hskip-10pt plus2pt minus2pt
\hbox to30pt{\hfill
\vbox to30pt{\nointerlineskip
\hbox to25pt{\rightarrowfill} \nointerlineskip
\hbox to25pt{\leftarrowfill} \nointerlineskip
\vfill}
\hfill}
\hskip-14pt plus2pt minus2pt
\hbox to15pt{\hfill \vbox to30pt{\vskip2pt $\otimes$ \vfill}\hfill}
\hskip8pt plus2pt minus2pt
\hbox to30pt{\hfill
\vbox to30pt{\nointerlineskip
\hbox to25pt{\rightarrowfill} \nointerlineskip
\hbox to25pt{\leftarrowfill} \nointerlineskip
\vfill}
\hfill}
\hskip-14pt plus2pt minus2pt
\hbox to15pt{\hfill \vbox to30pt{\vskip2pt $\bigcirc$ \vfill}\hfill}
\hskip8pt plus2pt minus2pt
\hbox to30pt{\hfill
\vbox to30pt{\nointerlineskip
\hbox to25pt{\rightarrowfill} \nointerlineskip
\hbox to25pt{\leftarrowfill} \nointerlineskip
\vfill}
\hfill}
\hskip-14pt plus2pt minus2pt
\hbox to15pt{\hfill \vbox to30pt{\vskip2pt $\bigcirc$ \vfill}\hfill}
\hskip8pt plus2pt minus2pt
\hbox to30pt{\hfill
\vbox to30pt{\nointerlineskip
\hbox to25pt{\rightarrowfill} \nointerlineskip
\hbox to25pt{\leftarrowfill} \nointerlineskip
\vfill}
\hfill}
\hskip-14pt plus2pt minus2pt
\hbox to15pt{\hfill \vbox to30pt{\vskip2pt $\otimes$ \vfill}\hfill}
\hfil
}
\hskip10pt}


\vskip10pt

The weight $\sigma_{AV} (\angle )=1$ if

(V) $S$ has $\le 4$ vertices.

In all other cases (I)--(IV) above, $\sigma_{AV} (\angle )=1/2$.

\smallpagebreak

Now we should check conditions of Vinberg's Lemma 3.7 with the
constants $C=3$ and $D=0$ which gives the desirable estimate:
$\dim \gamma <30$.

The proof of the condition (1) is very similar to Vinberg
\cite{V1, \S 6, Sect. 2} and is not difficult.

To prove the condition (2), one should also follow to Vinberg
\cite{V1, \S 6, Sects. 3---12}. For our case,
elliptic diagrams are only "crystallographic", i. e. of
the types
$$
{\bold A}_n, {\bold B}_n, {\bold C}_n, {\bold D}_n, {\bold E}_6,
{\bold E}_7, {\bold E}_8, {\bold F}_4, {\bold G}_2,
$$
and $E$-sets (or Lanner diagrams) have only types of Table 4.3.
One should follow to Vinberg step by step (we should say that
his considerations are long and very delicate, and it is
a hard work), and check that in fact for our "crystallographic
case" he uses only two arguments which work for our situation:
There do not exist a configuration of extremal rays of the type (II)
with one of extended Dynkin diagram of Table  4.2, because then
some linear combination of divisors of these rays gives a $nef$ element
with zero cube. Besides, two non-extremal sets of extremal rays of
the type (I) or (II) cannot be orthogonal (Lemma 3.4).

Only for the "non-crystallographic case", Vinberg uses
{\it superhyperbolicity} arguments which probably
do not work for our situation.

This finishes the proof of Theorem.
\enddemo
\enddemo

\subhead
5. ${\Bbb Q}$-factorial models of Calabi-Yau 3-folds
\endsubhead

Applying Diagram Method to an arbitrary $3$-dimensional
Calabi-Yau manifold $X$, we have to avoid two problems:
Mori cone $\nem (X)$ may not be finite polyhedral,
and $X$ may have small extremal rays.
Here we want to discuss one possibility
to avoid these problems and involve non-polyhedral case and
small extremal rays to the game.

One can make several transformations (i) and (ii) below:

(i) Contraction of a divisorial extremal ray.

(ii) Flop in a small extremal ray.

Repeating this operations, we get some $3$-fold $Y$.
The $3$-fold $Y$ still has the most important for
us properties: $Y$ has $\Bbb Q$-factorial (canonical)
singularities and belongs to the class $\LT$.
See Kawamata \cite{Ka2} and Shokurov \cite{Sh2}.

\definition{Definition 5.1}
A 3-fold $Y$ one can get
starting from a 3-dimensional Calabi-Yau manifold $X$ and
repeating transformations (i) and (ii) is
called a {\it $\Bbb Q$-factorial model of the Calabi-Yau manifold $X$}.
\enddefinition

Considering $\Bbb Q$-factorial models, we get a chance
to avoid cases when either Mori cone is not finite polyhedral or
there exists a small extremal ray.
If the $\Bbb Q$-factorial model $Y$ has a finite polyhedral Mori
cone and does not have a small extremal ray, we can apply Diagram
Method to $Y$ to find a rational $nef$ element $h$ on $Y$ with
cube $0$. By the way, using the statement (1) of
Basic Theorem 1.3.2, we can prove that transformations (i) and (ii)
do not make situation worse from Diagram Method point of view.

\proclaim{Lemma 5.2} Let $Y$ be a
$\Bbb Q$-factorial model such that we have properties (a), (b) and
(c) below:

(a) $\nem (Y)$ is finite polyhedral;

(b) $Y$ does not have a small
extremal ray;

(c) $Y$ does not have a rational $nef$ element $h$ with
$h^3=0$.

Then, starting from $Y$, repeating of transformations (i) and (ii)
preserves properties (a), (b) and (c).
\endproclaim

\demo{Proof} Let us consider a contraction $f:Y\to Y^\prime$
of a divisorial extremal ray $R$ on $Y$.
One easily can see that
$Y^\prime$ has properties (a) and (c) if $Y$ has these
properties. Besides, if $Y$ does not have a small extremal ray, then
$Y^\prime$ has a
small extremal ray only if the divisor $D(R)$ contains another divisorial
extremal ray $Q$. The image of $Q$ then gives a small extremal ray on
$Y^\prime$.
Thus, $R$ and $Q$ define a pair of the type $\bb_2$ on $Y$.
This is impossible by Theorem 1.3.2, (1).

If $Y$ does not have a small extremal ray,
then $Y$ does not have a flop and one
does not have the transformation (ii).
\enddemo

We suggest the following

\proclaim{Conjecture 5.3} There are absolute constants
$q$, $d$, $C_1$, $C_2$, $n_D$,
$n_C$, $n_A$, $d_A$, $C_A$ such that for any
$3$-dimensional Calabi-Yau manifold $X$ and
any its
${\Bbb Q}$-factorial model $Y$ we have estimates
$q(Y)\le q$, $d(Y)\le d$,
$C_1(Y)\le C_1$ and $C_2(Y)\le C_2$ for invariants of Definition 1.3.1,
and $n(Y)_D\le n_D$,
$n(Y)_C\le n_C$, $n(Y)_A\le n_A$, $d(X)_A\le d_A$ and
$C(X)_A\le C_A$ for invariants of Definition 3.1.
\endproclaim

If this conjecture does hold, applying Theorems 1.3.2 and 3.2, we get
an absolute estimate for $\rho (Y)$ if $Y$ satisfies conditions
(a), (b) and (c). Equivalently, if $X$ has a ${\Bbb Q}$-factorial
model $Y$ with a finite polyhedral Mori cone, without small extremal
rays and with big $\rho (Y)$, then $Y$ has a rational
$nef$ element with cube zero. One can think that existence of rational
$nef$ element with cube $0$ for a ${\Bbb Q}$-factorial model $Y$ of $X$
is of similar importance as for $X$.

Unfortunately, now we can prove Conjecture 5.3 only for very special
$\Bbb Q$-factorial models $Y$.

We can consider part of results of Sect. 4,
as the proof of Conjecture 5.3 for non-singular
models (i.e. for $3$-dimensional Calabi-Yau manifolds).

It is possible to extend these results for more large class of models
which are still very special.

\definition{Definition 5.4} A $\Bbb Q$-factorial model $Y$ of
a $3$-dimensional Calabi-Yau manifold is called {\it very good} if
there exists a $3$-dimensional
Calabi-Yau manifold $X$ and an extremal set
${\Cal R}=\{R_1,...,R_k \}$ of divisorial extremal rays on $X$ with
different divisors $D(R_1),...,D(R_k)$ such that we have the
following properties (i) and (ii):

(i) $Y$ is a contraction $f: X \to Y$ of the face
$\gamma =R_1+...+R_k$ of $\nem (X)$;

(ii) For any divisorial extremal ray $R$ of $Y$ there exists a
divisorial extremal ray $\tilde{R}$ of $X$ such that
$R=f(\tilde{R})$. This means that for the contraction
$f_R:Y\to Y^\prime$ of $R$, the composition
$f_R \cdot f:X\to Y^\prime$ is the contraction of the face
$R_1+...+R_k+\tilde{R}$.
\enddefinition

We have

\proclaim{Theorem 5.5} Let $Y$ be a very good $\Bbb Q$-factorial
model of a $3$-dimensional Calabi-Yau manifold. Then $Y$ has
constants
$$
q(Y)\le 3, d(Y)\le 8, C_1(Y) \le 16, C_2(Y) \le 18,
$$
and
$$
n(Y)_D\le 8, n(Y)_C\le 9, n(Y)_A\le 9, d(Y)_A\le 8, C(Y)_A\le 17.
$$
In particular, by Theorem 1.2.3,
$Y$ has $\le 3$ extremal rays of the type (I) and
$$
\rho (Y) \le (16/3)C_1(Y)+4C_2(Y)+6<164
$$
if $Y$ satisfies conditions (a), (b) and (c).
\endproclaim

Considering preimage of the $nef$ element on $Y$
with cube $0$, we get

\proclaim{Corollary 5.6} Let $X$ be a Calabi-Yau manifold $X$
which has a very good $\Bbb Q$-factorial model $Y$ such that
$Y$ has a finite polyhedral Mori cone $\nem (Y)$ and $Y$
does not have a small extremal ray.

Then $X$ has a rational $nef$ element $h$ with
$h^3=0$ if $\rho (Y)\ge 164$.
\endproclaim

\demo{Sketch of the proof of Theorem 5.5} We use notation of Definition 5.4.
For a set $Q$ of divisorial extremal rays on $Y$,
we denote $f^{-1}(Q) =
\{R_1,...,R_k \}\cup \tilde{Q}$ where
$\tilde{Q}=\{\tilde{R} \ \mid \ R\in Q\}$. The set $f^{-1}(Q)$ is called
the preimage of $Q$ and $\tilde{Q}$ is called the proper preimage of $Q$.

Next general statements will be very useful.

\proclaim{Lemma 5.7} Let $\E$ be an elliptic set of divisorial
extremal rays on $Y$. Then $f^{-1}(\E )$ is elliptic.
\endproclaim

\demo{Proof} Let $\E=\{Q_1,...,Q_t\}$. Since $\E$ is elliptic,
there are positive $a_1,...,a_t$ such that
$Q_j \cdot (a_1D(Q_1)+ \cdots +a_tD(Q_t))<0$ for all
$1\le j \le t$. By Lemma 1.2.2,
the set $\{R_1,...,R_k\}$ is elliptic. Hence, there are positive
$b_1,...,b_k$ such that
$R_i \cdot (b_1D(R_1)+ \cdots + b_kD(R_k))<0$ for all $1\le i \le k$.
Evidently,
$$
f^\ast (a_1D(Q_1)+ \cdots +a_tD(Q_t))=
c_1D(R_1)+ \cdots + c_k D(R_k)+
a_1D(\tilde{Q}_1)+ \cdots + a_tD(\tilde{Q}_t)
$$
with non-negative $c_1,...,c_k$. By projection formula,
$$
R_i \cdot
(c_1D(R_1)+ \cdots + c_k D(R_k)+
a_1D(\tilde{Q}_1)+ \cdots + a_tD(\tilde{Q}_t))=0
$$
for all $1\le i \le k$, and
$$
\split
&\tilde{Q}_j \cdot
(c_1D(R_1)+ \cdots + c_k D(R_k)+
a_1D(\tilde{Q}_1)+ \cdots + a_tD(\tilde{Q}_t))=\\
=&Q_j \cdot
(a_1D(Q_1)+ \cdots +a_tD(Q_t))) < 0.
\endsplit
$$
It follows that for a small $\epsilon >0$
the divisor
$$
\split
D(f^{-1}(Q)) &= \epsilon (b_1D(R_1)+ \cdots + b_kD(R_k))+\\
& +c_1D(R_1)+ \cdots + c_k D(R_k)+
a_1D(\tilde{Q}_1)+ \cdots + a_tD(\tilde{Q}_t)
\endsplit
$$
has positive coefficients and has
negative intersection with each element of $f^{-1}(Q)$.
Thus, $f^{-1}(Q)$ is elliptic.
\enddemo

\proclaim{Lemma 5.8} Let $\La$ be a Lanner set of divisorial
extremal rays on $Y$. Then $f^{-1}(\La )$ contains a quasi-Lanner
subset $\La^\prime$ such that
$\tilde{\La}\subset \La^\prime$ and such that for any non-empty
subset $U\subset \tilde{\La}$ the set
$\{R_1,...,R_k\}\cup (\La^\prime - U)$ is elliptic.
\endproclaim

\demo{Proof} Let $\La=\{Q_1,...,Q_t\}$.
Since $\La$ is Lanner, there are positive
$a_1,...,a_t$ such that
$R_j\cdot (a_1D(Q_1)+ \cdots +a_tD(Q_t))\ge 0$ for all $1\le j\le t$ and
one of these inequalities is strict.
Evidently, $f^\ast (a_1D(Q_1)+ \cdots +a_tD(Q_t))=
b_1D(R_1)+ \cdots +b_kD(R_k)+
a_1D(\tilde{Q_1})+ \cdots +a_tD(\tilde{Q_t})$
where $b_j\ge 0$. By projection formula,
$$
R\cdot (b_1D(R_1)+\cdots + b_kD(R_k)+
a_1D(\tilde{Q_1})+\cdots +a_tD(\tilde{Q_t})) \ge 0
$$
for any $R\in f^{-1}(\La)$, and one of these inequalities is strict.
Thus, the set $f^{-1}(\La)$ is not semi-elliptic. By Proposition
2.2.8, $f^{-1}(\La)$ contains a quasi-Lanner subset $\La^\prime$.
Since each proper subset of $\La$ is elliptic, by Lemma 5.7
any subset of $f^{-1}(\La)$ is elliptic if it does not contain
$\tilde{\La}$. It follows that
$\tilde{\La}\subset \La^\prime$ and
$\{R_1,...,R_k\}\cup (\La^\prime - U)$ is elliptic if
$U$ is a non-empty subset of $\tilde{\La}$.
\enddemo

Let us continue the proof of Theorem.

Let $\La$ be an $E$-set of divisorial extremal rays on $Y$
such that any proper subset of $\La$ is extremal
of Kodaira dimension $3$ and $\La$
satisfies the condition (iii) of Sect. 1.1.
By results of Sect. 1.2, the set $\La$ is Lanner. By Lemma 5.8,
$\tilde{\La}\subset \La^\prime$ where $\La^\prime$ is a
quasi-Lanner set of divisorial extremal rays on $X$. Considering
image of $\La^\prime$ by the morphism $f$, one sees
that
$$
\text{diam\ } \La \le \text{diam\ }\La^\prime.
$$
Any proper subset of $\La^\prime$ is either elliptic or
connected parabolic (see Proposition 2.2.8). Using
Theorem 2.2.6 and description in Sect. 4 of elliptic and
connected parabolic sets of divisorial extremal rays on
$3$-dimensional Calabi-Yau manifolds, one can describe possible
graphs of quasi-Lanner sets of divisorial extremal rays on $X$.
In particular, one can see that
$$
\sharp \La^\prime \le 10,\ \ \ \text{diam\ }\La^\prime \le 8
$$
for any quasi-Lanner set $\La ^\prime$ of divisorial extremal rays
on a Calabi-Yau $3$-dimensional manifold $X$
(compare with \cite{P, Fig. 1}).
It follows that $n(Y)_A\le 9$ and
$$
d(Y) \le 8.
$$
By Theorem 2.2.6, we also get
$n(Y)_D\le 8$ and $n(Y)_C\le 9$.

Now let us consider an extremal set $\E$ of Kodaira dimension $3$
of divisorial extremal rays on $Y$.
The set $f^{-1}(\E )$ is also extremal of Kodaira dimension $3$
and contains ${\Cal R}=\{R_1,...,R_k\}$. For
$Q_1,Q_2\in \E$ we have
$$
\rho (Q_1,Q_2)=\rho _{\Cal R}(\tilde{Q}_1,\tilde{Q}_2)
$$
where $\rho (Q_1,Q_2)$ is the distance in the oriented graph
$G(\E)$ and
$$
\rho_{\Cal R}(\tilde{Q}_1,\tilde{Q}_2)=
\min_{\gamma}{(\rho (\gamma)-v(\gamma\cap {\Cal R}))}
$$
where $\gamma$ is an oriented path in $G(f^{-1}(\E ))$
joining $\tilde{Q}_1$ and $\tilde{Q}_2$ and
$v(\gamma \cap {\Cal R})$ is the number of vertices of
$\gamma$ which belong to ${\Cal R}$.

Using description in Sect. 4 of graphs of elliptic
(in particular, extremal of Kodaira
dimension $3$) sets of divisorial extremal rays on $X$,
one obtains that
$q(Y)\le 3$, $C_1(Y) \le 16$, $C_2(Y) \le 18$, $C(Y) \le 17$.
This finishes the proof.
\enddemo

We hope to describe more precisely quasi-Lanner sets
of divisorial extremal rays on $3$-dimensional
Calabi-Yau manifolds
and give better estimates for Theorem 5.5 and Corollary 5.6
in more advanced variant of this preprint.

We give another interesting application of Theorem 5.5.

\proclaim{Corollary 5.9} Let $X$ be a $3$-dimensional
Calabi-Yau manifold and Mori cone $\nem (X)$ is generated
by a finite set of divisorial extremal rays. Let
$\gamma$ be a face of $nef$ cone $NEF(X)$ and
$R(\gamma^\perp)$ the set of all extremal rays orthogonal to
$\gamma$. Assume that $R(\gamma^\perp)$ does not contain an
extremal ray which belongs to a pair of the type $\bb_2$.
Then $\gamma$ contains a rational
$nef$ element $h$ with $h^3=0$ if $\dim \gamma > 163$.
\endproclaim

\demo{Proof} Let $\R=R(\gamma^\perp )$ and $f:X\to Y$ the
contraction of the face $\gamma ^\perp \subset \nem (X)$
generated by $\R$. Since $\R$ does not have pairs of
extremal rays of the type $\bb_2$,
the morphism $f$ is a sequence of
contractions of divisorial extremal rays which are images of
extremal rays from $\R$. Then $Y$ has $\Bbb Q$-factorial
singularities (canonical) and $\gamma= f^\ast (NEF(Y))$.
Since $\nem (X)$ is generated by a finite set of divisorial
extremal rays and $\R$ does not have an extremal ray
which belongs to a pair of the
type $\bb_2$, one sees that $Y$ does not have a small extremal ray
and is a very good $\Bbb Q$-factorial
model of $X$ with Mori cone $\nem (Y)$ generated by a finite set of
divisorial extremal rays.
If we additionally assume that $NEF(Y)$ (equivalently $\gamma$)
does not have a $nef$ element with cube $0$,
$\dim \gamma = \dim NEF(Y)<164$ by Theorem 5.5.
\enddemo

At last, we mention the following conjecture by D. Morrison which
is connected with the condition (a) above on Mori cone.

\proclaim{Conjecture 5.10} (by D. Morrison, \cite{Mor2}). For
a Calabi-Yau manifold $X$, the $nef$ cone $NEF(X)$ is rational
finite polyhedral modulo the group $\text{Aut}~X$
of biregular automorphisms of $X$. In particular,
$NEF(X)$ and Mori cone $\nem (X)$ are rational finite polyhedral
if $\text{Aut}~X$ is finite.
\endproclaim

For example, the automorphism group $\text{Aut}~X$ of a
$3$-dimensional Calabi-Yau manifold $X$ is finite if the cubic
intersection hypersurface ${\Cal W}_X$ is non-singular.

\subhead
6. Concluding remarks
\endsubhead

We want to give several remarks about importance of existence
of rational $nef$ elements with cube $0$ for Calabi-Yau 3-folds.

I.I. Piatetski-Shapiro and I.R. Shafarevich \cite{P\u S-\u S}
proved that a K3 surface $X$ has a rational $nef$ element
$h\in \text{Pic}~X$ with $h^2=0$ if and only if the Picard lattice
$\text{Pic}~X$ represents $0$, i. e. there exists
$0\not= x \in \text{Pic}~X$ such that $x^2=0$. In particular, this
is true if $\rho (X)=\dim \text{Pic}~X \ge 5$. Moreover, they
proved that the linear system $\vert h \vert$ defines
the elliptic fibration $\vert h \vert : X\to {\Bbb P}^1$ (i.e.
the general fiber is an elliptic curve) if $h\in \text{Pic}~X$ is
$nef$ and $h^2=0$.

One can ask about similar fact for Calabi-Yau $3$-folds:

\proclaim{Question 6.1}
Does exist a rational $nef$ element $h$ with $h^3=0$ for a
$3$-dimensional Calabi-Yau manifold $X$ if
$\rho (X)=\dim N_1(X)$ is sufficiently big?
\endproclaim

Affirmative answer to this question is important because
of two results below.

\proclaim{Theorem 6.2} (P.M.H. Wilson, \cite{W1,(3.2)'}).
Let $X$ be a $3$-dimensional Calabi-Yau manifold and
$h$ a rational $nef$ element with $h^3=0$. Assume that
$h^2\not\equiv 0$ and $h\cdot c_2(X)>0$. Then for a big $N$,
the linear system  $\vert Nh \vert$ defines an elliptic
fibration $\vert Nh \vert : X\to S$ (i.e. $S$ is a surface and
the general fiber is an elliptic curve).
\endproclaim

Here the condition $h^2\not\equiv 0$
is equivalent to that ${\Bbb C}h$ is not a singular point of the cubic
intersection hypersurface ${\Bbb P}{\Cal W}_X$.
By Y. Miyaoka \cite{Mi}, we have the inequality
$NEF(X)\cdot c_2(X)\ge 0$. Thus, for the "general case"
when the cubic intersection hypersurface
${\Bbb P}{\Cal W}_X$ is non-singular and
$NEF(X)\cdot c_2(X) > 0$, existence of a rational $nef$
element with cube $0$ gives existence of an elliptic fibration. Thus,
for this "general case", Theorems 4.3.1, 4.3.1' and Theorem 5.5 give
existence of an elliptic fibration on Calabi-Yau manifolds under
appropriate conditions.
If either cubic intersection hypersurface ${\Bbb P}{\Cal W}_X$
is singular or one only has the non-strict
condition $NEF(X)\cdot c_2(X)\ge 0$, Corollary 5.9 is
useful to satisfy conditions of Theorem 6.2, since Corollary 5.9
gives existence of rational $nef$ elements $h$
with cube $0$ in "many" faces of the $nef$ cone $NEF(X)$.
We mention that it is conjectured that the condition
$h\cdot c_2(X)>0$ is not essential for the existence of the
elliptic fibration $\vert Nh \vert$, see \cite{W3}. By K. Oguiso
\cite{O}, the linear system $\vert Nh \vert$ defines either
elliptic or K3 or Abelian surface fibration if
$\vert Nh \vert$ is not empty and $N$ is big.

Existence of an elliptic fibration for a Calabi-Yau $3$-fold
(this Calabi-Yau $3$-fold is called {\it elliptic} ) is
important because of the following result:

\proclaim{Theorem 6.3} (M. Gross, \cite{Gro}). There exists  a
finite set of families
${\Cal X}_1 \to {\Cal M}_1,...,{\Cal X}_n \to {\Cal M}_n$ of
elliptic Calabi-Yau $3$-folds with ${\Bbb Q}$-factorial terminal
singularities such that each elliptic Calabi-Yau $3$-fold $X$ with
${\Bbb Q}$-factorial terminal singularities is
birationally isomorphic to a fiber of one of these families.
(Here all families and isomorphisms preserve elliptic fiber
structure.)
\endproclaim

See also connected results by I. Dolgachev and
M. Gross \cite{D-Gro} and A. Grassi \cite{Gra}.

\smallpagebreak

Because of these Theorems 6.2. and 6.3 and results of this
paper, one can think that may be classification of Calabi-Yau
$3$-folds is simpler for the high Picard number.

\leftheadtext{Vyacheslav V. Shokurov}

\newpage

\rightheadtext{Anticanonical boundedness for curves}



\documentstyle{amsppt}
\magnification 1200
\define\ep{\varepsilon}  
\define\sn{\smallskip \noindent}

\define\bsq{$\blacksquare$}
\input cyracc.def

\document
\null
\vskip2cm

\centerline{\bf Appendix by V.V. Shokurov:
Anticanonical boundedness for curves}

\vskip1cm
The purpose of this note is to generalize slightly results of [K].

\medskip
\proclaim {Theorem} Let $f\colon X\to S$ be a projective
morphism of normal algebraic varieties over a field of
characteristic zero, and $D$ be an effective $\Bbb R$-divisor
on $X$ such that $K_X+D$ is $\Bbb R$-Cartier and
divisorially log terminal near the generic points of a subvariety $E$,
consisting of components of the {\rm degenerate locus\/}
$$\Cal Exc(f)\colon = \{x\in X\mid g \text{\rm\ is not finite at }x\}.$$
Then $E$ is covered by a family (possibly disconnected) of
effective 1-cycles
$\{C_{\lambda}\}/S$ with $(-K_X-D.C_{\lambda})\le 2n$ where
$n=\text{\rm\ dim } E/S$ (and even $<2n$ if
$K_X+D$ is Kawamata log terminal in the generic points of $E$ and
$X\not=E$).
Moreover, we could assume that the generic 1-cycles $C_{\lambda}$
are {\rm curves\/}, i.e., reduced and irreducible, when $K_X+D$ is
numerically definite, and the curves $C_{\lambda}$
with $(K_X+D.C_{\lambda})<0$ (resp. $\le 0$) are rational.
\endproclaim

The Theorem implies the following results.

\medskip
\proclaim {Corollary 1} If, in addition, $F$ is
an $\Bbb R$-Cartier divisor such that $K_X+D+F$
is also divisorially log terminal near the generic points of $E$,
$D+F$ is effective,
and $K_X+D$ is numerically semi-negative with respect
to $f$, then $E$ is covered by a family (possibly disconnected)
ofCKWARD
E effective 1-cycles $\{C_{\lambda}\}/S$ (resp. with the curves
as generic members when $K_X+D+F$ is numerically definite) with
$$(F.C_{\lambda})\ge -2n$$ (resp. $>-2n$ if $K_X+D+F$ is
Kawamata log terminal in the generic points of $E$, and $X\not=E$).
\endproclaim

For example, let $f$ be an extremal divisorial contraction of
a normal variety $X$ with the exceptional $\Bbb Q$-Cartier
divisor $F=E=\Cal Exc(f)$, and $K_X$ numerically semi-negative for $f$.
Then we cover $F$ by a family of rational
curves $\{C_{\lambda}\}/S$ with
$$(F.C_{\lambda})\ge -2\text{ dim }X+2.$$
In the low dimensional case, i.e., when $\text{ dim }X\le 3$
we have a more sharp inequality $\ge -\text{ dim }X.$
This is the best bound but it is known only in that case.
(See Remark 3 below.)

\sn
{\bf Proof.} According to the Theorem, take a family of 1-cycles
$\{C_{\lambda}\}/S$ with respect to the second log divisor
$K_X+D+F$. Then
$$(-K_X-D.C_{\lambda})-(F.C_{\lambda})=(-K_X-D-F.C_{\lambda})\le 2n$$
(resp. $<2n$ if $K_X+D+F$ is
Kawamata log terminal in the generic points of $E$, and $X\not=E$).
But we assume that $(-K_X-D.C_{\lambda})\ge 0$ which
gives the required inequality.\bsq

\medskip
\proclaim {Corollary 2 ([Ko1], cf. also [K])} Let $f\colon X\to S$
be a projective morphism of normal algebraic varieties over a field of
characteristic zero, and $D$ be an effective $\Bbb R$-divisor
on $X$ such that $K_X+D$ is $\Bbb R$-Cartier,
let $H$ be an $f$-ample $\Bbb R$-Cartier divisor and $\ep>0$.
Then the number of extremal contractions $\text{cont}_R$ and
corresponding rays $R$ such that $K_X+D$ is

\sn
(*) divisorially log terminal near
a generic point of the degenerate fibers of
$\text{cont}_R$,

\sn
and such that $(K_X+D+\ep H.R)<0$, is finite.

Thus the half-cone $(K_X+D.R)<0$ of the Kleiman-Mori
cone $\overline{\text\rm NE}(X/S)$ is locally polyhedral
when (*) holds for all extremal rays in it,
including the {\rm existence} of the extremal contractions
$\text{cont}_R$.
\endproclaim

In particular, if $X$ is normal projective with isolated singularities,
and $\Bbb Q$-Goren-\break stein, we have at most finite number of extremal
contractions which are negative for $K_X+\ep H$ (cf. [Ko2, Th.]).
Note that this does not mean that the half-cone $(K_X+D.R)<0$
is always locally polyhedral.
This means only, that whenever this fails, there exists
an extremal {\sl non-contracted\/} ray $R$ with $(K_X+D.R)<0$.

\sn
{\bf Proof.} See [K].\bsq

\medskip
\proclaim {Conjecture} In the Theorem we may replace
the divisorially log terminal property near by
the log canonical one in and always find a
family $\{C_{\lambda}\}/S$ of curves.
{\rm Note that, even in the definite case either we
consider a covering of the generic points of $E$, or
we admit effective 1-cycles as curves $C_{\lambda}$
for {\sl degenerations}.
Nevertheless in this case we abuse terminology and
say on a covering by curves.

We should have even more.
\endproclaim

\medskip
\proclaim {Theorem$'$} The Conjecture implies Corollaries 1-4
with the divisorially log terminal property near
replaced by the log canonical one in
and effective 1-cycles by curves.

If $\text{\rm dim }X\le 3$ the Conjecture and the improved
corollaries holds.
\endproclaim

\sn
{\bf Heuristic Arguments.} Here we deduce the Conjecture and,
in particular, the Theorem using
the LMMP (the Log Minimal Model Program).
Thus this proves the Theorem and Theorem$'$ for $\text{ dim }X\le 3$ [Sh3].
For all dimensions, a proof of the Theorem refines [K] and [MM],
and will be given later.

First, we discuss what means the log canonicity (or
terminality) in a generic point of a subvariety
$E\subseteq X$.
A generic point $P$ is one of its irreducible components.
The log divisor $K_X+D$ is log canonical (resp. log
terminal, purely log terminal, or Kawamata log terminal)
in $P$ if $K_X+D$ possesses that property in
a generic point of $P$ in a naive sense.
A more rigorous point of view means that
for the log discrepancies $a_i$ of the divisors $E_i$
such that the center $c(E_i)=P$, we have
the log canonical property $a_i\ge 0$
(resp. in addition, $a_i>0$ for such exceptional divisors
$E_i$ of one log nonsingular resolution over
a neighborhood of the generic point $P$;
$a_i>0$ for all such exceptional divisors
$E_i$; or $a_i>0$ for all such divisors $E_i$) [Sh2].
Note that the purely log terminal property
coincides with Kawamata one if $P$ has codimension
$\ge 2$.

A property near means in a neighborhood.

We may assume that $E$ and $S$ are irreducible.
For this add also that $\text{\rm dim } E/S=\
\text{\rm dim } E-\text{\rm\ dim } f(E)$
if $E$ is irreducible, and the maximum of
such dimensions for irreducible components
of $E$ in general.

Now we reduce the Theorem to the case
when $S$ is a point, i.e., $X$ is
a projective variety with a trivial
morphism $f\colon X \to \text{\rm\ pt. }$.
Taking a generic hyperplane section $H$ on $S$
and replacing $X$ by its inverse image $f^{-1}H$,
we could reduce the problem to the case when
$f(E)=p$ is a point.
Note that then $n=\text{\rm\ dim } E$.

If $S\not= p$ we can add a non-negative multiple
of $f^*H$, where $H$ is a generic hyperplane section
on $S$ through $f(E)=p$, in such a way that $K_X+D$
will be maximally log canonical in the generic point of $E$.
This means that $K_X+D$ is log canonical
in the generic point of $E$, and for any $\ep>0$,
$K_X+D+\ep f^*H$ is not so.
In other words, there exists a divisor $F$ with
the log discrepancy $0$ and with the center $E$.

In addition, if $S\not= p$, we may move and split $D$ into a boundary
preserving the log canonical property along $E$.
So, $K_X+D$ will be log canonical in a neighborhood
of the generic point of $E$.

If $F=E$ is a non-exceptional divisor on $X$,
then $D$ has the multiplicity 1 in $E$ and
we use the Adjunction Formula and Effectiveness
[Sh2, 3.1 and 3.2.2].
This reduces to the case when $X=E^{\nu}$ and $S=p$.

Using the LMMP we can do the same in general.
For this, we should replace $X$ by a strict log minimal
model $g\colon Y\to X$ of a neighborhood of $E$ with
respect to the log divisor $K_X+B$ where a boundary
$B$ has multiplicities $b_i=\text{\rm\ min }\{1,d_i\}$
(cf. [Sh2, 3.4]).
According to the above perturbation of $D$, $D=B$ in
a neighborhood of the generic point of $E$.
Perhaps after an additional monoidal transform, we
can also assume that $F=E$ is not exceptional on
$X\colon=Y$.
Thereafter we replace $K_X+D$ by  $f^*(K_X+D)=K_Y+D'$,
i.e., $D$ is replaced by an divisor $D'$.
$D'$ is effective by a construction of log models
and by the Negativity of Birational Contractions
[Sh2, 1.1].
The multiplicity of $E$ in $D$ is one by the construction.
Here we meet one annoyance:
$\text{\rm\ dim } E/S$ may be higher than
$n=\text{\rm\ dim } g(E)/S$ when $E$ is
exceptional for a model $g$.
However in that case we have an additional
structure, namely, a projective {\sl Iitaka}
morphism $i=g|_E\colon E\to g(E)$.
This means that $i$ is a projective contraction,
and $K_E+D_E=(K_X+D)|E$ is numerically trivial
for $i$.
Note that
$n=\text{\rm\ dim } E/S = \text{\rm\ dim } i(E)$.
The divisor $E$ itself is a projective variety
with a trivial morphism
$f\colon E \to \text{\rm\ pt. }$.
According to the Adjunction and Projection  formula,
the above properties of $i$ follows from
the assumption that $K_X+D$ is maximally log
canonical in $E$, i.e., $D$ has the multiplicity
$1$ in $E$.
These imply also that a required family of
curves in $E$ induces that of on a subvariety $g(E)$
of the original $X$.

Take now $X=E$. By the construction $X$ is again
a normal, projective variety with a trivial morphism
$f\colon X \to \text{\rm\ pt. }$.
In addition, we have an Iitaka morphism
$i\colon X\to E$ with respect to $K_X+D$.
Also by the construction and the Adjunction formula
$K_X+D$ is log canonical over the generic point
of $E$.
After additional blow ups, we can assume that
$X$ is $\Bbb Q$-factorial and $K_X+D$ is log terminal
over the generic point of $E$.
Since $D$ has multiplicities $d_i>1$ only for
prime divisors $D_i$ with $i(D_i)\not=E$, we
can drop these components and suppose that
$K_X+D$ is strictly log terminal, in particular,
$D$ is a boundary on $X$.
However, this can spoil the Iitaka morphism $i$,
namely, it may appear curves $C$ in fibers of $i$
such that $(K_X+D.C)\not=0$.
Nevertheless, $i$ will be again the Iitaka morphism
over the generic point of $E$.
This is enough and coherent to the following
arguments.

So, we should cover $X$ by a family of curves
$\{C_{\lambda}\}$ with $(-K_X-D.C_{\lambda})\le 2n$ where
$n=\text{\rm\ dim } E$.
Since $K_X+D$ is strictly log terminal and $X$ is projective, we
can apply the LMMP to $X$ with the log divisor $K_X+D$.
If $K_X+D$ is nef, then we take any family of curves
$\{C_{\lambda}\}$ covering $X$.
Otherwise we have an extremal contraction $g\colon X\to Z$.
If $g$ has the fiber type, then the generic fiber $F$ of
$g$ intersects the generic fiber of $i$ in a finite set.
Therefore, $\text{\rm dim } F\le n=\text{\rm\ dim } E$
and we can reduce the problem to fibers in this case.
Note also that divisorial contractions and flips
only decrease the intersection $(K_X+D.C_{\lambda})$
for the generic curve of a covering family of $X$.
Thus we can consider later only birational contractions.

If $g$ is birational, we make a divisorial contraction
or a flip of $X$.
However, the exceptional locus $F$ of $g$ intersects
the generic fiber of $i$ at most in a finite set.
If it does not intersect, we can make such transform
and preserve birationally $i$.
Otherwise $i(F)=E$ and $F$ is finite over the generic
point of $E$.
In this case, we use above arguments for $E=F$ and
the induction on the dimension of $X$, since the
above arguments restrict the proof to a divisor.
The correspondent Iitaka morphism is induced by $i$.

According to the termination of the LMMP,
this completes the reduction to the case
when $X=E$ and $S=p$ is a point.
Thus $X$ is normal and projective.
After an additional blow ups, we can assume that
$X$ is also $\Bbb Q$-factorial.
Now $n=dim X$, and  we should cover $X$ by a family
of curves $\{C_{\lambda}\}$ with
$(-K_X-D.C_{\lambda})\le 2n$.
Dropping $D$, and after an additional blow ups,
we can assume that $D=0$ and $X$ has only log
terminal or even terminal singularities.
In the last case we use the MMP.

If $K_X$ is nef, then we take any family of curves
$\{C_{\lambda}\}$ covering $X$.
Otherwise we apply the MMP.
According to the termination, after a finite number
of extremal transformation we reduce to the
case when $X$ possesses a Fano fibering.
Taking the generic fiber, we get the case when
$X$ is a Fano variety having only terminal
singularities.
Then the required covering family exists
according to Miyaoka and Mori [MM].

If $K_X+D$ is Kawamata log terminal in the generic
points of $E$ and $X\not=E$, we could replace
$f^*H$  by an effective Cartier divisor which
is numerically negative on $E$, for instance,
we could take a generic anti-hyperplane section
of $E$ over a neighborhood of $p=f(E)$.
Then subsequent reductions give
a stronger inequality $<2n$.
\bsq

Below in the proof of the Theorem, [MM]
plays the same role.
The previous heuristic arguments can be
replaced by the following refinement of
[K, Lemma].
We lose some properties stated in the conjecture.
This is the price of homological methods (cf. [Sh3]).

\medskip
\proclaim {Lemma} Let $f\colon X\to Y$ be a projective
morphism of normal algebraic varieties over a field of
characteristic zero,
$H$ an $\Bbb R$-Cartier divisor on $X$,
$D$ an effective $\Bbb R$-divisor,
$E$ an irreducible component of the degenerate locus
$\Cal Exc(f)$ of $f$, $n=\text{\rm\ dim } E$, and
$\nu \colon E^{\nu}\to E$ the normalization.
Suppose that $f$ is finite over the generic point of
$Y$, $H$ is nef with respect to $f$, $K_X+D$ is
$\Bbb R$-Cartier and divisorially log terminal near the generic
point of $E$  and $f(E)$ is a point.
Then
$$(H^{n-1}.K_X+D.E)\ge ((\nu^*H)^{n-1}.K_{E^{\nu}})$$
(and even $>$ if $H$ is ample and $K_X+D$ is Kawamata
log terminal in the generic point of $E$).
\endproclaim

\sn
{\bf Proof.} For the non-strong inequality we
can use an approximation of $H$
by $\Bbb Q$-Cartier ample divisors.
Moreover, we prove the following
polylinear inequality
$$(H_1.\ ...\ .H_{n-1}.K_X+D.E)\ge
(\nu^*H_1.\ ...\ .\nu^*H_{n-1}.K_{E^{\nu}})$$
where all $H_i$'s are nef (and even $>$ if all $H_i$'s
are ample and $K_X+D$ is Kawamata
log terminal in the generic points of $E$.
Since any ample $\Bbb R$-divisor is a sum
of an ample $\Bbb R$-divisor and an ample
$\Bbb Q$-divisor, we may assume that all
ample $H_i$'s are $\Bbb Q$-Cartier and even Cartier.)

Thus as in [K], we suppose that
$Y$ is affine, $H_i$ are very ample, and $n=1$,
i.e., we have no $H_i$'s and $E$ is a curve.

Since $K_X+D$ is divisorially log terminal near
the generic point of
$E$, then $\text{\rm Supp }D$ will be log nonsingular
near the generic point of $E$
unless $K_X+D$ is purely log terminal
near the generic point of $E$.

In the log nonsingular case $D$ is $\Bbb R$-divisor
in the generic point of $E$.
Then adding some small effective Cartier divisor
trivial near the generic point of $E$
we preserve the statement and may assume that
$D-\ep F\ge 0$ where $\ep>0$ and $F$ is
an effective Cartier divisor that coincides
with $D$ near the generic point of $E$.
This makes $K_X+D$ Kawamata log terminal
near the generic point of $E$ and in
particular purely log terminal near
the generic point of $E$.
As a limit for $\ep \to 0$ we get the required inequality.

Thus perturbing $D$ we may assume also that $K_X+D$ is
purely and even Kawamata log terminal near the generic point of $E$
except the case when $E$ is a divisor, i.e. $X$ is
a surface.
The above arguments also work for the last exception.

Since we are working over a field of characteristic
zero, we can replace it by $\Bbb C$ and replace
$f$ by an analytic contraction over a small analytic
neighborhood of $p=f(E)$ in $Y$.

Suppose that $(K_X+D.E)\le \text{\rm\ deg }K_{E^{\nu}}$.
Then there exists a Cartier divisor $A_0$ on $E$
such that $\text{\rm deg }A_0\ge (K_X+D.E)$ and
$H^0(E^{\nu},K_{E^{\nu}}-\nu^*A_0)\not=0$.
As in [K], we have $H^0(E,\omega(-A_0))\not=0$,
and we can extend $A_0$ to a Cartier divisor $A$ on
$X$.
Moreover, we may assume that $A-K_X-D$ is nef$/Y$
and $A$ is enough ample$/Y$ on components of $\Cal Exc(f)$
except $E$.
Enough means enough for vanishings below.

According to our assumptions
$K_X+D$ has nonpositive log discrepancies only for
divisors $F$ the centers of which intersects $Y$
in finite sets.

Since $A-K_X-D$ is nef, we
have $R^1f_*\Cal I_X(A)=0$ where $\Cal I_X$ is an ideal
sheaf of a subscheme $S\subset X$ which intersects $E$
in a finite set.
This is by the proof
of [KMM, 1.2.5] and the Kawamata log terminality of
$K_X+D$ near the generic point of $E$.
It is also important that $D$ is effective,and
$X$ is normal.
So, any function regular in codimension $2$ on $X$
will be regular everywhere.

According to the construction $A$ is enough ample
on $S/Y$, i.e., for $i\ge 1$,
$R^if_*\Cal O_S(A)=0$ by Serre.
Therefore the exact sequence
$$0\to \Cal I_X(A)\to \Cal O_X(A)\to \Cal O_S(A)\to 0$$
implies the vanishing $R^1f_*\Cal O_X(A)=0$.

Now the base change gives that
$\text{H}^1(f^{-1}p,A|_{f^{-1}p})=0$.
Again the restriction sequence
$$0\to \Cal J(A|_{f^{-1}p})\to
\Cal O_{f^{-1}p}(A|_{f^{-1}p})\to \Cal O_E(A_0)\to 0$$
and Serre vanishing imply the vanishing
$\text{H}^1(E,A_0)=0$
and we obtain the required inequality as in [K].
For this note that the support of the ideal sheaf $\Cal J$
intersects $E$ at most in one-dimensional set $E$ and we may choose $A$
such that it is
enough ample$/Y$, i.e., $\text{H}^2(f^{-1}p,\Cal J(A|_{f^{-1}p}))=0$.
\bsq

\sn
{\bf Proof of the Theorem.}
If $M=-\nu^*(K_X+D)$ is not $f$-nef on $E^{\nu}$ we
have a curve $C/S$ on $E^{\nu}$ with $(M.C)<0$.
Moreover, for 1-cycle $N\nu(C)/S$ with $N\gg 0$,
$(-K_X-D.N\nu(C))=(M.NC)\ll 0$.
Then it is easy to add to $N\nu(C)$ an effective
family of curves$/S$ covering $E$ and such
that for the obtained 1-cycles $C_{\lambda}/S$,
$(-K_X-D.C_{\lambda})<0\le 2n$.

Therefore, we suppose later that $M=-\nu^*(K_X+D)$
is $f$-nef on $E$.
Then we can use the arguments of [K, Proof of the
Theorem 1].
The only difference that now $M$  is nef.
Note that we can use [MM] even for a nef divisor $M$,
because the required inequalities are not strong.
Indeed, if $M$ is a limit of ample
$\Bbb Q$-Cartier divisors $M_n$ for $n\to \infty$.
Then the results of [MM] for $M_n$ gives in a limit
the same for $M$.\bsq

\medskip
\proclaim {Corollary 3 (cf. [K, Th. 2])} For $f\colon X\to S$
and $D$ as in the Theorem, let $E$ be a subvariety,
consisting of components of
$$\eqalign{\Cal Exc(f)\colon =& \{x\in X\mid \text {\rm an irreducible
component of a fiber of $f$ through $x$}\cr &\text{\rm  having
the dimension greater than $d$\ = dim }X/S \}\cr}.$$
Then $E$ is covered by a family (possibly disconnected) of effective
1-cycles
$\{C_{\lambda}\}/S$ with $(-K_X-D.C_{\lambda})\le 2(n-d)$ where
$n=\text{\rm\ dim } E/S$ (and even $<2(n-d)$ if
$K_X+D$ is Kawamata log terminal in the generic points of $E$).
Moreover, we could assume that the generic 1-cycles $C_{\lambda}$
are curves when $K_X+D$ is numerically definite, and these curves
$C_{\lambda}$ with $(K_X+D.C_{\lambda})< 0$
(resp. $\le0$) are rational.
\endproclaim

\sn
{\bf Proof} Take general hyperplane sections of $S$ (cf. the
Heuristic Arguments) and then use the Theorem.\bsq

\medskip
\proclaim {Corollary 4} If, in addition to the Corollary 3,
$F$ is
an $\Bbb R$-Cartier divisor such that $K_X+D+F$
is also divisorially log terminal near the generic points of $E$,
and $K_X+D$ is numerically semi-negative with respect
to $f$, then $E$ is covered by a family (possibly disconnected)
of effective 1-cycles (curves when $K_X+D+F$ numerically definite)
$\{C_{\lambda}\}/S$ with
$$(F.C_{\lambda})\ge -2(n-d)$$ (resp. $>-2(n-d)$ if $K_X+D+F$ is
Kawamata log terminal in the generic points of $E$).
\endproclaim

\sn
{\bf Remarks.}
(1) Conjecturally, we my suppose that the generic members of
$\{C_{\lambda}\}/S$ in the Corollaries 3 and  4 should
always be curves.
This is true when $n-d\le 3$.

(2) If $f$ is a projective morphism as in the Theorem we can
introduce its $D$-length as the maximum of $(-K_X-D.C_{\lambda})$
for families $\{C_{\lambda}\}/S$ covering $X$.
Of course, the length is 0  when $f$ is finite over the generic
point of $S$, and by the Theorem the length is not higher than
$2n$, if $K+D$ has in generic mild singularities,
where $n=\text{\rm\ dim } X/S$.

According to [MM], a smooth projective variety $X$ is
uniruled if and only if its length (=0-length) is positive.

Here we could state problems for the length similar to that of
for Fano indexes and minimal discrepancies [Sh1].

(3) The estimate $\le 2n$ (respectively $<2n$) is essentially
derived from that of [MM] for Fano $n$-folds with
terminal and even $\Bbb Q$-factorial singularities.
For $n\le 3$ and in the nonsingular case,
according to the classification of such Fanos,
we can improve the estimate up to $\le n+1$ (cf. [M]).
Since terminal singularities are nonsingular for $n\le 2$
this gives the same improvement in general for $n\le 2$.

So, we have the following problems here.
Does the estimate $\le 2n $ sharp for $n\ge 3$?
If not try to find the right one.
Why not $n+1$?

(4) Perhaps, the proven results and the Conjecture hold
in any characteristic.

\bigskip
{\bf Acknowledgment\/}
\medskip

Partial financial support was provided by the NSF under
the grant number DMS-9200933.

\bigskip
\noindent
{\bf Bibliography}

\medskip
\frenchspacing
\item {[K]} Kawamata Y., {\sl On the length of an extremal rational curve},
Invent. Math {\bf 105} (1991), no. 3, 609-611.
\item {[KMM]} Kawamata Y., Matsuda K. and Matsuki, {\sl
Introduction to the minimal model problem\/}, in Algebraic
Geometry Sendai 1985, Adv. St. Pure Math. {\bf 10} (1987),
283-360.
\item {[Ko1]} Koll\'ar J., {\sl The Cone theorem: Note to Kawamata's
``The cone of curves of algebraic varieties"\/}, Ann. of Math.
{\bf 120} (1984), 1-5.
\item {[Ko2]} Koll\'ar J., {\sl Cone theorems and cyclic covers\/},
preprint.
\item {[M]} Mori S., {\sl Threefolds whose canonical bundles are
not numerically effective\/}, Ann. of Math. {\bf 116} (1982),
133-176.
\item {[MM]} Miyaoka Y. and Mori S., {\sl A numerical criterion of
uniruledness}, Ann. of Math. {\bf 124} (1986), 65-69.
\item {[Sh1]} Shokurov V.V., {\sl Problems about Fano varieties\/},
In Birational Geometry of Algebraic Varieties, Open Problems,
The XXIII International Symposium, Division of Mathematics,
The Taniguchi Foundation, August 22 -- August 27, 1988,
Katata, 30-32.
\item {[Sh2]} Shokurov V.V., {\sl 3-fold log flips\/}, Russian Acad.
Sci. Izv. Math., Vol. {\bf 40} (1993), No. 1, 95-202.
\item {[Sh3]} Shokurov V.V., {\sl 3-fold log models\/}, preprint.

\end